\documentclass[preprint,12pt]{elsarticle}
\usepackage[hmargin=0.8in,vmargin=0.9in]{geometry}
\usepackage{fancyhdr}
\usepackage{graphicx}
\usepackage{subfigure}
\usepackage{amssymb}
\usepackage{amsmath}
\usepackage{amsfonts}
\usepackage{theorem}
\usepackage{mathrsfs}
\usepackage{mathtools}
\usepackage{makecell}
\usepackage{multirow}
\usepackage{empheq}
\usepackage{algorithm}
\usepackage{algpseudocode}
\usepackage{bm}
\usepackage{color}
\algnewcommand\algorithmicforeach{\textbf{for each}}
\algdef{S}[FOR]{ForEach}[1]{\algorithmicforeach\ #1\ \algorithmicdo}
\usepackage{algpseudocode}
\usepackage{varwidth}
\usepackage{setspace}
\usepackage[labelfont=bf]{caption}
\usepackage{epstopdf}
\usepackage{float}
\usepackage{cite}

\makeatletter
\def\ps@pprintTitle{%
	\let\@oddhead\@empty
	\let\@evenhead\@empty
	\let\@oddfoot\@empty
	\let\@evenfoot\@oddfoot
}
\makeatother

\usepackage{booktabs} 
\usepackage{array} 
\usepackage{paralist} 
\usepackage{verbatim} 
\usepackage{subfigure} 
\usepackage{natbib}
\usepackage{hyperref}
\usepackage{xstring}

\usepackage{hyperref}
\hypersetup{
	colorlinks   = true,
	citecolor    = blue,
	linkcolor=blue  
}

\usepackage{xstring}
\makeatletter

{\theorembodyfont{\rmfamily}\newtheorem{remark}{Remark}[section]}

\allowdisplaybreaks
\numberwithin{equation}{section}
\numberwithin{figure}{section}
\numberwithin{table}{section}

\usepackage{ae,aecompl}

\def\Chi {{\Large{\mbox{$\chi$}}}}

\newcommand{\mF}{\bm{F}}
\newcommand{\mS}{\bm{S}}

\newcommand{\mG}{\bm{G}}

\newcommand{\mH}{\bm{H}}

\newcommand{\mU}{\bm{U}}

\newcommand{\Ups}{\bm{\Upsilon}}

\newcommand\eref[1]{(\ref{#1})}%
\newcommand\tref[1]{Table \ref{#1}}%
\newcommand\fref[1]{Fig.~\ref{#1}}%
\makeatletter
\newcommand{\refcheckize}[1]{%
	\expandafter\let\csname @@\string#1\endcsname#1%
	\expandafter\DeclareRobustCommand\csname relax\string#1\endcsname[1]{%
		\csname @@\string#1\endcsname{##1}\wrtusdrf{##1}}%
	\expandafter\let\expandafter#1\csname relax\string#1\endcsname
}
\makeatother

\newcommand*\xbarsh[1]{%
	\hbox{%
		\vbox{%
			\hrule height 0.5pt 
			\kern0.4ex
			\hbox{%
				\kern-0.22em
				\ensuremath{#1}%
				\kern-0.05em
			}%
		}%
	}%
}

\newcommand*\xbar[1]{%
	\hbox{%
		\vbox{%
			\hrule height 0.5pt 
			\kern0.4ex
			\hbox{%
				\kern-0.05em
				\ensuremath{#1}%
				\kern-0.05em
			}%
		}%
	}%
}

\begin{document}
	\begin{frontmatter}
		
		\author{Thuong Nguyen\fnref{label3}}
		\fntext[label3]{Department of Mathematics, The University of Utah, Salt Lake City, UT 84112, USA;
			{\tt tnguyen@math.utah.edu} }
		
		\title{Adaptive Central-Upwind Scheme on Triangular Grids for the Shallow Water Model with variable density}
		
		
		\begin{abstract}
			In this paper, we construct a robust adaptive central-upwind scheme on unstructured triangular
			grids for two-dimensional shallow water equations with variable density. The method is well-balanced,
			positivity-preserving, and oscillation free at the curve where two types of fluid merge. The proposed approach is an extension of the adaptive  well-balanced,
			positivity-preserving scheme developed in Epshteyn and Nguyen (arXiv preprint arXiv:2011.06143, 2020). In particular, to  preserve ``lake-at-rest'' steady states, we utilize the Riemann Solver with appropriately rotated coordinates to obtain the point values in neighborhood of the fluid interface.  In addition,  to improve the efficiency of an adaptive method in the multifluid flow, the curve of density discontinuity is reconstructed by using the level set method and volume fraction method.  To demonstrate the accuracy, high-resolution, and efficiency of the new adaptive central-upwind scheme, several challenging tests for Shallow water models with variable density are performed.
		\end{abstract}
		
	\end{frontmatter}

	\section{Introduction}\label{sec:introduction}
The main goal of this paper is to develop an adaptive well-balanced positivity-preserving
central-upwind scheme on triangular grids for shallow water equations with variable density
(SWEDs). The two-dimensional (2-D) system of SWEDs can be written as,
	
\begin{subequations}\label{eq:swed}
	\begin{align}
		&w_t+(hu)_x+(hv)_y=0,\label{eq:swed1}\\
		&(hu)_t+\Big(hu^2+\frac{g}{2\rho_0}h^2\rho\Big)_x+(huv)_y=-\frac{g}{\rho_0}h\rho B_x,\label{eq:swed2}\\
		&(hv)_t+(huv)_x+\Big(hv^2+\frac{g}{2\rho_0}h^2\rho\Big)_y=-\frac{g}{\rho_0}h\rho B_y,\label{eq:swed3}\\
		&(h\rho)_t+(hu\rho)_x+(hv\rho)_y=0,\label{eq:swed4}
	\end{align}
\end{subequations}	
	where $t$ is the time, $x$ and $y$ are spatial coordinates ($(x,y)\in\Omega$), $h(x,y,t)$ is the water height, $\rho(x,y,t)$ is the density, $u(x,y,t)$ and
	$v(x,y,t)$ are the $x$- and $y$-components of the flow velocity, $B(x,y)$ is the bottom topography, $g$ is the constant gravitational
	acceleration, and $\rho_0$ is the reference density.  The system \eref{eq:swed1}--\eref{eq:swed4} was  proposed in \citep{dellar2003common, ripa1993conservation, ripa1995improving, GHAZIZADEH2020104633} as a variation of the Saint-Venant equations to model multi-phase flows
	in estuaries or deep ocean currents. The derivation of the system is based on hydrostatic approximation which eliminates the variability in the z-direction. The
	design of robust and accurate numerical algorithms for computing the solutions of SWEDs system is an important and challenging problem that has been
	extensively studied in the recent years.
	
A number of numerical schemes for balance laws have been introduced in recent years, \citep{MR1763829,KNP,KLab,SDBL,KP1,KM,CEHK,MR3229988,MR3187922,KurLiu,MR3071176,KP05,AIPBEKP,MR2804645,MR3440159,LAEK}. Most of them utilize a Riemann problem solver for the upwind evolution of the calculated
solution. However, as discussed in \citep{chertock2014central}, the eigensystem of the system  \eref{eq:swed1}--\eref{eq:swed4}  may be incomplete due to the resonance phenomenon. Hence, it may be very difficult to design a reliable upwind scheme for the SWEDs. In our paper, we therefore use central-upwind schemes which are Riemann-problem-solver free methods, \citep{kurganov2007reduction, kurganov2001semidiscrete, LAEK}. Central-upwind schemes have been referred to ``black-box'' solvers for general multidimensional systems of hyperbolic systems of conservation laws. In our prior work \citep{epshteyn2020adaptive}, we have derived a successful adaptive central-upwind method for  Saint-Venant system on triangular grids.

Similar to the Saint-Venant system, a good method for  SWEDs system should preserve the non-negativity of $h$
and $\rho$, which is called the positivity-preserving property. In addition, the scheme must ensure a well-balanced property obtained when the numerical method preserve ``lake-at-rest'' steady-state solutions. Otherwise, the numerical method may lead to significant oscillations. Note that, the system \eref{eq:swed1}--\eref{eq:swed4} admits the following two ``lake-at-rest'' steady-state solutions, \citep{chertock2014central}:

\begin{equation}
	w=\max\big\{C,B(x,y)\big\},\quad C=\rm{Const}, \quad\rho =P\equiv Const, \quad u\equiv v\equiv0,  \label{eq:lara}
\end{equation}

and

\begin{equation}
	B\equiv \rm{Const}, \quad h^2\rho \equiv Const, \quad u\equiv v\equiv0.  \label{eq:larb}
\end{equation}
Preserving the solution \eref{eq:larb} is a big challenge for numerically solving the system \eref{eq:swed1}--\eref{eq:swed4} since using the conventional central-upwind methods may not ensure the variable $h^2\rho$, so-called variable pressure,  to be constant at the contact waves. Therefore, for our adaptive scheme of SWEDs, we will utilize the central-upwind scheme which is derived from \citep{chertock2014central} for shallow water model with horizontal temperature variable (the system in \citep{chertock2014central}  has similar properties with the SWEDs). In \citep{chertock2014central}, the proposed second-order semi-discrete central-upwind scheme is capable of preserving the “lake at rest” steady state \eref{eq:lara} and \eref{eq:larb}  as well as the positivity of the water depth and the temperature (the temperature variable is equivalent to the variable density in our work). In particular, to preserve the second type of “lake
at rest” steady state \eref{eq:larb} and suppress the pressure oscillations across the interface, an efficient interface
tracking method is performed. The main idea of the
interface approach in \citep{chertock2014central} is to completely avoid to use the information from the cells where two types of fluids are numerically mixed, so-called ``mixed'' cells when evolving the solution in the neighborhood of the interface. The data in the “mixed”
cells is replaced by the interpolated values that are calculated using the reliable information
from the nearby “single fluid” cells. Namely, the point values in ``mixed'' cells are obtained by using the approximated solution
of the 1-D Riemann problems between the reliable “single fluid” cell averages. However, the central-upwind method and the interface tracking in  \citep{chertock2014central} are designed on structured rectangular grids. In practice, one needs to deal with complicated geometries, where the use of triangular grids could be advantageous or even unavoidable. Hence, in this study, we extend the interface tracking method from \citep{chertock2014central} to unstructured triangular grids by utilizing the idea of the rotated coordinates proposed in \citep{balsara2014multidimensional}. The 1-D Riemann solver \citep{chertock2014central} is then performed in the direction of normal vectors on each edge of ``mixed'' cells. This extended central-upwind scheme maintains the well-balanced and positivity-preserving properties.

	In addition to achieving the accuracy of the solution, minimizing computational cost is also one of the major challenge in the modeling and numerical
	analysis of hydrodynamics. The traditional numerical schemes  for the system
	\eref{eq:swed1}--\eref{eq:swed4} are based on the use of very fine fixed meshes to reconstruct
	delicate features of the solution. However, this can lead to high computational cost, as well as poor resolution of all small scale features of the problem. In many engineering and scientific applications, it is beneficial to use adaptive meshes for improving the accuracy of the approximation at a much lower
	cost. Therefore, another goal of this work is to design an adaptive numerical algorithm for shallow water equations with variable density. 
 There	is some very recent effort on the design of adaptive well-balanced and
	positivity-preserving central-upwind schemes on quad-tree grids for
	shallow water models \citep{MR3315267,GHAZIZADEH2020104633, ghazizadeh2020adaptive}, but no research has been
	done for the development of such adaptive schemes for the shallow water with variable density on unstructured triangular grids. 	
	
	As a part of the mesh reconstruction, we will need to project the data from the old mesh onto a new adaptive mesh. Since we avoid using the cell averages of mixed cells, the data of a new cell which is reconstructed from an old mixed cell is calculated based on the location of the density jumps and the data from the nearby ``reliable'' cells. Therefore, the interface separating different fluid phases must be accurately tracked or captured. In each ``mixed'' cell, we will reconstruct
	an approximate interface, which is called ``interface reconstruction''. Many methods have been developed for this purpose in the last two decades, \citep{osher2001level, sethian2003level, hirt1981volume, rider1998reconstructing, unverdi1992front, tryggvason2001front}. In our work, we consider the interface reconstruction method derived in \citep{yang2006adaptive} due to the great advantages of this approach such as high-resolution, efficiency and simplicity. 
	
	\par This paper is organized as follows. In Section \ref{sec:centralupwind}, we present a well-balanced positivity-preserving central-upwind scheme
	on unstructured triangular grids for SWEDs which serves as the
	underlying discretization for the developed adaptive algorithm. Next, in Section \ref{sec:mixcelltrack}, we present the procedure to detect mixed cells. The interface approximation will be briefly reviewed in Section \ref{sec:interface}. In Section \ref{sec:correction}, we discuss the cell averages correction, which is used to prevent the density diffusion around the density jumps. We
	summarize the adaptive central-upwind method in Section \ref{sect3a0}. We
	discuss a strategy of adaptive mesh refinement in
	Section \ref{sect3a}. In Section \ref{sec:adaptivetime},  we present an adaptive
	second-order strong stability preserving Runge-Kutta method, employed
	as a part of time evolution for the adaptive central-upwind scheme. We
	develop a local posteriori error estimator in Section \ref{sect3c} which is
	used as a robust indicator for the adaptive mesh refinement in our
	work. Finally, in Section \ref{sec:numerical}, we illustrate the high accuracy and
	efficiency of the developed adaptive central-upwind scheme on a number
	of challenging tests for multi-fluid shallow water models.

	\section{The Central Upwind Scheme for SWEDs in Triangular Mesh}\label{sec:centralupwind}
 
 In this section, we focus on developing the central-upwind method for the SWEDs system \eref{eq:swed1}-\eref{eq:swed4},  \citep{GHAZIZADEH2020104633}. In particular, we will extend the adaptive scheme from the Saint-Venant system in \citep{epshteyn2020adaptive} to the SWEDs system. The developed scheme will eliminate the oscillation appearing at the interface and ensure the well-balanced and positivity-preserving properties. 

		
	We first rewrite the SWEDs system \eref{eq:swed1}-\eref{eq:swed4} into the vector form as,
	
	\begin{equation}
		\mathbf{U}_t+	\mathbf{F}(\mathbf{U},B)_x+	\mathbf{G}(\mathbf{U},B)_y=\mathbf{S}(\mathbf{q},B),\label{eq8}
	\end{equation}
	where 
		$$U=(w,hu,hv,h\rho),$$
	and the fluxes and source term are: 
	
	\begin{equation}
		\begin{aligned}
		\mathbf{F}(\mathbf{U},B)&=(hu,\dfrac{(hu)^2}{w-B}+\dfrac{g}{2\rho_0}\rho(w-B)^2,\dfrac{(hu)(hv)}{w-B},hu\rho)^T,\\
			\mathbf{G}(\mathbf{q},B)&=(hv,\dfrac{(hu)(hv)}{w-B},\dfrac{(hv)^2}{w-B}+\dfrac{g}{2\rho_0}\rho(w-B)^2,hv\rho)^T,\\\mathbf{S}(\mathbf{U},B)&=(0,-g\frac{\rho}{\rho_0}(w-B) B_x ,-g\frac{\rho}{\rho_0}(w-B) B_y,0)^T.
		\end{aligned}
	\end{equation}
	
\begin{figure}[ht!]
		\centering
		\includegraphics[width=8cm]{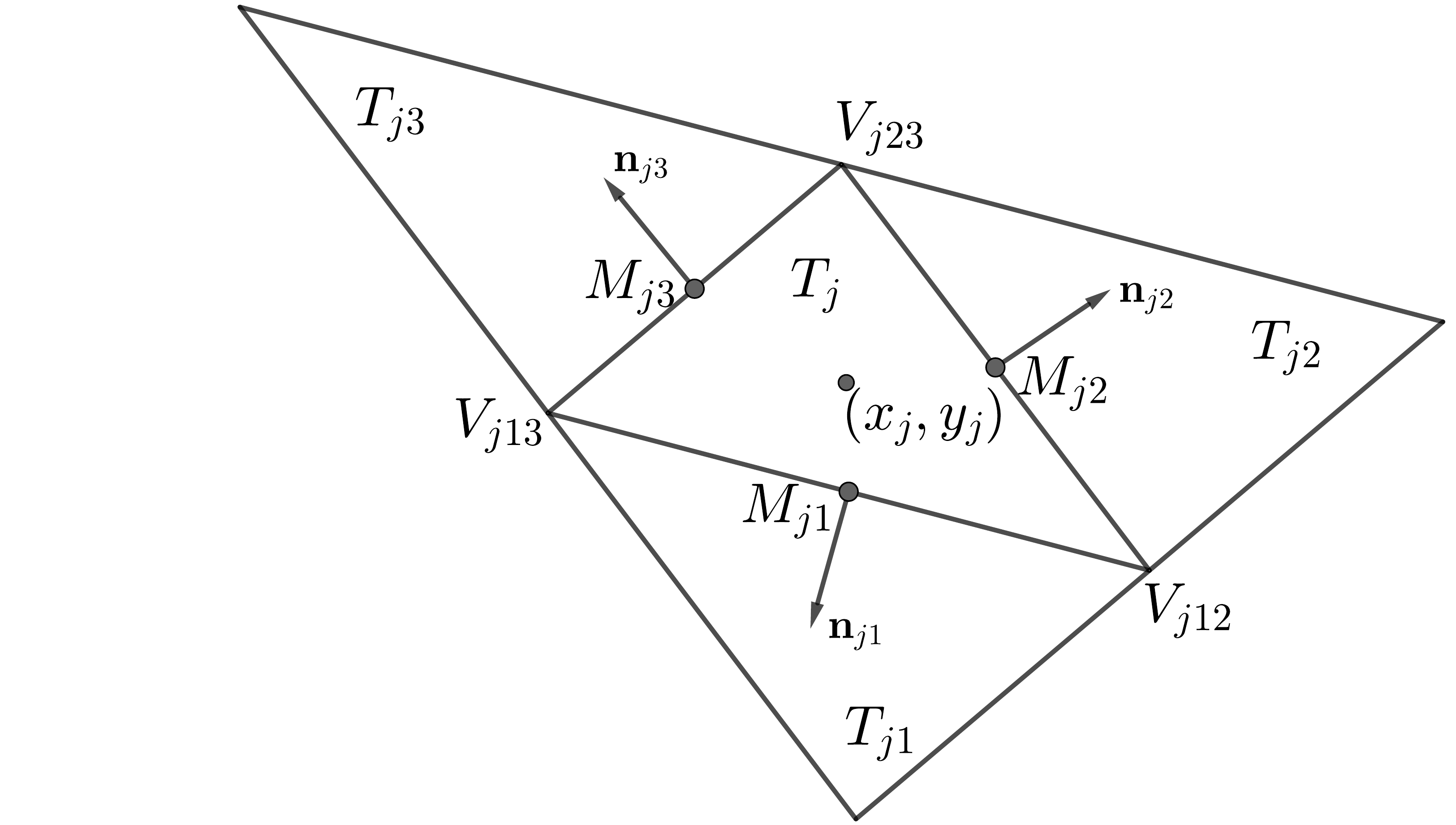}\\
		\vspace*{5mm}
		\caption{{ A typical triangular cell with three neighbors.\label{fig:trinei}}}
		\vspace*{2mm}
	\end{figure}

	In our research, we consider the triangular mesh illustrated in \fref{fig:trinei} with the following notations.
	
	\noindent
	${\mathcal {T}}:=\bigcup_j T_j$ is an unstructured triangulation of the computational domain $\Omega$;
	
	\noindent
	$T_j\in \mathcal{T}$ is a triangular cell of size $|T_j|$ with the barycenter $(x_j,y_j)$;
	
	\noindent
	$V_{j\kappa}=(\widetilde{x}_{j\kappa},\widetilde{y}_{j\kappa}),~\kappa=12,23,31$ are the three vertices of $T_j$;
	
	\noindent
	$T_{jk},~k=1,2,3$ are the neighboring triangles that share a common side with $T_j$;
	
	\noindent
	$\ell_{jk}$ is the length of the common side of $T_j$ and
	$T_{jk}$,  and $M_{jk}$ is its midpoint;
	
	\noindent
	$\bm{n}_{jk}:=(\cos(\theta_{jk}),\sin(\theta_{jk}))^\top$ is the outer unit normal to the $k$th side of $T_j$.

	Next, the bottom topography $B$ is replaced with its continuous piecewise linear approximation $\widetilde{B}$ given by
\begin{equation*}
	\left|\begin{array}{ccc}x-\widetilde{x}_{j12}&y-\widetilde{y}_{j12}&{\widetilde B}(x,y)-\widehat{B}_{j12}\\[0.5ex]
		\widetilde{x}_{j23}-\widetilde{x}_{j12}&\widetilde{y}_{j23}-\widetilde{y}_{j12}&\widehat{B}_{j23}-\widehat{B}_{j12}\\[0.5ex]
		\widetilde{x}_{j13}-\widetilde{x}_{j12}&\widetilde{y}_{j13}-\widetilde{y}_{j12}&\widehat{B}_{j13}-\widehat{B}_{j12}\end{array}\right|=0,
	\quad(x,y)\in T_j,
\end{equation*}
where, in the case of continuous bottom topography,
$\widehat{B}_{j\kappa}:=B(V_{j\kappa})~\kappa=12,23,31$. Then, denote:
$$
B_{jk}:=\widetilde{B}(M_{jk}),\quad B_j:=\widetilde{B}(x_j,y_j)=\frac{1}{3}(\widehat{B}_{j12}+\widehat{B}_{j23}+\widehat{B}_{j13}).
$$

At time $t$, define by $\xbar\mU_j(t)$ the approximation of the cell averages of the solution,
\begin{equation*}
	\xbar\mU_j(t)\approx\frac{1}{|T_j|}\iint\limits_{T_j}\mU(x,y,t)\,dxdy.
\end{equation*}

From now on, in order to shorten the formula,  we suppress the time-dependence in the algorithm. All
indexed quantities used in the following formula will be computed at time $t$. Then, as shown in \citep{LAEK,MR2804645}, the semi-discrete
second-order central-upwind scheme for the Saint-Venant system \eref{eq:swed1}-\eref{eq:swed4}  on triangular grids is written as the following system of ODEs,
\begin{equation}
	\begin{aligned}
		\frac{d\,\xbar{\mU}_j}{dt}=&-\frac{1}{|T_j|}\big[\mH_{j1}+\mH_{j2}+\mH_{j3}\big]+\xbar{\mS}_j,
		\label{10eq8}
	\end{aligned}
\end{equation}
where the numerical fluxes through the edges of the triangular
cell $T_j$ are

	\begin{align}
		\mH_{jk}=&\frac{\ell_{jk}\cos(\theta_{jk})}{a_{jk}^{\rm in}+a_{jk}^{\rm out}}
		\Big[{a_{jk}^{\rm in}\mF(\mU_{jk}(M_{jk}),B_{jk})+a_{jk}^{\rm out}\mF(\mU_j(M_{jk}),B_{jk})}\Big]\\
		+&\frac{\ell_{jk}\sin(\theta_{jk})}{a_{jk}^{\rm in}+a_{jk}^{\rm out}}
		\Big[{a_{jk}^{\rm in}\mG(\mU_{jk}(M_{jk}),B_{jk})+a_{jk}^{\rm out}\mG(\mU_j(M_{jk}),B_{jk})}\Big]\\
		-&\ell_{jk}\frac{a_{jk}^{\rm in}a_{jk}^{\rm out}}{a_{jk}^{\rm in}+a_{jk}^{\rm out}}\big[\mU_{jk}(M_{jk})-\mU_j(M_{jk})\big],\quad k=1,2,3.
		\label{flux}
	\end{align}

Here, $\bm{U}_{j}(M_{jk})$ and $\bm{U}_{jk}(M_{jk})$ are the
reconstructed point values of $\mU$ at the middle points of
the edges $M_{jk}$. To obtain these values, first, 
a piecewise linear reconstruction of the variables
$\Ups:=(w,u,v,\rho)^\top$ is computed as,
\begin{equation}
	\widetilde\Ups(x,y)=\sum\limits_j\Ups_j(x,y)\Chi_{T_j},\quad\Ups_j(x,y):=\Ups_j+(\widehat{\Ups}_x)_j(x-x_j)+(\widehat{\Ups}_y)_j(y-y_j),
	\label{eq:pwlinear}
\end{equation}
where $\Chi_{T_j}$ is the characteristic function of the cell $T_j$, $\Ups_j$ are the point values of $\Ups$ at the cell centers, and
$(\widehat{\Ups}_x)_j$ and $(\widehat{\Ups}_y)_j$ are the limited partial derivatives (see Section 3 in \citep{LAEK} for more details of the computation). In order to compute the cell center point values of the density, $\rho_j\approx \rho(x_j, y_j, t)$ in cell $T_j$ and velocities $u_j\approx u(x_j, y_j, t)$ and $v_j\approx v(x_j, y_j, t)$, we use the desingularization procedure presented in \citep{LAEK}. After that, the second and third components of the point values $\bm{U}_{j}(M_{jk})$ and $\bm{U}_{jk}(M_{jk})$ are obtained from,
$\Ups_j(M_{jk})$ and $\Ups_{jk}(M_{jk})$,
$$
\begin{aligned}
	&(hu)_j(M_{jk})=(w_j(M_{jk})-B_{jk})\,u_j(M_{jk}),&&(hu)_{jk}(M_{jk})=(w_{jk}(M_{jk})-B_{jk})\,u_{jk}(M_{jk}),\\
	&(hv)_j(M_{jk})=(w_j(M_{jk})-B_{jk})\,v_j(M_{jk}),&&(hv)_{jk}(M_{jk})=(w_{jk}(M_{jk})-B_{jk})\,v_{jk}(M_{jk}).\\
	&(h\rho)_j(M_{jk})=(w_j(M_{jk})-B_{jk})\,\rho_j(M_{jk}),&&(h\rho)_{jk}(M_{jk})=(w_{jk}(M_{jk})-B_{jk})\,\rho_{jk}(M_{jk}).
\end{aligned}
$$

However, in some numerical examples, the oscillation may appear when we use the piecewise linear approximation \eref{eq:pwlinear} to obtain the point values in some cells. It may occur due to the appearance of local extrema values at middle points or as a consequence of the density discontinuity. To prevent this oscillation, we will use some techniques that will be discussed later in Sections \ref{sec:scalelimiter} and \ref{sec:riemannsolver}. 

\par	In \eref{flux}, $a_{jk}^{\rm in}$ and $a_{jk}^{\rm out}$ are the one-sided local speeds of propagation in the directions $\pm\bm{n}_{jk}$.
These speeds are related to the largest and smallest eigenvalues of the Jacobian matrix
$J_{jk}=\cos(\theta_{jk})\,\frac{\partial {\bm F}}{\partial {\bm U}}+\sin(\theta_{jk})\,\frac{\partial {\bm G}}{\partial {\bm U}}$,
denoted by $\lambda_+[J_{jk}]$ and $\lambda_-[J_{jk}]$, respectively, and are defined by
\begin{equation}
	\begin{aligned}
		&a^{\rm in}_{jk}=-\min\{\lambda_-[J_{jk}(\mU_j(M_{jk}))],\,\lambda_-[J_{jk}(\mU_{jk}(M_{jk})],\,0\},\\
		&a^{\rm out}_{jk}=\max\{\lambda_+[J_{jk}(\mU_j(M_{jk}))],\,\lambda_+[J_{jk}(\mU_{jk}(M_{jk})],\,0\},
	\end{aligned}\label{eq:speeds}
\end{equation}
where
$$
\begin{aligned}
	&\lambda_\pm[J_{jk}(\mU_j(M_{jk}))]=\cos(\theta_{jk})u_j(M_{jk})+\sin(\theta_{jk})v_j(M_{jk})\pm\sqrt{\frac{g}{\rho_0}h_j(M_{jk})\rho_j(M_{jk})},\\
	&\lambda_\pm[J_{jk}(\mU_{jk}(M_{jk}))]=\cos(\theta_{jk})u_{jk}(M_{jk})+\sin(\theta_{jk})v_{jk}(M_{jk})\pm\sqrt{\frac{g}{\rho_0}h_{jk}(M_{jk})\rho_{jk}(M_{jk})}.
\end{aligned}
$$
\begin{remark}
	In order to avoid division by 0 (or by a very small positive number), the numerical flux \eref{flux} is replaced with
	$$
	\begin{aligned}
		\mH_{jk}=&\frac{\ell_{jk}\cos(\theta_{jk})}{2}\left[\mF(\mU_{jk}(M_{jk}),B_{jk})+\mF(\mU_j(M_{jk}),B_{jk})\right]\\
		+&\frac{\ell_{jk}\sin(\theta_{jk})}{2}\left[\mG(\mU_{jk}(M_{jk}),B_{jk})+\mG(\mU_j(M_{jk}),B_{jk})\right]
	\end{aligned}
	$$
	wherever $a_{jk}^{\rm in}+a_{jk}^{\rm out}<\sigma$. In
	all of the reported numerical examples in Section \ref{sec:numerical}, we have taken $\sigma=10^{-6}$.
\end{remark}	 
\par 	Finally, the cell average of the source term $\mS_j$ in \eref{10eq8},
\begin{equation*}
	\xbar\mS_j(t)\approx\frac{1}{|T_j|}\iint\limits_{T_j}\mS\big(\mU(x,y,t),B(x,y)\big)\,dxdy,
\end{equation*}
has to be discretized in a well-balanced manner presented later in Section \ref{sec:source}\\

\subsection{The Limiter for Piecewise Linear Approximation}\label{sec:scalelimiter}

As discussed above, we may observe the oscillation when using piecewise linear approximation \eref{eq:pwlinear} due to the appearance of local extrema at midpoints $M_{jk}$, see \citep{LAEK}. Hence, in order to maintain the numerical stability of the scheme, the following monotonicity condition must be satisfied.

\begin{equation}
	\min\big(\Ups^{(i)}_j,\Ups^{(i)}_{jk}\big)\le\Ups^{(i)}_j(M_{jk})\le\max\big(\Ups^{(i)}_j,\Ups^{(i)}_{jk}\big),\quad k=1,2,3,\label{eq:localextrema}
\end{equation}
where 	$\Ups^{(i)}_j(M_{jk})$ is the $i$-th component of the point values at midpoint $M_{jk}$ which is obtained by the linear reconstruction \eref{eq:pwlinear}. Here, $\Ups^{(i)}_j$ and $\Ups^{(i)}_{jk}$ are the center point value in cell $T_j$ and the neighboring cell $T_{jk}$, see \citep{kurganov_2018}.  To ensure that the reconstructed point value $\Ups^{(i)}_j(M_{jk})$ is between two cell averages $\Ups^{(i)}_j$ and $\Ups^{(i)}_{jk}$, in \citep{LAEK, AIPBEKP}, the gradients $\widehat{\Ups}_x$ and $\widehat{\Ups}_y$ are set to zero for cells violating at least one of the inequalities \eref{eq:localextrema}. However, using zero gradients may reduce the convergence rate in numerical simulation. Hence, in our research, instead of zero gradients, we consider the idea of the positivity-preserving limiter  developed in \citep{qin2016bound,xing2013positivity} to correct the piecewise linear reconstruction.
The details for this correction are presented as follows. 

Consider an arbitrary cell $T_j$ with the polynomial $\Ups^{(i)}_j(x,y)$ that does not satisfy at least one of local extrema conditions \eref{eq:localextrema}. In that case, we will replace $\Ups^{(i)}_j(x,y)$ by

\begin{equation}
	(\Ups^{(i)}_j)^*(x,y)=\theta(\Ups^{(i)}_j(x,y)-\Ups^{(i)}_j)+\Ups^{(i)}_j,\label{eq:eq1}
\end{equation}
where $\theta\in[0,1]$ such that $(\Ups^{(i)}_j)^*(x,y)$ satisfies the inequalities \eref{eq:localextrema} for all $k=1,2,3$. We then need to solve for the constant $\theta$. Plugging \eref{eq:eq1} into the inequalities \eref{eq:localextrema}, we have

\begin{equation}
	\theta\in\left[\min(\alpha_{jk},\beta_{jk}),\max(\alpha_{jk},\beta_{jk})\right], \quad  k=1,2,3\end{equation}
or 
\begin{equation}
	\theta\in[\max_k(\min(\alpha_{jk},\beta_{jk})),\min_k(\max(\alpha_{jk},\beta_{jk}))]
\end{equation}
where $\alpha_{jk}=\dfrac{\min\big(\Ups^{(i)}_j,\Ups^{(i)}_{jk}\big)-\Ups^{(i)}_j}{\Ups^{(i)}_j(M_{jk})-\Ups^{(i)}_j}, \quad \beta_{jk}=\dfrac{\max\big(\Ups^{(i)}_j,\Ups^{(i)}_{jk}\big)-\Ups^{(i)}_j}{\Ups^{(i)}_j(M_{jk})-\Ups^{(i)}_j}.$

Hence we can take \begin{equation}
	\theta=\min(\min_k(\max(\alpha_{jk},\beta_{jk})),1),\label{eq2}\end{equation}
Note that $\theta\geq 0$ since $\theta=0$, which corresponds to the constant approximation for $(\Ups^{(i)}_j)^*(x,y)$, is one solution of the inequalities \eref{eq:localextrema}. Hence, with this optimized limiter, the point values stay within the given range. In addition, the approximated values of water depth and density at middle points are also positive as $\Ups^{(i)}_j(M_{jk})\ge\min\big(\Ups^{(i)}_j,\Ups^{(i)}_{jk}\big)>0$. Therefore, the stability is ensured without using zero-gradients. 

\subsection{Riemann Solver for Mixed Cells}\label{sec:riemannsolver}

As mentioned in Section \ref{sec:introduction} and \citep{chertock2014central}, a good scheme should not develop
spurious pressure oscillations in the neighborhood of density jumps. The oscillations appear when we numerically solve the compressible multifluid systems  by using conventional Godunov-type finite-volume methods. We define the cells which contains the curves of density discontinuity as ``mixed'' cells. Note that the quantities in mixed cells are calculated as an artificial numerical mixture of two different fluids. Hence, the cell averages of mixed cells may have no or very little physical sense and become 'unreliable' \citep{chertock2014central}. Consequently, using the cell averages in mixed cells to obtain the linear approximation \eref{eq:pwlinear} in the neighborhood of the interface will lead to unexpected pressure oscillations. In \citep{chertock2014central}, to eliminate the oscillation, instead of the linear approximation \eref{eq:pwlinear}, the point values in the ``mixed'' cells are calculated by using the approximated solution of 1-D Riemann problems. Namely, according to this approach, the 1-D Riemann Solver is applied in the $x$-direction to approximate the point values on the vertical boundaries of each rectangular mixed cell. Similarly, for the point values on horizontal sides of the rectangular, 1-D Riemann problem is considered in $y$-direction. This leads to the idea of applying Riemann Solver approach in the direction of the normal vector, so called ``normal direction'', to get the middle point values on each side of the mixed cells in triangular meshes.  We can understand this idea as rotating the reference coordinate. Recently, a method using rotated coordinate frame to solve multi-dimensional Riemann problems has been proposed in \citep{balsara2014multidimensional}. In such approach, the 2-D Riemann Problem is converted to 1-D problem in a particular direction in order to easily obtain the corresponding solution in 2-D. Therefore, in our research, we will calculate the point values in triangles based on the intermediate states of 1-D Riemann problem in the normal direction as follows.

\begin{figure}[h!]
	\centering
\includegraphics[width=0.7\linewidth]{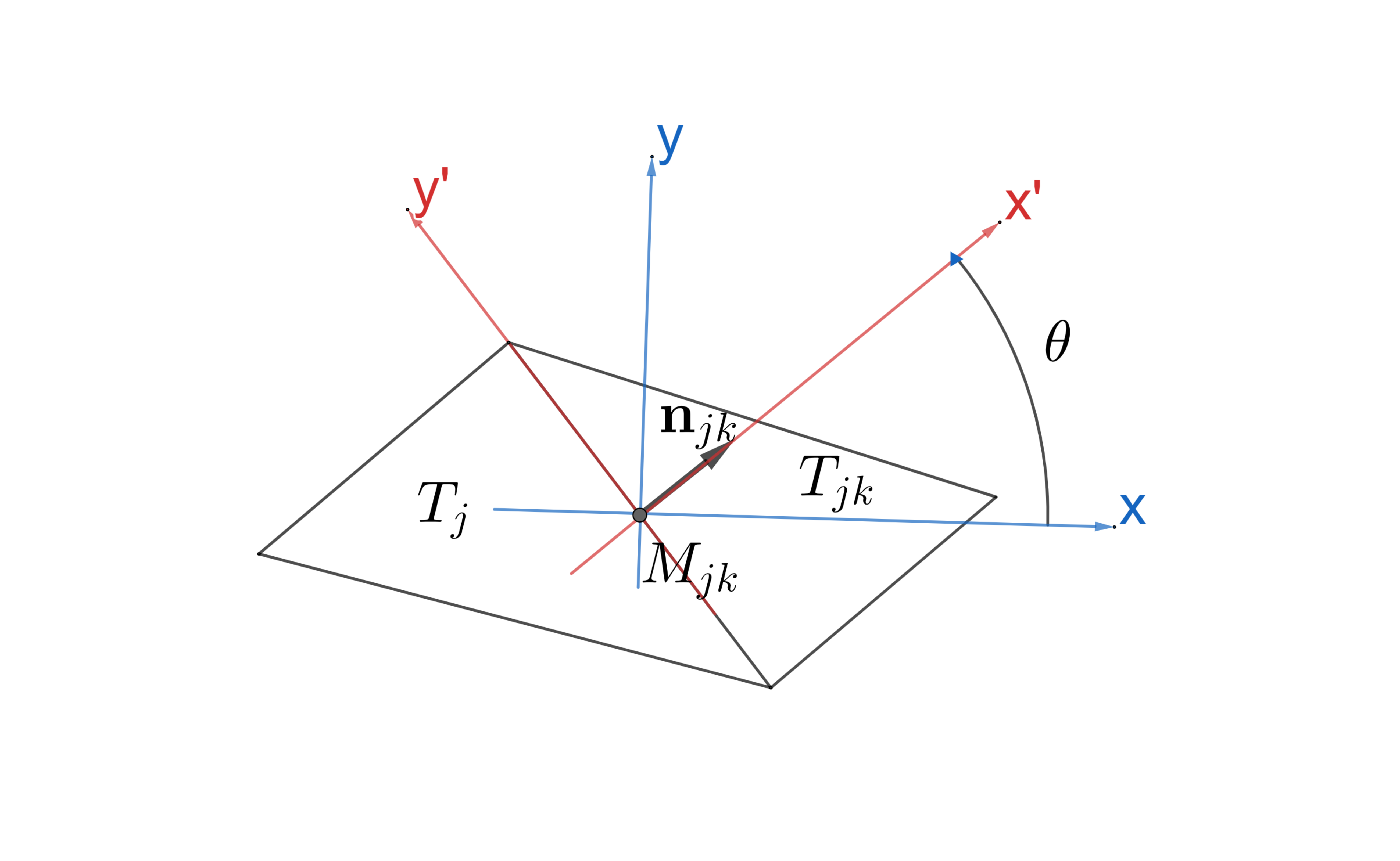}
	\vspace{-1cm}
	\caption{An example of the rotated coordinates $(x',y')$ where $x'$-axis is in the direction of normal vector  $\mathbf{n}_{jk} $ and the new origin is at midpoint $M_{jk}$ of side $k$.}\label{fig:coordinator}
\end{figure}

Assume that triangle $T_j$ is a mixed cell which has the outward unit normal vector $\mathbf{n}_{jk}=(\cos(\theta_{jk}),\sin(\theta_{jk}))$ on side $k$. We will consider a new reference frame $(x',y')$ such that the new horizontal axis is in the direction of the outward normal vector $\mathbf{n}_{jk}$ and the  origin $(0,0)$ is at the middle point $M_{jk}$ of side $k$ as shown in \fref{fig:coordinator}.

Suppose $T_{L}$ and $T_{R}$ to be the two closest single-fluid cells to $T_j$ such that $T_L$ and $T_R$ stay on different sides of the $k$-th edge of cell $T_j$. We assume that  $T_{R}$ stays on the side where the outward normal vector $\mathbf{n}_{jk}$ is pointing. Let $\Ups_L=(w_L,u_L,v_L,\rho_L)$ and $\Ups_R=(w_R,u_R,v_R,\rho_R)$ be the center point values of single-fluid cells $T_{L}$ and  $T_{R}$. Note that $(u_L,v_L)$ and $(u_R,v_R)$ are the values of velocities in the Cartesian coordinate system $(x,y)$. The projections of velocities $(u_L,v_L)$ and $(u_R,v_R)$ onto the rotated frame $(x',y')$ are

 $$u'_L=u_L\cos(\theta_{jk})+v_L\sin(\theta_{jk}),\quad v'_L=-u_L\sin(\theta_{jk})+v_L\cos(\theta_{jk}),$$ $$u'_R=u_R\cos(\theta_{jk})+v_R\sin(\theta_{jk}),\quad  v'_R=-u_R\sin(\theta_{jk})+v_R\cos(\theta_{jk}).$$
 
  Similar to \citep{balsara2014multidimensional, chertock2014central}, to compute the point value at midpoint $M_{jk}$ on side $k$ of triangle $T_j$, we have to solve the following 1-D Riemann problems between states $\Ups_L$ an $\Ups_R$ in direction $\mathbf{n}_{jk}$. 

\begin{equation}
	\begin{cases}
		w_t+(hu)_{x'}=0\\
		(hu)_t+\Big(\dfrac{(hu)^2}{w-B}+\frac{g}{2}\frac{\rho}{\rho_0}(w-B)\Big)_{x'}=-g\dfrac{\rho}{\rho_0}(w-B) B_{x'},\\
		(hv)_t+\left(\dfrac{(hu)(hv)}{w-B}\right)_{x'}=0,\\
		(h\frac{\rho}{\rho_0})_t+\left(uh\dfrac{\rho}{\rho_0}\right)_{x'}=0,	
	\end{cases}\label{eq:Rie}
\end{equation}
subject to the following initial condition

\begin{equation}
	(w,u,v,\rho,B)(x',0)=	\begin{cases}
		(\Ups_L)':=(w_L,u'_L,v'_L,\rho_L,B_L) \quad \mbox{if} \quad x'<0,\\ 
		(\Ups_R)':=(w_R,u'_R,v'_R,\rho_R,B_R) \quad \mbox{if} \quad x'>0.
	\end{cases}\label{eq:intRie}
\end{equation}

The details of the approximate Riemann solver for the above Riemann problem can be seen in \citep{chertock2014central}. Suppose $(\Ups^*_L)':=(w^*_L,(u^*_L)',(v^*_L)',\rho^*_L)$ and $(\Ups^*_R)':=(w^*_R,(u^*_R)',(v^*_R)',\rho^*_R)$ are the intermediate states calculated by Riemann solver for \eref{eq:Rie}-\eref{eq:intRie} in $(x',y')$ coordinates. We recall that in \citep{chertock2014central}, the point values at midpoints on left and right boundaries of a rectangular mixed cell are obtained respectively based on the left and right Riemann intermediate states,$(\Ups^*_L)'$ and $(\Ups^*_R)'$. Here, in triangular grids, we can choose either $(\Ups^*_L)'$ or $(\Ups^*_R)'$ to calculate the point values at the middle point $M_{jk}$ of cell $T_j$. Note that  $T_L$ and $T_R$ are both single-fluid cells nearby mixed cell $T_j$. However, if the neighboring cell $T_{jk}$ shown in \fref{fig:coordinator} is a single-fluid cell, then $T_R\equiv T_{jk}$ and $M_{jk}\in T_R$, which means $T_R$ is closer to $M_{jk}$ compared to $T_L$. Hence, in our work, we select  the right intermediate values $(\Ups^*_R)'$ and the information from $T_R$ to compute the point values $\Ups_j(M_{jk})$. From \citep{chertock2014central}, we first calculate the sound of speed $c^*_R=h^2\dfrac{g\rho}{2\rho_0}$. If $h^*_R>0, \rho^*_R>0$ and $(u^*)'_R-c^*_R<0$, then we replace the piecwise linear approximation in \eref{eq:pwlinear} by $\Ups_j(M_{jk})=(\Ups^*_R)'$, otherwise $\Ups_j(M_{jk})=(\Ups_R)'$ (more details of the point values correction can be seen in \citep{chertock2014central}). Finally, we convert the computed solution back to original Cartesian coordinate $(x,y)$ as 

\begin{equation}
\Ups_j(M_{jk})=\begin{cases}
    \Ups^*_R,\quad &\mbox{if}\quad h^*_R>0, \rho^*_R>0,(u^*)'_R-c^*_R<0,\\
    \Ups_R,\quad &\mbox{otherwise},
\end{cases}
		\label{eq:Riesol}
\end{equation}
where $\Ups^*_R:=(w^*_R,(u^*_R)'\cos(\theta_{jk})-(v^*_R)'\sin(\theta_{jk}),(u^*_R)'\sin(\theta_{jk})+(v^*_R)'\cos(\theta_{jk}),\rho^*_R)$.

From \citep{chertock2014central}, in ``lake at rest'' solutions \eref{eq:lara} and \eref{eq:larb}, we have $\Ups^*_R=\Ups_R$ and the steady state solutions are then preserved by using the discretization of the source term presented in Section \ref{sec:source}. 
\subsection{Positivity-preserving condition}\label{sec:pos}

Next,  we will discuss the restriction of time step to preserve the positivity of water depth $h$ and the variable $h\rho$. 	Consider the Forward Euler (FE) equation 

$$\xbar{{U}}_j^{n+1}=\xbar{{U}}_j^{n}-\frac{1}{|T_j|}\sum_{k=1}^3\Delta
t_{jk}{H}_{jk}+\Delta
t\,\bar{{S}}_j.$$

From \citep{AIPBEKP}, the time step condition to guarantee the positivity of water height $h$ is

\begin{equation}
	\Delta t<\dfrac{1}{6a}\min_{jk}(r_{jk}).\label{eq:CFLSWE}
\end{equation}

One can clearly see  that the time step size \eref{eq:CFLSWE} will also achieve nonnegative water height $h$ in the density shallow water model \eref{eq:swed1}-\eref{eq:swed4}, see \citep{AIPBEKP, LAEK}. We now use the idea from \citep{AIPBEKP} to find the CFL-type condition for the positivity-preserving property of the last variable $h\rho$.
To this end, we apply the forward Euler discretization to the last component of the scheme. 
\begin{equation}
	\begin{aligned}
		(\xbar{h\rho})^{n+1}_j=(\xbar{h\rho})^{n}_j-&\dfrac{\Delta t}{|T_j|}\sum_{k=1}^3\frac{\ell_{jk}\cos(\rho_{jk})}{a_{jk}^{\rm in}+a_{jk}^{\rm out}}
		\Big[a_{jk}^{\rm in}(h\rho u)_{jk}(M_{jk})+a_{jk}^{\rm out}(h\rho u)_j(M_{jk})\Big]\\
		-&\dfrac{\Delta t}{|T_j|}\sum_{k=1}^3\frac{\ell_{jk}\sin(\rho_{jk})}{a_{jk}^{\rm in}+a_{jk}^{\rm out}}
		\Big[a_{jk}^{\rm in}(h\rho v)_{jk}(M_{jk})+a_{jk}^{\rm out}(h\rho v)_j(M_{jk})\Big]\\
		+&\dfrac{\Delta t}{|T_j|}\sum_{k=1}^3\ell_{jk}\frac{a_{jk}^{\rm in}a_{jk}^{\rm out}}{a_{jk}^{\rm in}+a_{jk}^{\rm out}}\big[(h\rho)_{jk}(M_{jk})-(h\rho)_j(M_{jk})\big].
	\end{aligned}\label{eq:pos1}
\end{equation}

Note that

\begin{equation}
	(	\xbar{h\rho})^n_j=h_j \rho_j=\left(\dfrac{1}{3}\sum_{m=1}^3 h_j(M_{jm})\right)\left(\dfrac{1}{3}\sum_{s=1}^3 	\rho_j(M_{js})\right)
	=\dfrac{1}{9}\sum_{m,s}h_j(M_{jm})\rho_j(M_{js}).\label{eq:hpj}
\end{equation}

Plug \eref{eq:hpj} to \eref{eq:pos1}, we have

\begin{align*}	
	(\xbar{h\rho})^{n+1}_j=&\dfrac{1}{9}\sum_{m,s}h_j(M_{jm})\rho_j(M_{js})\\&-\dfrac{\Delta t}{|T_j|}\sum_{k=1}^3\frac{\ell_{jk}\cos(\rho_{jk})}{a_{jk}^{\rm in}+a_{jk}^{\rm out}}
	\Big[a_{jk}^{\rm in}h_{jk}(M_{jk})\rho_{jk}(M_{jk}) u_{jk}(M_{jk})+a_{jk}^{\rm out}h_j(M_{jk})\rho_j(M_{jk}) u_j(M_{jk})\Big]\\
	&-\dfrac{\Delta t}{|T_j|}\sum_{k=1}^3\frac{\ell_{jk}\sin(\rho_{jk})}{a_{jk}^{\rm in}+a_{jk}^{\rm out}}
	\Big[a_{jk}^{\rm in}h_{jk}(M_{jk})\rho_{jk}(M_{jk}) v_{jk}(M_{jk})+a_{jk}^{\rm out}h_j(M_{jk})\rho_j(M_{jk}) v_j(M_{jk}) \Big]\\
	&+\dfrac{\Delta t}{|T_j|}\sum_{k=1}^3\ell_{jk}\frac{a_{jk}^{\rm in}a_{jk}^{\rm out}}{a_{jk}^{\rm in}+a_{jk}^{\rm out}}\big[h_{jk}(M_{jk})\rho_{jk}(M_{jk})-h_j(M_{jk})\rho_j(M_{jk})\big]\\
	=&\dfrac{1}{9}\sum_{m\neq s}h_j(M_{jm})\rho_j(M_{js})\\
	&+\sum_{k=1}^3h_j(M_{jk})\rho_j(M_{jk}) \Big[\dfrac{1}{9}-\dfrac{\Delta t}{|T_j|}\frac{\ell_{jk}a_{jk}^{\rm out}}{a_{jk}^{\rm in}+a_{jk}^{\rm out}}\left(a_{jk}^{\rm in}+u^\perp_j(M_{jk})\right)\Big]\\
	&+\sum_{k=1}^3h_{jk}(M_{jk})\rho_{jk}(M_{jk}) \dfrac{\Delta t}{|T_j|}\frac{\ell_{jk}a_{jk}^{\rm in}}{a_{jk}^{\rm in}+a_{jk}^{\rm out}}\left(a_{jk}^{\rm out}-u^\perp_{jk}(M_{jk})\right).
\end{align*}\label{eq:sum}

From the definitions of the local speeds (\ref{eq:speeds}) we obtain that $u^\perp_{j}(M_{jk})=\cos(\theta_{jk})u_j(M_{jk})+\sin(\theta_{jk})v_j(M_{jk})\leq a_{jk}^{\rm out}$, $0\leq a_{jk}^{\rm in}$, and $0\leq a_{jk}^{\rm out}$. Besides, the corrected reconstruction in Section \ref{sec:scalelimiter} guarantees that $h_{jk}(M_{jk}) \geq 0$ and $\rho_{jk}(M_{jk}) \geq 0$ for all $j$ and $k = 1, 2, 3$. Therefore, all terms in the first and third sum on the RHS
of (\ref{eq:sum}) are nonnegative. To enforce the second sum on the RHS
of (\ref{eq:sum}) to be nonegative, we then need:

$$\dfrac{\Delta t}{|T_j|}\frac{\ell_{jk}a_{jk}^{\rm out}}{a_{jk}^{\rm in}+a_{jk}^{\rm out}}\left(a_{jk}^{\rm in}+u^\perp_j(M_{jk})\right)\leq \dfrac{\Delta t}{|T_j|}\frac{\ell_{jk}a_{jk}^{\rm out}}{a_{jk}^{\rm in}+a_{jk}^{\rm out}}\left(a_{jk}^{\rm in}+a_{jk}^{\rm out}\right)=\dfrac{\Delta t}{|T_j|}l_{jk}a_{jk}^{\rm out}< \dfrac{1}{9},\quad \forall j \mbox{ and } k = 1, 2, 3.$$

Hence we conclude that all terms in the second sum on the RHS of (\ref{eq:sum}) are also nonnegative under the following CFL-type condition.

\begin{equation}
	\Delta t<\dfrac{1}{18a}\min_{jk}(r_{jk})=\dfrac{\min_{jk}\left(\dfrac{2|T_j|}{l_{jk}}\right)}{18\max_{jk}(a_{jk}^{\rm in},a_{jk}^{\rm out})}<\dfrac{1}{9}\min_{jk}\left(\dfrac{|T_j|}{l_{jk}a^{\rm out}_{jk}}\right)\label{eq:CFL},
\end{equation}
where $a=\max_{jk}(a_{jk}^{\rm in},a_{jk}^{\rm out})$ and $r_{jk}=\dfrac{2|T_j|}{l_{jk}}$ is the $k$-th altitude of triangle $T_j$. This completes the proof. This proof is still valid if one uses a higher-order SSP ODE solver (either the Runge-Kutta or the
multistep one), because such solvers can be written as a convex combination of several forward Euler
steps.

\subsection{Well-balanced discretization of source term.}\label{sec:source}

In this section, we first develop a well-balanced discretization of the source terms, which guarantees the first type of ``lake at rest'' steady-state
solutions,
\begin{equation}
	w=\max\big\{C,B(x,y)\big\},\quad C=\rm{Const}, \quad\rho =P\equiv Const, \quad u\equiv v\equiv0, \label{lar}
\end{equation}
are exactly preserved by the resulting central-upwind scheme. This means that the source discretization $\xbar{\mathbf S}_j$ should exactly balance the numerical fluxes so that the right-hand side (RHS) of
\eref{eq:swed1}-\eref{eq:swed4} vanishes at ``lake at rest'' steady states.

To this end, we substitute the ``lake at rest'' state \eref{lar} into FE scheme of the system (it is still correct with the adaptive Runge-Kutta scheme) and conclude that a well-balanced quadrature
for $\xbar{\mathbf S}_j$ should satisfy the following two conditions:
\begin{equation}
	-\frac{g}{2|T_j|}\sum_{k=1}^3\ell_{jk}\cos(\theta_{jk})\cdot\frac{\Delta t^{(\rho)}_{jk}}{\Delta t}\cdot\dfrac{P}{\rho_0}\big[C-B(M_{jk})\big]^2+
	\bar S_j^{\,(2)}=0,
	\label{sp1}
\end{equation}

\begin{equation}
	-\frac{g}{2|T_j|}\sum_{k=1}^3\ell_{jk}\sin(\theta_{jk})\cdot\frac{\Delta t^{(\rho)}_{jk}}{\Delta t}\cdot\dfrac{P}{\rho_0}\big[C-B(M_{jk})\big]^2+
	\bar S_j^{\,(3)}=0,
	\label{sp2}
\end{equation}
where
$$
\bar S_j^{\,(2)}\approx-\frac{g}{|T_j|\rho_0}\iint\limits_{T_j}P(C-B(x,y))B_x(x,y)\,dxdy,\quad
\bar S_j^{\,(3)}\approx-\frac{g}{|T_j|\rho_0}\iint\limits_{T_j}P(C-B(x,y))B_y(x,y)\,dxdy.
$$

In order to derive a well-balanced quadrature, similar to \citep{AIPBEKP, LAEK}, we first apply Green's formula,
$\iint_{T_j}{\rm div}\,\mathbf{{\cal G}}\,dxdy=\int_{\partial T_j}\mathbf{{\cal G}}\cdot\mathbf{n}\,ds$, to the vector field
$\mathbf{{\cal G}}=(\dfrac{1}{2}\rho(x,y)(w(x,y)-B(x,y))^2,0)$ and obtain
\begin{equation}
	\begin{aligned}
		-\iint\limits_{T_j}\rho(x,y)(w(x,y)-B(x,y))B_x(x,y)\,dxdy=
		&\sum_{k=1}^3\int\limits_{(\partial{T_j})_k}\rho(x,y)\frac{(w(x,y)-B(x,y))^2}{2}\cos(\theta_{jk})\,ds\\
		&-\iint\limits_{T_j}\rho(x,y)(w(x,y)-B(x,y))w_x(x,y)\,dxdy\\
		&-\iint\limits_{T_j}\rho_x(x,y)\frac{(w(x,y)-B(x,y))^2}{2}\,dxdy,
		\label{ppp}
	\end{aligned}
\end{equation}
where $(\partial{T_j})_k$ is the $k$-th side of the triangle $T_j$, $k=1,2,3$. The
double integrals are approximated using the trapezoidal rule. This results in the following quadrature for $\bar S^{\,(2)}_j$:
\begin{equation}
	\begin{aligned}
		\bar S^{\,(2)}_j&=\frac{g}{2|T_j|}\sum_{k=1}^3\ell_{jk}\cos(\theta_{jk})\cdot
		\frac{\Delta t^{(\rho)}_{jk}}{\Delta t}\cdot\dfrac{\rho(M_{jk})}{\rho_0}\big[w(M_{jk})-B(M_{jk})\big]^2\\
		&-\frac{g}{3\rho_0}\Big[\rho_{j12}(w_{j12}-\widehat B_{j12})\,w_x(V_{j12})+\rho_{j23}(w_{j23}-\widehat B_{j23})\,w_x(V_{j23})+
		\rho_{j13}(w_{j13}-\widehat B_{j13})\,w_x(V_{j13})\Big].\\
		&-\frac{g}{6\rho_0}\Big[\rho_x(V_{j12})(w_{j12}-\widehat B_{j12})^2\,+\rho_x(V_{j23})(w_{j23}-\widehat B_{j23})^2\,+
		\rho_x(V_{j13})(w_{j13}-\widehat B_{j13})^2\,\Big].
		\label{eq:qra1}
	\end{aligned}
\end{equation}
A similar quadrature for $\bar S^{\,(3)}_j$ is
\begin{equation}
	\begin{aligned}
		\bar S^{\,(3)}_j&=\frac{g}{2|T_j|}\sum_{k=1}^3\ell_{jk}\sin(\theta_{jk})\cdot
		\frac{\Delta t^{(\rho)}_{jk}}{\Delta t}\cdot\dfrac{\rho(M_{jk})}{\rho_0}\big[w(M_{jk})-B(M_{jk})\big]^2\\
		&-\frac{g}{3\rho_0}\Big[\rho_{j12}(w_{j12}-\widehat B_{j12})\,w_y(V_{j12})+\rho_{j23}(w_{j23}-\widehat B_{j23})\,w_y(V_{j23})+
		\rho_{j13}(w_{j13}-\widehat B_{j13})\,w_y(V_{j13})\Big].\\
		&-\frac{g}{6\rho_0}\Big[\rho_y(V_{j12})(w_{j12}-\widehat B_{j12})^2\,+\rho_y(V_{j23})(w_{j23}-\widehat B_{j23})^2\,+
		\rho_y(V_{j13})(w_{j13}-\widehat B_{j13})^2\,\Big].
		\label{eq:qra2}
	\end{aligned}
\end{equation}

Notice that the piecewise linear reconstruction  procedure ensures that at the steady state (\ref{lar}), $\nabla\rho(V_{j\kappa})=0, u_x=v_y=0$ and $\nabla w(V_{j\kappa})=0$ throughout the entire computational domain. This implies that
$(w_{j\kappa}-\widehat B_{j\kappa})\,w_x(V_{j\kappa})\equiv(w_{j\kappa}-\widehat B_{j\kappa})\,w_y(V_{j\kappa})\equiv0$ and $\rho_x(V_{j\kappa})(w_{j\kappa}-\widehat B_{j\kappa})^2\equiv \rho_y(V_{j\kappa})(w_{j\kappa}-\widehat B_{j\kappa})^2\equiv 0$. Therefore, the
quadratures (\ref{eq:qra1}) and (\ref{eq:qra2}) satisfy the desired well-balanced requirements (\ref{sp1}) and (\ref{sp2}).

Next, we consider the ``lake at rest'' situation  

\begin{equation}
	B\equiv\rm{Const}, \quad h^2\rho=Q \equiv Const, \quad u\equiv v\equiv0.  \label{lar1}
\end{equation}
Note that in the region occupied by only one fluid, the density  $\rho$ is a constant and water surface $w$ is then also a constant based on \eref{lar1}. Hence, in single-fluid cells, the steady state \eref{lar1} is equivalent to the solution \eref{lar} which is maintained by the discretization of source term \eref{eq:qra1}-\eref{eq:qra2}. In addition, in mixed cells, the point values are obtained by the data from nearby single-fluid cells thanks to the Riemann Solver approximation presented in Section \ref{sec:riemannsolver}, see also \citep{chertock2014central}. Therefore, the discretization of source term \eref{eq:qra1}-\eref{eq:qra2}, is also capable of preserving the steady state solution \eref{lar1}.

\section{Interface Tracking}\label{sec:interfacetracking}

To achieve an efficient adaptive scheme for multifluid flows, it is essential to accurately capture  the curve where two types of fluid join simultaneously 
with the flow field evolution. The interface tracking is important in  preventing excessive numerical diffusion of variable density, see Section \ref{sec:correction}. We also need the location of density jumps to exactly update the cell averages in the adaptive mesh reconstruction, Section \ref{sect3}.  Therefore, we desire a simple and effective interface reconstruction for the system \eref{eq:swed1}-\eref{eq:swed4}.  Over the last two decades, many methods have been proposed for this purpose such as the level set method \citep{osher2001level, sethian2003level}, the volume of fluid method \citep{hirt1981volume, rider1998reconstructing}, and the front tracking method \citep{unverdi1992front, tryggvason2001front}. Among various versions of interface tracking, we have considered the approach described in  \citep{yang2006adaptive, yang2006analytic} for our work due to its simplicity, robustness, and high efficiency. In particular, the interface is obtained by using  both the level set and the volume fraction functions. 
The details of the interface reconstruction is described as follows.
\subsection{Mixed Cell Detecting by Level Set Function}\label{sec:mixcelltrack}

In the first part of this section, we will discuss the method of detecting the mixed cells in the triangular grids. For this end, we consider the level set function $\phi$ which is defined such that it is positive in one fluid, is negative in the other fluid, and has zero value at the interface. Very often, the level set value of each grid point is initialized by the signed distance from that point to the curve of density discontinuities, see \citep{LAEK, osher2001level, sethian2003level, yang2006adaptive}. As discussed in \citep{LAEK, osher2001level, sethian2003level, yang2006adaptive}, level
set function is evolved by the velocity $(u,v)$ of the flow field as
\begin{equation}
	\phi_t + u\phi_x + v\phi_y = 0.\label{eq:phi}
\end{equation}

 The equation \eref{eq:phi} can be rewritten in conservation form as follows.

\begin{equation}
	\phi_t + (u\phi)_x + (v\phi)_y = (u_x + v_y)\phi.\label{eq:phis}
\end{equation}

 We have used the central-upwind scheme presented in Section \ref{sec:centralupwind} to solve for $\phi$, where the point values $\phi(M_{jk}), k=1,2,3$ in triangle $T_j$ are approximated by the piecewise linear reconstruction \eref{eq:pwlinear}. The integral of the source term in  \eref{eq:phis} can be approximated by midpoint rule. 
 Note that using the central-upwind scheme to solve \eref{eq:phi} only give us the cell averages of $\phi$ in each cell and does not point out which cells contain the interface $\phi=0$. We have to provide a numerical method for detecting mixed cells. However, a triangle is a single-fluid cell, also called as ``reliable'' cell, if the level set values at its vertices are either all positive or all negative. Hence, we will locate  each node in the grid based on its point value of level set as presented below.
 
  \par Without the loss of generality, we assume that $\phi(x,y)$ is positive in fluid 1 and $\phi(x,y)$ is negative in fluid 2. The point value of level set function $\phi^*_{j\kappa}$ at each vertex $V_{j\kappa}$ of cell $T_j$ is approximated by extrapolation \citep{chertock2014central}
\begin{equation}
\phi^*_{j\kappa}=\dfrac{\displaystyle\sum_{i=0}^{m^*_\kappa} c_i\phi^i_{j\kappa}}{\displaystyle\sum_{i=0}^{m^*_\kappa} c_i},\label{eq:phipoint}
\end{equation}
where $\phi^i_{j\kappa}, $ is the cell averages of level set function in cells $T^i_{j\kappa}$  which have common vertex at $V_{j\kappa}$, and the weight $c_i$ is inversely proportional to the distance between the center of cells $T^i_{j\kappa}, i=0,1,...,m^*_\kappa$ and vertex $V_{j\kappa}$. If $\phi^*_{j\kappa}>0$, we have the vertex $V_{j\kappa}$ staying in fluid 1. Otherwise, if $\phi^*_{j\kappa}<0$, the vertex $V_{j\kappa}$ is in fluid 2. Finally, based on the physical meaning of the level set function, a triangle $T_j$ is naturally marked as a mixed cell if it has two vertices staying in different fluids.

\subsection{Interface Reconstruction }\label{sec:interface}

In most works of interface tracking, see \citep{osher2001level, sethian2003level, hirt1981volume, rider1998reconstructing, unverdi1992front, tryggvason2001front}, the curve of density discontinuity in a mixed cell is approximated as a linear line segment of the form $\mathbf{n}\cdot x=\alpha$, where $\mathbf{n}$ is the normal
vector of the interface, $x$ is the location of a point, and $\alpha$ is the line constant. A variety of methods for computing normal vector $\mathbf{n}$ and  parameter $\alpha$ have been briefly introduced and compared in \citep{yang2006adaptive, rider1995stretching}. Based on the advantages and disadvantages of conventional interface tracking methods, an innovative approach has been proposed in \citep{yang2006adaptive}. The approach employs both the level set and volume of fluid methods to denote the density segment in each flagged cell by its two endpoints. Namely, the interface normal vector is calculated from the
level set function while the exact location of two endpoints of the segment is determined by enforcing mass conservation based on the volume fraction. In our work, we will apply this interface tracking method  in the adaptive triangular mesh due to its simplicity, accuracy, and versatility. The details of this method can be seen in \citep{yang2006adaptive}. In the following, we will briefly present  the main steps of the interface reconstruction. 

Once the mixed cells are detected, we will first compute the interface normal vector in each flagged cell based on the level set function and a least square problem. 
 Namely, for each mixed cell $T_j$ with center at $(x_j, y_j)$, we consider a stencil that consists of all centers $(x_i,y_i), i=1,2,...,N$ of $N$ cells that share at least a vertex with $T_j$ (to simplify the computation, we do not place the vertices of  cell $T_j$ in the stencil as in \citep{yang2006adaptive}). A quadratic function for function $\phi$ of a point $(x,y)$ is then given in the generic form

$$\phi=ax'^2 +bx'y'+ cy'^2 + dx' +ey' + g,$$
where $(x':=x-x_j,y':=x-x_j)$ is the location of point $(x,y)$ in the local coordinate frame $x'$ and $y'$ obtained by shifting the original coordinate such that the new origin is at the center of $T_j$. The coefficients $a, b, c, d, e,$ and $g$ are determined by using the least squares method with the linear system 
$$Qs=r,$$
where 	\begin{equation}
	Q=\begin{bmatrix}
		x'^2_1 &x'_1y'_1 &y'^2_1& x'_1 &y'_1& 1\\	
		x'^2_2 &x'_2y'_2 &y'^2_2& x'_2 &y'_2& 1\\
		\vdots&\vdots&\vdots&\vdots&\vdots&\vdots\\
		x'^2_N &x'_Ny'_N &y'^2_N& x'_N &y'_N& 1\\
	\end{bmatrix},\quad s=\begin{bmatrix}
		a\\b\\c\\d\\e\\g
	\end{bmatrix}, \quad \phi=\begin{bmatrix}
		\phi_1\\\phi_2\\\vdots\\\phi_n
	\end{bmatrix},  \label{eq:leastsquare}
\end{equation}
$(x'_i=x_i-x_j,y'_i=y_i-y_j), i = 1,. . . ,N,$ are the local coordinates of the $i-$th node in the stencil and $\phi_i, i = 1,. . . ,N,$ is the cell average of level set function in cell $T_i$ corresponding to node $i$-th. Once the coefficients are known, we compute the unit normal vector by
$\mathbf{n}=\left(\dfrac{d}{\sqrt{d^2+e^2}},\dfrac{e}{\sqrt{d^2+e^2}}\right).$
The system \eref{eq:leastsquare} is solvable as explained in \citep{yang2006adaptive}. In addition, since the level set function is continuous, the normal calculation is second-order accurate.
\par
\begin{figure}[h!]
	\centering
	\includegraphics[width=8cm]{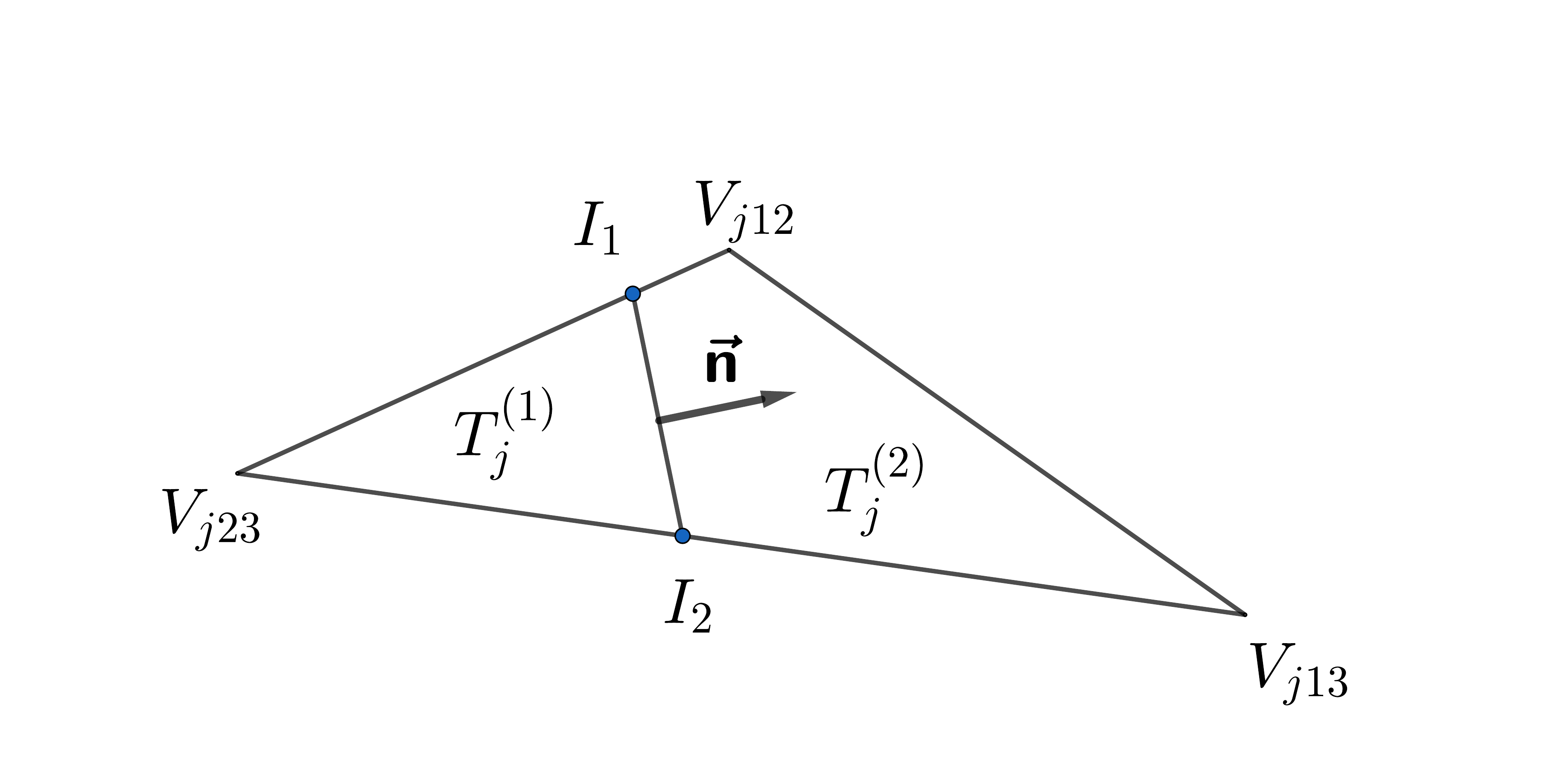}	\includegraphics[width=8cm]{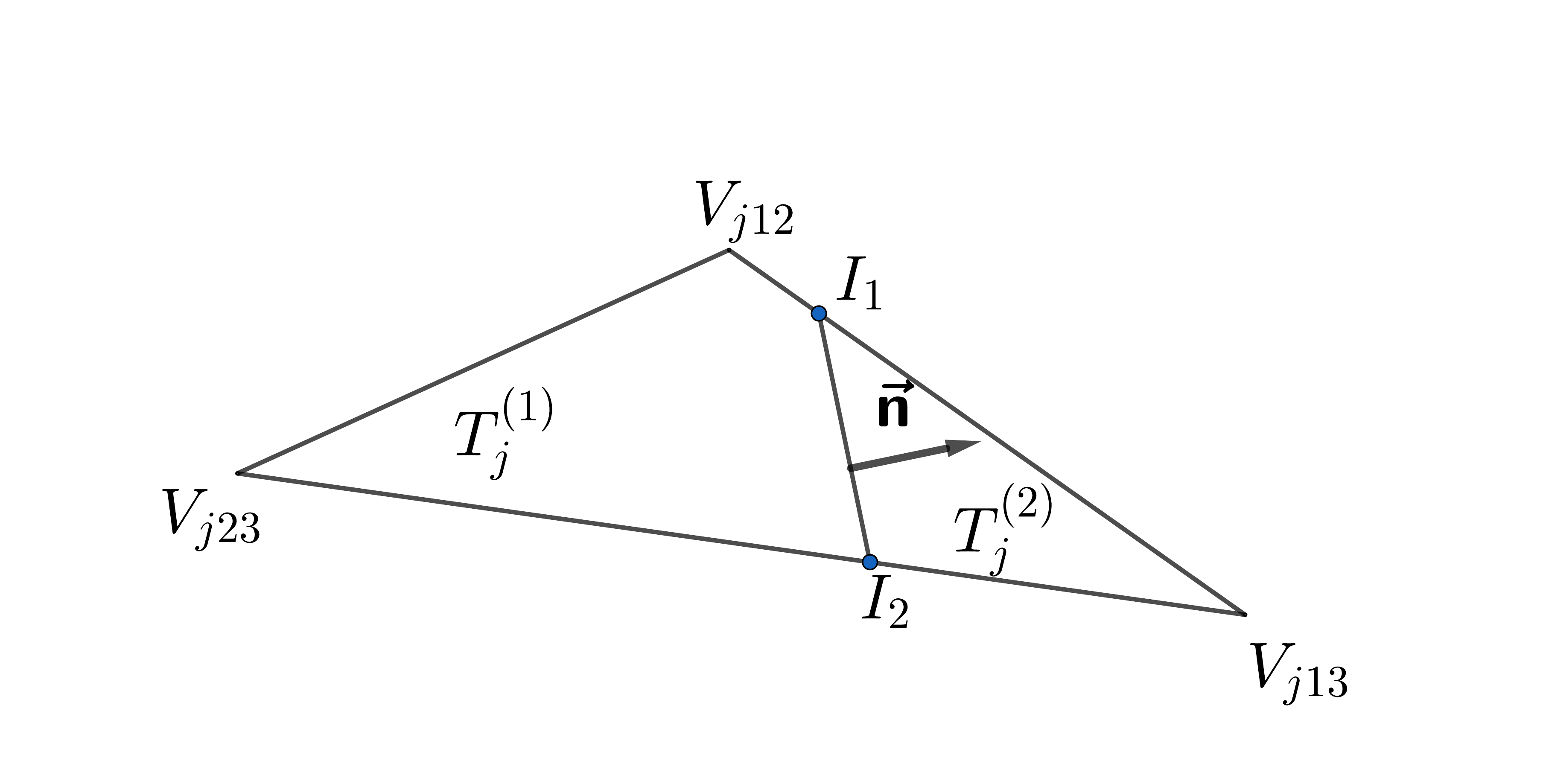}\\
	\caption{{An example of mixed cell $T_j=V_{j12}V_{j23}V_{j13}$ with the interface segment $I_1I_2$ and the normal vector $\vec{\mathbf{n}}$.\label{fig:mixedcell}}}
\end{figure}

 Next, we continue to determine the exact coordinates of two endpoints of the interface by using the volume of fluid approach (VOF). The VOF method has been widely used and developed for interface tracking in \citep{hirt1981volume, rider1998reconstructing, yang2006analytic, yang2006adaptive}. In the VOF method, the volume fraction  function, denoted by $f$, is defined as the ratio of the volume of one fluid, called fluid 1, in each cell to the total volume of the cell. Hence, $f$ is unity if the cell is a single-fluid cell in fluid 1, and is zero if the cell only contains fluid 2. For mixed cells, we have $0<f<1$. The function $f$ is advected by 
\begin{equation}
	f_t + uf_x + vf_y = 0.\label{eq:f}
\end{equation}
The main idea of the VOF method proposed in \citep{yang2006analytic, yang2006adaptive} is to reconstruct the interface segment by determining the coordinates of its endpoints. In particular, the two endpoints must form a segment which has the normal vector $ \mathbf{n}$ as computed above and the interface truncates the cell with the given volume fraction. \fref{fig:mixedcell} is an illustration for two cases of a mixed cell $T_j$ where the interface $I_1I_2$ splits $T_j$ into two parts, $T^{(1)}_j$ and $T^{(2)}_j$, respectively occupied by fluid 1 and fluid 2.

 We then have \begin{equation*}
I_1I_2\cdot \mathbf{n}=0, \quad   \dfrac{T^{(1)}_j}{T_j}=\xbar{f}_j,
\end{equation*}
where $\mathbf{n}$ is the normal vector of the interface computed by \eref{eq:leastsquare} and $\xbar{f}_j$ is the cell average of volume fraction in mixed cell $T_j$ obtained by using the central-upwind method, see Section \ref{sec:centralupwind}, to solve \eref{eq:f}.
This interface reconstruction  is explicit, accurate, and capable of conserving the volume of the fluid. The details of endpoint calculation  can be seen in \citep{yang2006adaptive}. Finally, the reconstructed interface will be used in the adaptive mesh reconstruction as discussed in the following sections.
\subsection{The Cell Averages Correction}\label{sec:correction}

The idea of correcting the solution in the neighborhood of the interface has been used in numerical methods for multifluid flows \citep{chertock2014central, yang2006adaptive} to improve the computed solution. This technique prevents the diffusion of density emerging when we use the central-upwind method to solve  compressible systems. In our work, we will also perform the correction procedure  such that the local conservation is ensured based on the location of density jumps. Namely, suppose we have two types of fluid, fluid 1 and fluid 2, in the flow. After determining the single-fluid cells and mixed cells by using the point values of the level set function, see
Section \ref{sec:mixcelltrack}, the cell average correction will proceed as follows.
\begin{itemize}
    \item If a cell $T_j$ is a single fluid cell in fluid $i=1,2$, we correct the cell averages of $h\rho$ by $\xbar{h\rho}^{new}_j=h_j\rho^{(i)}_j$, where $\rho^{(i)}_j$ is the value of density in fluid $i=1,2$. $h_j=\xbar{w}_j-B_j$ is the cell average of water depth in $T_j$. To preserve the mass of $h\rho$ in the domain, we equally split  the change of $\xbar{h\rho}$ in $T_j$ into the nearby mixed cells. Namely, if the neighboring cell $T_{jk}, k=1,2,3$ of $T_j$ is a mixed cell, we obtain the new cell average of $h\rho$ in $T_{jk}$ by
\begin{equation}
\xbar{h\rho}^{new}_{jk}=\xbar{h\rho}_{jk}+\dfrac{\xbar{h\rho}_j-\xbar{h\rho}^{new}_j}{n_{mix}},\label{eq:corr}  
\end{equation}
where $\xbar{h\rho}_{j}$ and $\xbar{h\rho}_{jk}$ are the old cell averages computed by using the central-upwind scheme in triangles $T_j$ and $T_{jk}$, and $n_{mix}$ is the number of mixed cells surrounding the cell $T_j$.

\item If $T_j$ is a mixed cell, except the update \eref{eq:corr} from its nearby single-fluid cells, no further correction is needed.
\end{itemize}
Due to the fact that the interface normally moves from one cell to its adjacent cells, only cells surrounding the interface are updated  and the  loss of mass from this correction is negligible. In addition, for cells in  the ``lake at rest'' area,  this correction will not change the existing solution therefore maintaining the well-balanced properties.

\section{Adaptive Central-Upwind Scheme}\label{sect3}

 The traditional numerical methods for system
 \eref{eq:swed1}--\eref{eq:swed4} consider very fine
 fixed meshes to reconstruct delicate features of the
 solution. However, this can lead to high computational cost, as well
 as to a poor accuracy of small scale characteristics of the problem. Therefore, we will use adaptive meshes to improve the accuracy of
the approximation at a much lower cost. In this section, we will review an efficient and accurate adaptive central-upwind algorithm from \citep{epshteyn2020adaptive}  which we adapt to the system \eref{eq:swed1}-\eref{eq:swed4}. 

\subsection{Adaptive Central-Upwind Algorithm}\label{sect3a0}
From \citep{epshteyn2020adaptive}, the adaptive central-upwind algorithm for the system \eref{eq:swed1}- \eref{eq:swed4} is described briefly by the following steps.

\textbf{Step 0.} At time $t=t^0$, generate the initial uniform
grid $\mathcal{T}^{0, 0}$.

\textbf{Step 1.} On mesh $\mathcal{T}^{n, \mathcal{M}_n}$, evolve  the cell averages $\xbar{\mU}^n$ of
the solution from time $t^n$ to $\xbar{\mU}^{n+1}$ at the
next time level $t^{n+1} $ using adaptive central-upwind
scheme \eref{eq:rk1}-\eref{eq:rk2}, see Section \ref{sect3}: 	
\begin{itemize}
	\item  At time $t^n$,  determine the level $l=0, 1,
	..., L$ of each triangle
	$T^{n, \mathcal{M}_n}_j \in \mathcal{T}^{n, \mathcal{M}_n}$,
	\eref{eq:lev}, Section \ref{sec:adaptivetime}.
	\item At each time level $t^{n,p}_l, p=0, 1, ..., \mathcal{P}_l-1$,  perform the piecewise polynomial reconstruction
	\eref{eq:pwlinear} for single-fluid cells and apply the Riemann Solver \eref{eq:Rie} on the mixed cells to calculate the point values, Section
	\ref{sec:centralupwind}, Section \ref{sec:correction}. 
	\item At each time level $t^{n,p}_l, p=0, 1, ..., \mathcal{P}_l-1$,  compute the one-sided local speeds of
	propagation using \eref{eq:speeds},
	Section \ref{sec:centralupwind}.
	\item At time $t^n$, calculate the reference time step
	$\Delta t$ using \eref{eq:CFLdtref}, Section \ref{sec:adaptivetime}.
	\item  At each time level $t^{n,p}_l, p=0, 1, ..., \mathcal{P}_l-1$, compute the local time step for
	each cell level, \eref{eq:dtl},  Section \ref{sec:adaptivetime}.
	\item  At each time level $t^{n,p}_l, p=0, 1, ...,
	\mathcal{P}_l-1$, compute numerical fluxes and source term in
	the adaptive central-upwind scheme
	\eref{eq:rk1}-\eref{eq:rk2}, \eref{flux},
	\eref{eq:qra1}-\eref{eq:qra2}, Section \ref{sec:centralupwind},  Section
	\ref{sec:adaptivetime}.			
\end{itemize} 

\textbf{Step 2.} On mesh $\mathcal{T}^{n, \mathcal{M}_n}$, compute WLR error using \eref{ej} in Section
\ref{sect3c} and determine
the refinement/de-refinement status for each cell,
Section \ref{sect3c}.

\textbf{Step 3.} Generate the new adaptive mesh $\mathcal{T}^{n+1, \mathcal{M}_{n+1}}$ at  $t^{n+1}$,
Section \ref{sect3a}. This step includes coarsening of some cells, refinement of
some cells, and the appropriate projection of
the cell averages from the mesh $\mathcal{T}^{n,
	\mathcal{M}_n}$ at $t^n$ onto a new adaptive mesh  $\mathcal{T}^{n+1, \mathcal{M}_{n+1}}$
at time $t^{n+1}$,
Section \ref{sect3a}.

\textbf{Step 4.}  Repeat \textbf{Step 1} - \textbf{Step 3}
until final time.

	\subsection{Adaptive Mesh Refinement/Coarsening}\label{sect3a}

\par
\begin{center}
	\begin{figure}[h!]
		\centering
		\subfigure[Triangle $T^{n,m}_j$ (left) is split into four ``children'' cells $T^{n,m+1}_{js}, s=1,2,3,4$  (right).]{\includegraphics[width=0.8\textwidth]{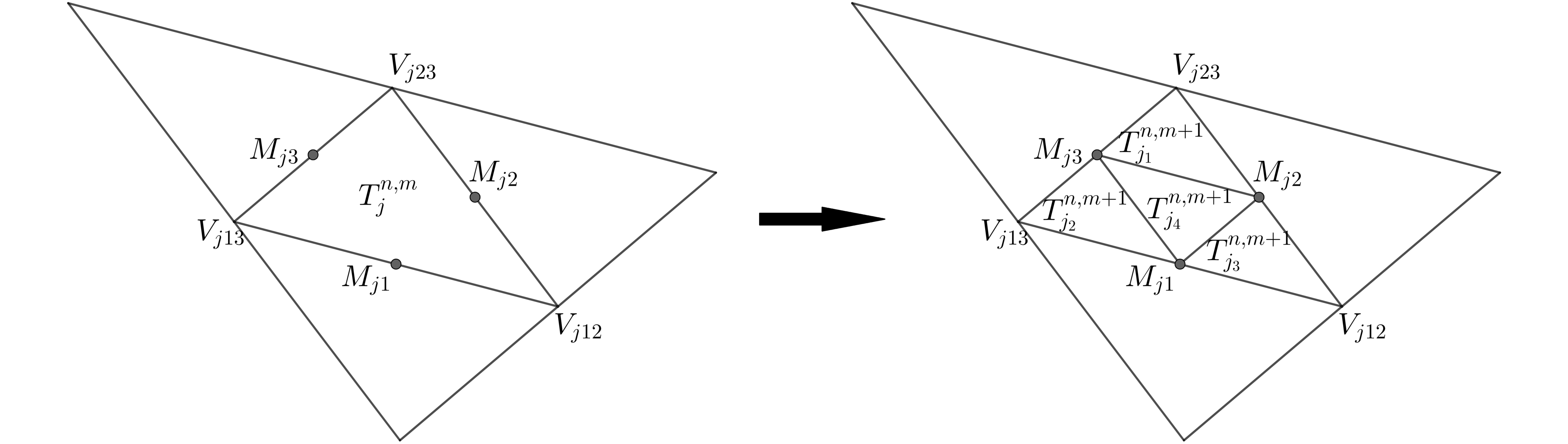}}\label{fig:ref1b}	
		\subfigure[Refinement in the neighboring cells  of $T^{n,m}_j$.]{\includegraphics[width=0.4\textwidth]{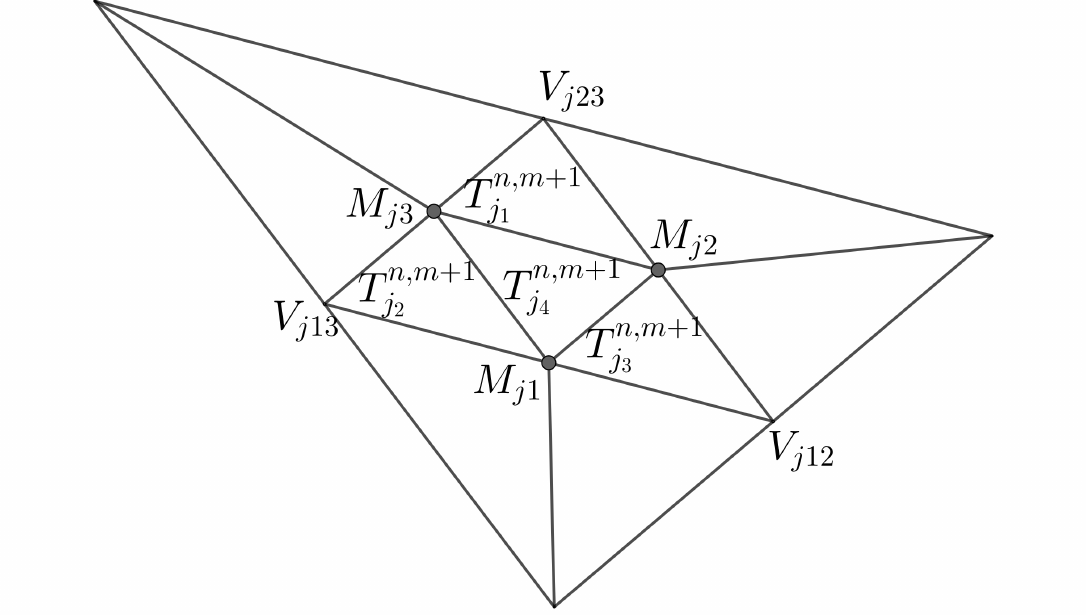}}\label{fig:ref1c}	
		\vspace*{2mm}
		\caption{An outline of the ``regular refinement''.}\label{fig:ref1}
	\end{figure}
\end{center}

 The purpose of mesh reconstruction is to obtain adaptive grids which delivers higher accuracy with a lower computational cost. An efficient adaptive mesh should  have small cells in the regions with large errors and large cells in the other regions. In particular, at time $t^n$, we start with the given mesh, denoted  as
$\mathcal{T}^{n,m}=\{ T^{n,m}_j\}$, where $T^{n,m}_j$ is a
triangular cell with the barycenter $(x^{n, m}_j, y^{n,m}_j)$,  and index $m=0,1,2 ...$ is the level of refinement ($\mathcal{T}^{n,0}\equiv\mathcal{T}^{0,0}$ for all $n$ and $\mathcal{T}^{0,0}$ represents the initial mesh with no refinement). The triangular cells in the mesh
$\mathcal{T}^{n,m}$ are flagged for  refinement/de-refinement based on the weak local residual (WLR)
error estimator,  see Section \ref{sect3c}. On grid $\mathcal{T}^{n,m}$, we apply ``regular
refinement'' described in \citep{epshteyn2020adaptive} on the triangles flagged for refinement  to obtain a new mesh
$\mathcal{T}^{n,m+1}$ with the refinement level $m+1$. Namely, each flagged triangle (``parent'' triangle)  is split into four smaller
triangles (``children'' triangles) by inserting a new node at the
mid-point of each edge of the ``parent'' triangle. \fref{fig:ref1} (a) is an illustration of the ``regular
refinement'', where a flagged cell $T^{n,m}_j$ is refined to obtain the ``children'' cells $T^{n,m+1}_{j_s}, s=1,2,3,4$  by
using the mid-points of the sides. In addition, due to the insertion of new nodes on the edges of the  non-flagged triangles adjacent to refined
triangles, we must also refine these neighboring cells by creating a new edge between the hanging node and the opposite
corner as illustrated in \fref{fig:ref1} (b). 	
\vspace*{-10mm}

In practice, there may be some cells having very large WLR error \eref{ej}, and we need to reach a higher level of refinement for those cells to improve the accuracy. This can be done by repeating the refinement for the flagged triangles in the refined mesh  $\mathcal{T}^{n, m+1}$ to get the mesh with higher level $\mathcal{T}^{n, m+2}, m=0,1,2,....$. \fref{3fig1} is the illustration of the ``regular
refinement'' procedure with two levels of refinement.

\begin{center}\label{3fig1}
	\begin{figure}[h!]
		\centering
		\includegraphics[scale=0.65]{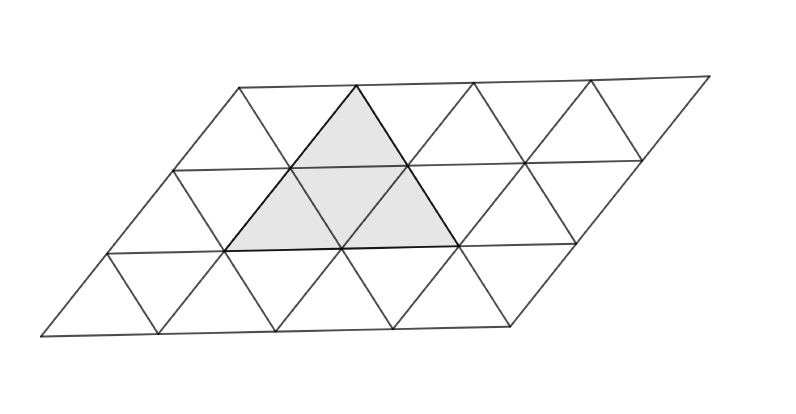}
		\includegraphics[scale=0.65]{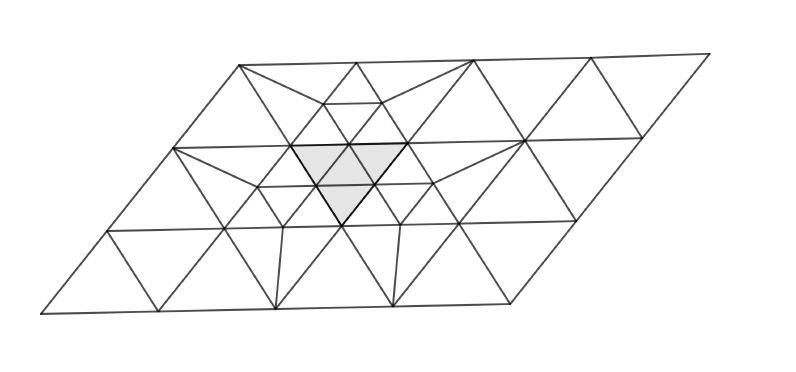}
		\includegraphics[scale=0.65]{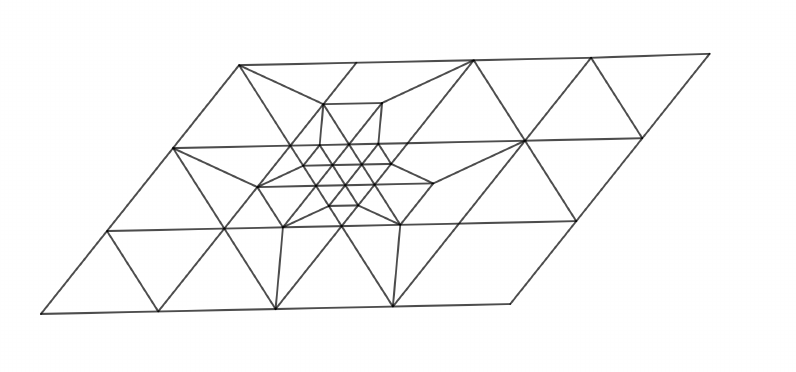}
		\caption{An example of the ``regular
			refinement'' procedure with two levels of refinement. The left figure is the initial coarse mesh
			$\mathcal{T}^{n,0}$ with the region
			flagged for the refinement (gray). The middle figure is the ``first''
			level mesh
			$\mathcal{T}^{n,1}$  with the region
			flagged for higher level of the refinement
			(gray). The right figure is the
			``second'' level mesh $\mathcal{T}^{n,2}$.}\label{fig:mull}
	\end{figure}
\end{center}		

\vspace*{-10mm}
Note that in the numerical simulations of the wave phenomena,
the regions of the domain flagged for refinement changes over time evolution and the refinement in some cells may become no longer
needed. Hence, in \citep{epshteyn2020adaptive}, the de-refinement/coarsening procedure is
performed to deactivate unnecessarily fine cells in the grid. Namely, at time $t^n$, we deactivate
``children'' cells in the mesh $\mathcal{T}^{n,m+1}, m=0,1,...,M_n-1$ based on the WLR and activate the corresponding ``parent'' cell from the mesh
$\mathcal{T}^{n,m}$ back. More details of refinement/de-refinement process can be seen in \citep{epshteyn2020adaptive}.

Finally, at time $t^n$, a
hierarchical system of grids
$\mathcal{S}^n=\{\mathcal{T}^{n,0}, \mathcal{T}^{n,1},
\mathcal{T}^{n,2},..., \mathcal{T}^{n,\mathcal{M}_n}\}$ is obtained, where
$\mathcal{T}^{n,m}$, $m=1,2,...,\mathcal{M}_n$ is the grid
with the level of refinement $m$ reconstructed from the grid
$\mathcal{T}^{n,m-1}$. We will use the final mesh $ \mathcal{T}^{n, \mathcal{M}_n} \in
\mathcal{S}^n$ for
adaptive central-upwind scheme \eref{eq:rk1}-\eref{eq:rk2} to evolve the numerical solution
from time $t^{n}$ to time $t^{n+1}$. Next, at $t^{n+1}$, we generate a new
adaptive grid $\mathcal{T}^{n+1, \mathcal{M}_{n+1}} \in
\mathcal{S}^{n+1}$ from the mesh $\mathcal{T}^{n,
	\mathcal{M}_n}$. After the mesh reconstruction at time $t^{n+1}$, the obtained cell
averages $\xbar{\mU}^{n+1}$ on the mesh $ \mathcal{T}^{n,
	\mathcal{M}_n}$ need to be accurately projected  on the new
mesh $\mathcal{T}^{n+1, \mathcal{M}_{n+1}}$, using the ideas
as summarized briefly below.

\textbf{Case 1.} At $t^{n+1}$, a triangle $T^{n+1,
	\mathcal{M}_{n+1}}_j\in \mathcal{T}^{n+1,
	\mathcal{M}_{n+1}}$ is  the same cell as in the grid
$\mathcal{T}^{n,
	\mathcal{M}_{n}}$, we will maintain the cell
averages for that triangle at $t^{n+1}$.

\textbf{Case 2.} A cell $T^{n+1, \mathcal{M}_{n+1}}_j\in
\mathcal{T}^{n+1, \mathcal{M}_{n+1}}$ is obtained by
de-refining some finer cells $T^{n, \mathcal{M}_n}_{j_s} \in
\mathcal{T}^{n, \mathcal{M}_n}, s=1,2,.., S$. The cell average of solution, $\xbar{U}^{n+1}_j$ in
the cell $T^{n+1, \mathcal{M}_{n+1}}_j$,  is computed as
$$\begin{aligned}
	\xbar{U}^{n+1}_j&=\dfrac{1}{|T^{n+1,\mathcal{M}_{n+1}}_j|}\sum \limits_{s=1}^{S}\xbar{U}^{n}_{j_s}|T^{n,\mathcal{M}_n}_{j_s}|, 
\end{aligned}$$
where $\xbar{U}^{n}_{j_s}$ is the solution in $T^{n,\mathcal{M}_n}_{j_s}$.

\textbf{Case 3.} A triangle $T^{n+1,
	\mathcal{M}_{n+1}}_j\in \mathcal{T}^{n+1,
	\mathcal{M}_{n+1}}$ is obtained from the refinement of the cell  $T^{n,\mathcal{M}_n}_i\in \mathcal{T}^{n,\mathcal{M}_n}$. If $T^{n,\mathcal{M}_n}_i$ is a single-fluid cell, the
cell averages at $t^{n+1}$ in $\mathcal{T}^{n+1,
	\mathcal{M}_{n+1}}$ are approximated by using the  the piecewise linear reconstruction
\eref{eq:pwlinear} of the
solution at $t^{n+1}$ in the triangle
$T^{n,\mathcal{M}_n}_i$. If $T^{n,\mathcal{M}_n}_i$ contains the density discontinuity, we will compute the cells averages in the ``son'' cell $T^{n+1,\mathcal{M}_{n+1}}_j$ based on the interface tracking, see Section \ref{sec:interface}, and the information of the nearby single-fluid cells. For example, suppose that from the reconstructed interface in ``dad'' cell $T^{n,\mathcal{M}_n}_i$, the ``child'' cell $T^{n+1,	\mathcal{M}_{n+1}}_j$ completely lies  in one fluid, called fluid 1. Hence, triangle  $T^{n+1,	\mathcal{M}_{n+1}}_j$  is a single-fluid cell and the cell averages $\xbar{U}^{n+1}_j$ are set to be equal to the cell averages of another ``reliable'' cell  which is in fluid 1 and is the closest cell to the ``dad'' cell $T^{n,\mathcal{M}_n}_i$. 
\begin{remark}
The update of cell averages in case 3 does not ensure the mass conservation which means the total mass of ``son'' cells in the adaptive grid is not equal to the mass of their ``dad'' cell. However, the volume of fluid and the exact location of the interface are conserved. In addition, the reconstructing method also maintains the steady state solution for ``lake at rest'' situations, see example 2 Section \ref{sec:numerical}, since the values of the cell averages are obtained by using the information from nearby ``reliable'' cells. 

\end{remark}
		\subsection{Second-order Adaptive Time Evolution}\label{sec:adaptivetime}
\begin{center}
	\begin{figure}[ht!]
		\centering
		\includegraphics[width=0.7\textwidth]{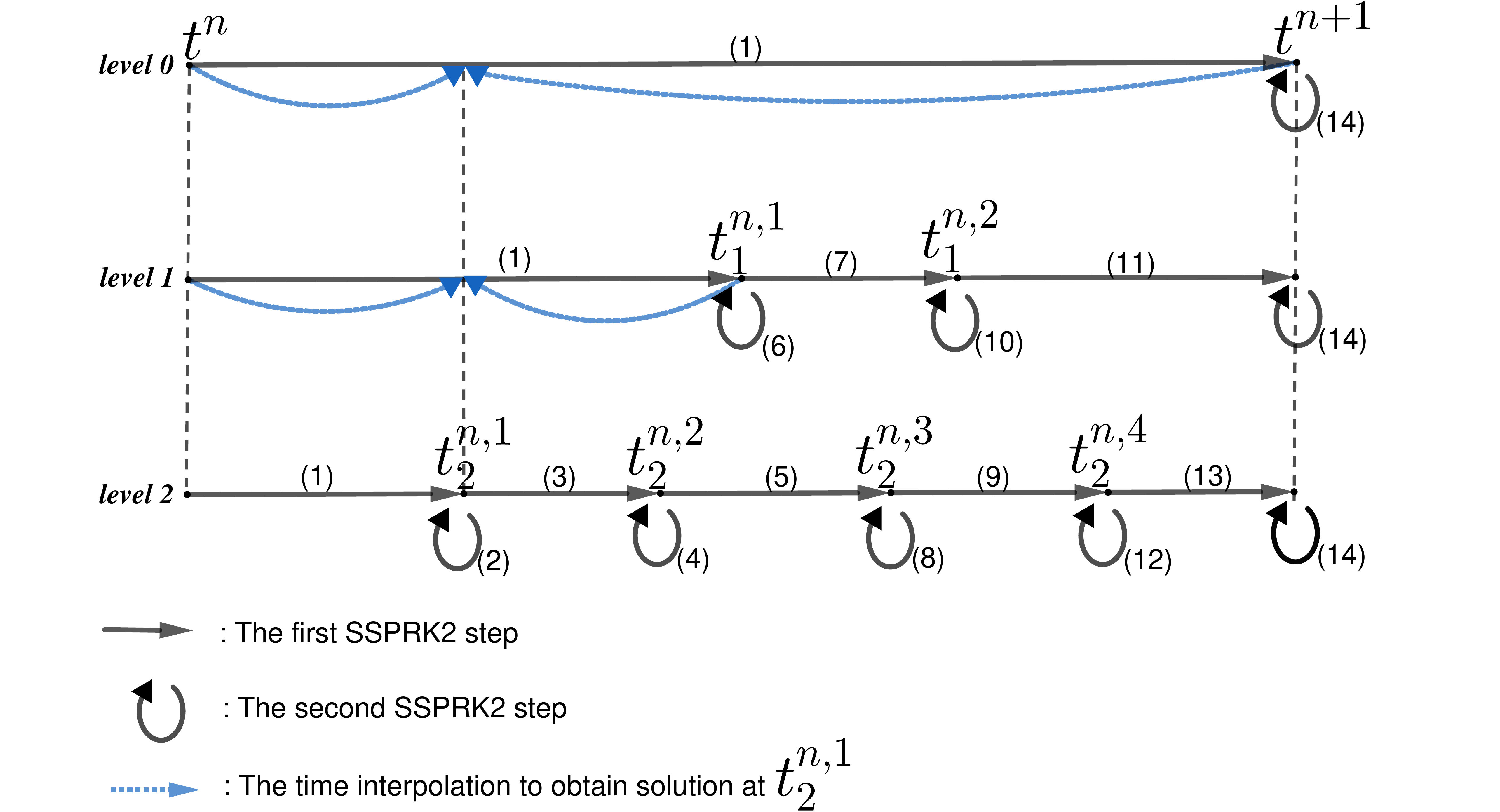}	
		\caption{The example of SSPRK2 on mesh with
			three cell levels,  $l=0,1,2$.}\label{fig:order}
	\end{figure}
\end{center}

Note that using a global time step in the adaptive
algorithm may lead to 
a very small time step due to the presence of much
finer cells in the mesh. To reduce the computational cost, we consider the approach based on the adaptive
time step from \citep{epshteyn2020adaptive, domingues2008adaptive, Donat, MORE}. 
The main idea of this approach is to group cells into different levels based on the cell sizes and applying local time step on each level to evolve from $t^n$ to $t^{n+1}$. Recently, in our work on the shallow water model \citep{epshteyn2020adaptive}, we have derived a simple and efficient
 adaptive time evolution algorithm based
        on  the second-order strong stability preserving Runge-Kutta methods (SSPRK2) in \citep{GST,domingues2008adaptive, Donat, MORE}. This method is also capable to perform on the multi-fluid system \eref{eq:swed1}-\eref{eq:swed2}. The algorithm can be briefly described  by an example as illustrated in \fref{fig:order}. 
\vspace*{-8mm}

First, we group all cells in the grid
$\mathcal{T}^{n,\mathcal{M}_n}$ at time $t^n$
in cell levels $l=0, 1,
.., L$ based on their sizes. Namely, a cell
$T^{n,\mathcal{M}_n}_j\in\mathcal{T}^{n,\mathcal{M}_n}$ belongs to the
level $l$, if $l$ is the smallest positive integer satisfying,	
\begin{equation}
	2^l\geq \dfrac{\displaystyle  \max_j\left(\displaystyle
		\min_{k}\left(r_{jk}\right)\right)}{\displaystyle\min_{k}(r_{j,k})}, \label{eq:lev}
\end{equation}	
where $r_{jk}, k=1,2,3$ are three altitudes of triangle
$T^{n,\mathcal{M}_n}_j$. Next, at $t^n$, we define the
reference time step $\Delta t$ as the local time step on the
coarsest level $l=0$  of cells in the mesh
$\mathcal{T}^{n,\mathcal{M}_n}$ by considering the CFL-type
condition \eref{eq:CFL} locally on level $l=0$. 
\begin{equation}
	\begin{aligned}
		\Delta t\equiv \Delta t^{n,0}_0 =\frac{0.9\displaystyle \max_j\left(\displaystyle \min_{k}\left(r_{jk}\right)\right)}{18a_{max}},
	\end{aligned}
	\label{eq:CFLdtref}
\end{equation}
where
\begin{equation}
	a_{max}:= \displaystyle \max_{j,k}(a^{\rm in}_{jk},a^{\rm out}_{jk}),\label{eq:amax}
\end{equation}	
and $(a^{\rm in}_{jk},a^{\rm out}_{jk})$ are the local
one-sided speeds of propagation \eref{eq:speeds} at $t^n$ for
sides $k=1,2,3$ in the triangle
$T^{n,\mathcal{M}_n}_j\in \mathcal{T}^{n,\mathcal{M}_n}
$.	We set, $t^{n+1}=t^n+\Delta t$.

 Next, assume that $\mathcal{P}_l$ is the
number of steps needed for the higher levels $l=1, .. L$ to evolve from $t^n$ to
$t^{n+1},$ namely $[t^n, t^{n+1}]=\cup[t_l^{n,
	p}, t_l^{n, p+1}], p=0, ..., \mathcal{P}_l-1$ with $t^{n,0}_l\equiv
t^n$, $t_l^{n, \mathcal{P}_l}\equiv t^{n+1} \quad \forall l$.  At $t^{n, p}_l$, the local time step for cells on these levels $l=1, ..., L$ is calculated by
\begin{equation}\Delta t^{n, p}_l=\dfrac{2^{-l} \Delta
		t}{\max(\mu^{n, p}_l,1)},\label{eq:dtl}
\end{equation}
where parameter $\mu^{n,p}_l$ takes into account change in the local one-sided
speeds of the propagation, 
$$\mu^{n, p}_l=\dfrac{	\displaystyle\max_{j,k}(a^{\rm in}_{jk},a^{\rm
		out}_{jk})^{n, p}_l}{a_{max}},$$ where $(a^{\rm
	in}_{jk},a^{\rm out}_{jk})^{n, p}_l$ are the local one-sided
speeds of propagation of the cell
$T^{n,\mathcal{M}_n}_j$ in the level $l$ at $t^{n, p}_l$. Therefore, on
each cell $T^{n,\mathcal{M}_n}_j$ of level $l$,  for each substep $[t_l^{n,p},t_l^{n, p+1}]\equiv [t_l^{n,p},t_l^{n, p}+\Delta
t_l^{n, p}],
p=0,1,2,3,...,\mathcal{P}_l-1$ of the evolution from $t^n$ to  $t^{n+1}$, we apply the following two adaptive steps of the SSPRK2 method, see  \citep{GST,domingues2008adaptive, Donat, MORE},
\begin{subequations}
	\begin{align}
		\xbar{\bm{U}}_j^{(1)}&=\xbar{\bm{U}}_j^{n,
			p}-\Delta
		t^{n, p}_{l}\left(\frac{1}{|T^{n,\mathcal{M}_n}_j|}\sum_{k=1}^3{\bm
			H}^{n,p}_{jk}-\xbar{\bm{S}}^{n,
			p}_j\right):={\bm{R}}(\xbar{\bm{U}}_j^{n,
			p},\Delta t_l^{n, p}),\label{eq:rk1}\\
		\xbar{\bm{U}}_j^{n,
			p+1}&=\dfrac{1}{2}\xbar{\bm{U}}_j^{n,
			p}+\dfrac{1}{2}{\bm{R}}(\xbar{\bm{U}}_j^{(1)},\Delta
		t_l^{n, p}).\label{eq:rk2}
	\end{align}
\end{subequations}
 The flux term ${\bm
	H}^{n,p}_{jk}$ in \eref{eq:rk1}
-\eref{eq:rk2} is the flux
\eref{flux} computed at $t=t^{n,
	p}_l$. The source term $\xbar{{\bm
		S}}^{n,p}_{j}$ in \eref{eq:rk1}
-\eref{eq:rk2} is the source
\eref{eq:qra1} - \eref{eq:qra2}
computed at $t=t^{n, p}_l$ with
the time step $\Delta t_l^{n, p}$.
Note that,
$\xbar{\bm{U}}_j^{n, 0}\equiv
\xbar{\bm{U}}_j^{n}$ and 	$\xbar{\bm{U}}_j^{n,
	\mathcal{P}_l}\equiv \xbar{\bm{U}}_j^{n+1}$.
\begin{remark}  If cells from different cell levels are neighbors, we
	use linear interpolation in time to match the time levels of such
	cells, see also Fig. 3.6, for the illustration of the interpolation. 
\end{remark}	
\subsection{A Posteriori Error Estimator}\label{sect3c}
\par To create a
robust indicator for the adaptive mesh refinement, in our prior work \citep{epshteyn2020adaptive},
we have derived local error estimator from the idea of Weak Local Residual (WLR) presented in
\citep{MR1111445, MR2126235}. This error indicator has shown its advantages in accurately capturing the regions with large error in numerical simulation for Saint-Venant system of shallow water model. Hence, for the adaptive central-upwind scheme for SWEDs, we will extend the error estimator by applying the computation performed in \citep{epshteyn2020adaptive} for the last equations in the system \eref{eq:swed4}.

\par Let us recall that from the weak form of the mass conservation equation \eref{eq:swed1}, in \citep{epshteyn2020adaptive}, the WLR errors
$E^{w,n+\frac1{2}}_{i}$ at each node $N_i$ on mesh $\mathcal{T}^{n, \mathcal{M}_{n}}$ is given by the formula,

\begin{equation}\label{wlrerr}
	\begin{aligned}
		E^{w,n+\frac 1{2}}_{i}&=\frac1{\Delta}(\mathcal{U}^{w,n+\frac 1{2}}_{i}+\mathcal{F}^{w,n+\frac 1{2}}_{i}+\mathcal{G}^{w,n+\frac 1{2}}_{i}),\\\mathcal{ \mU}^{w,n+\frac 1{2}}_{i}&=\sum \limits_{c=1}^{C_i} \frac1{3}|T^{n, \mathcal{M}_n}_{j_c}| (\xbar{w}^{n}_{j_c}-\xbar{w}^{n+1}_{j_c}),\\
		\mathcal{ F}^{w,n+\frac 1{2}}_{i}&=\sum \limits_{c=1}^{C_i}
		a^{(i)}_c\frac{\Delta t}{2}|T^{n, \mathcal{M}_n}_{j_c}| ((\xbar{hu})^{n}_{j_c}+(\xbar{hu})^{n+1}_{j_c}),\\
		\mathcal{ G}^{w,n+\frac 1{2}}_{i}&=\sum \limits_{c=1}^{C_i}
		b^{(i)}_c\frac{\Delta t}{2}|T^{n, \mathcal{M}_n}_{j_c}| ((\xbar{hv})^{n}_{j_c}+(\xbar{hv})^{n+1}_{j_c}),
	\end{aligned}
\end{equation}
where $C_i$ is the number of triangles $T^{n, \mathcal{M}_n}_{j_c}$ having common vertex at node $N_i$. Here, the quantity $(a^{(i)}_c, b^{(i)}_c)$ is the gradient of the
linear piece restricted to $T^{n,\mathcal{M}_n}_{j_c}$,	
\begin{equation}
	\begin{aligned}
		a^{(i)}_c&=\dfrac{\widetilde{y}_2-\widetilde{y}_3}{(\widetilde{y}_3-\widetilde{y}_i)(\widetilde{x}_2-\widetilde{x}_i)-(\widetilde{y}_2-\widetilde{y}_i)(\widetilde{x}_3-\widetilde{x}_i)},\\
		b^{(i)}_c&=\dfrac{\widetilde{x}_3-\widetilde{x}_2}{(\widetilde{y}_3-\widetilde{y}_i)(\widetilde{x}_2-\widetilde{x}_i)-(\widetilde{y}_2-\widetilde{y}_i)(\widetilde{x}_3-\widetilde{x}_i)},\\
	\end{aligned}\label{gradfor}
\end{equation}
where $N_i=(\widetilde{x}_i, \widetilde{y}_i), (\widetilde{x}_2, \widetilde{y}_2), $ and $(\widetilde{x}_3, \widetilde{y}_3) $ are the three vertices of triangle $T^{n,\mathcal{M}_n}_{j_c}$.

Now, by applying the same calculation in \citep{epshteyn2020adaptive} on the weak form of the last equations in the system \eref{eq:swed4}, we then define the WLR error of variable $h\rho$ at node $N_i$ in the grid as,

\begin{equation}\label{wlrerrrho}
	\begin{aligned}
		E^{h\rho,n+\frac 1{2}}_{i}&=\frac1{\Delta}(\mathcal{U}^{h\rho,n+\frac 1{2}}_{i}+\mathcal{F}^{h\rho,n+\frac 1{2}}_{i}+\mathcal{G}^{h\rho,n+\frac 1{2}}_{i}),\\\mathcal{ \mU}^{h\rho,n+\frac 1{2}}_{i}&=\sum \limits_{c=1}^{C_i} \frac1{3}|T^{n, \mathcal{M}_n}_{j_c}| (\xbar{h\rho}^{n}_{j_c}-\xbar{h\rho}^{n+1}_{j_c}),\\
		\mathcal{ F}^{h\rho,n+\frac 1{2}}_{i}&=\sum \limits_{c=1}^{C_i}
		a^{(i)}_c\frac{\Delta t}{2}|T^{n, \mathcal{M}_n}_{j_c}| ((\xbar{hu\rho})^{n}_{j_c}+(\xbar{hu\rho})^{n+1}_{j_c}),\\
		\mathcal{ G}^{h\rho,n+\frac 1{2}}_{i}&=\sum \limits_{c=1}^{C_i}
		b^{(i)}_c\frac{\Delta t}{2}|T^{n, \mathcal{M}_n}_{j_c}| ((\xbar{hv\rho})^{n}_{j_c}+(\xbar{hv\rho})^{n+1}_{j_c}).
	\end{aligned}
\end{equation}

Hence, the error in a cell $T^{n,\mathcal{M}_n}_j\in \mathcal{T}^{n.\mathcal{M}_n}$ takes into account both WLR errors of water surface $w$ and of variable $h\rho$ as, 
\begin{equation}
	e_j=\max_{k}\left(\left|E^{w,n+\frac 1{2}}_{jk}\right|,\left|E^{h\rho,n+\frac 1{2}}_{jk}\right|\right), \quad
	k=1, 2, 3, \label{ej}
\end{equation}
where $E^{w,n+\frac 1{2}}_{jk}$ and $E^{h
	\rho,n+\frac 1{2}}_{jk}$  are the  WLR errors
computed in (\ref{wlrerr}) and (\ref{wlrerrrho}) at a node $k$ of triangle $T_j$.

The
error $e_j$ in each cell $T^{n.\mathcal{M}_n}_j \in \mathcal{T}^{n,
	\mathcal{M}_{n}}$ is compared to the error tolerance
computed by \begin{equation}
\omega=\sigma \max_j(e_j),\label{eq:tol}    
\end{equation}
where $\sigma<1 $ is a given problem-dependent constant. Based on the error comparison, the cell
is then either ``flagged'' for refinement/de-refinement or ``no-change''. 
 
	\section{Numerical Experiments}\label{sec:numerical}
	
\begin{figure}[h!]	
		\centering
		\includegraphics[width=0.3\linewidth]{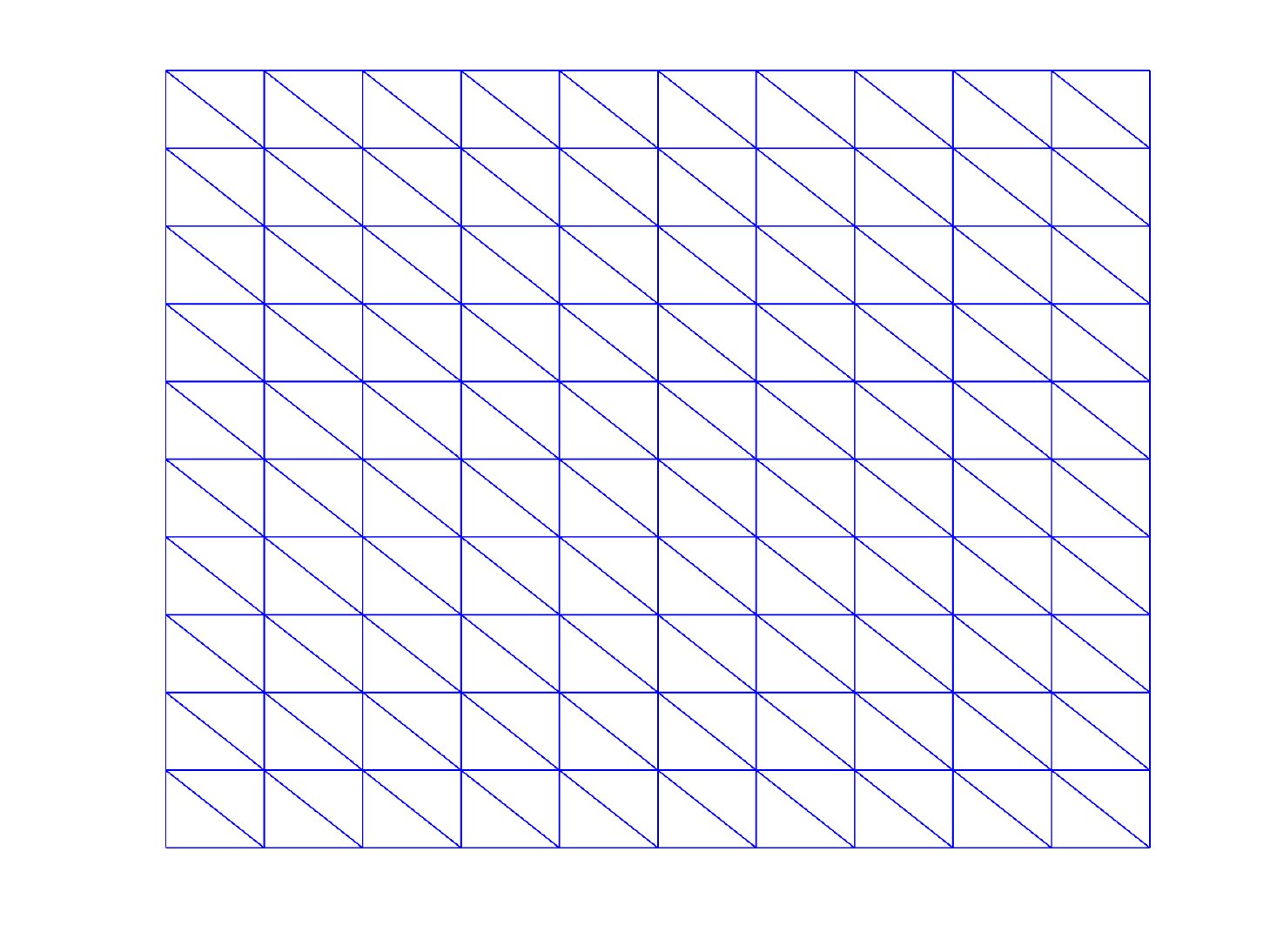}\\
		\vspace*{2mm}
		\caption{An outline of uniform triangular mesh.}
		\label{fig:ex1}
	\end{figure}

	In this section, we will verify the computational efficiency of the
	designed adaptive central-upwind scheme. We compare the results of the adaptive method with the results
	of the central-upwind scheme without the adaptivity, see Section \ref{sect3}, on uniform
	triangular meshes (example of such uniform triangular mesh is
	outlined  in \fref{fig:ex1}). To this end, in all examples,
	we calculate the $L^1$-errors and a ratio,
	$\mathcal{R}_{CPU}=\frac{CPU_{uniform}}{CPU_{adaptive}}$,
	which is the ratio of the CPU times of the central-upwind
	algorithm without the mesh reconstruction to the CPU time of the adaptive
	central-upwind algorithm. To  compare $L^1$-errors, as well as to compare  the CPU times and to
	compute $\mathcal{R}_{CPU}$, we consider uniform mesh and the
	adaptive mesh with the same size of the smallest cells. Namely, in Tables \ref{tab:ex1a}-\ref{tab:ex3t},  $L^1$-errors  and   $\mathcal{R}_{CPU}$ are computed by using the uniform meshes $2\times N\times N, N=100, 200, 400$ and the corresponding adaptive meshes which are reconstructed from the coarser uniform mesh  $2\times N/2^\mathcal{M}\times N/2^\mathcal{M}$ ($\mathcal{M}=1,2$ is the highest level of refinement in the adaptive mesh). To compute the  $L^1$-errors, the reference solution is obtained by applying the central-upwind method without implementing adaptivity techniques on the
	uniform mesh with 
	$2\times 800\times 800$ triangles. In all experiments, we consider a zero-order
	extrapolation at all boundaries. In addition, we use the gravitational acceleration, $g = 1.0$ and the reference density, $\rho_0=997$ \citep{GHAZIZADEH2020104633} for all examples. The
	desingularization parameters $\tau$ and $\varepsilon$ for
	calculations of the velocity components $u$ and $v$ are set
	$\tau=\max_j\{|T^{n,\mathcal{M}_n}_j|^2\}$ and $\varepsilon=10^{-4}$ (see Section 2.1 formula (2.6) in \citep{LAEK}).
	
	\subsection{Example 1:}
	
	\par In the first example taken from \citep{chertock2014central}, we will compare the
	performance of the adaptive central-upwind scheme
	and the central-upwind scheme without adaptivity on uniform triangular meshes. We also verify experimentally the advantages of  the interface reconstruction in compressing the diffusion of variable density at the interface of fluids and improving the accuracy of computed solutions.
	
	\par  We consider the bottom topography $B(x,y,0):=0$ and the following initial condition,
	
	\begin{equation}
	(w,u,v,\rho)^T(x,y,0)=\begin{cases}
	    (2,0,0,1.5\rho_0), \quad \mbox{if}\quad  & x^2+y^2<0.5,\\
	    (1,0,0,\rho_0), \quad &  \mbox{otherwise}.
	\end{cases}\label{eq:ex1}
	\end{equation}
	The data is simulated in the domain $[-1,1]\times[-1,1]$. The error tolerance
	\eref{eq:tol} for the
	mesh refinement in this example is set to $\omega = 0.01\max_j(e_j)$.

\begin{figure}[htp!]
	\centering
	\subfigure[Uniform mesh $2\times 100\times 100$.]{\includegraphics[width=1\textwidth]{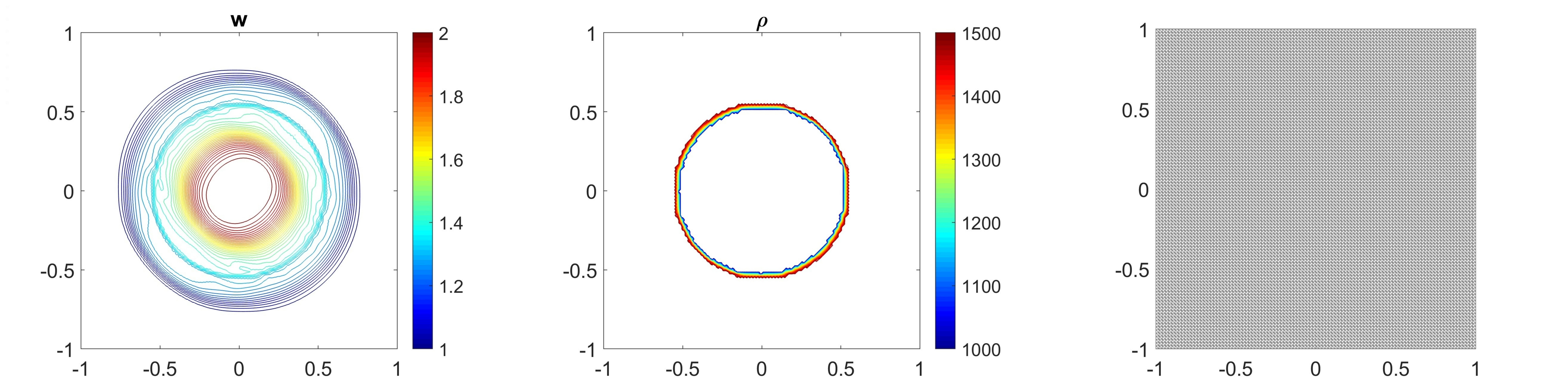}}\label{fig:1a}\\	
	\subfigure[Uniform mesh $2\times 200\times 200$.]{\includegraphics[width=1\textwidth]{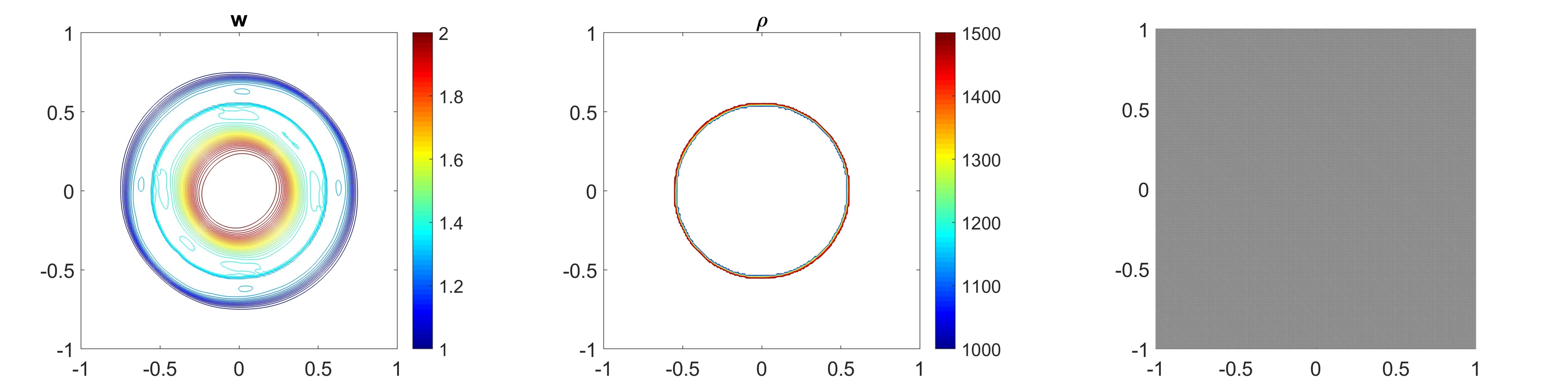}}\label{fig:1b}\\
	\subfigure[Adaptive mesh with $\mathcal{M}=1$.] {\includegraphics[width=1\textwidth]{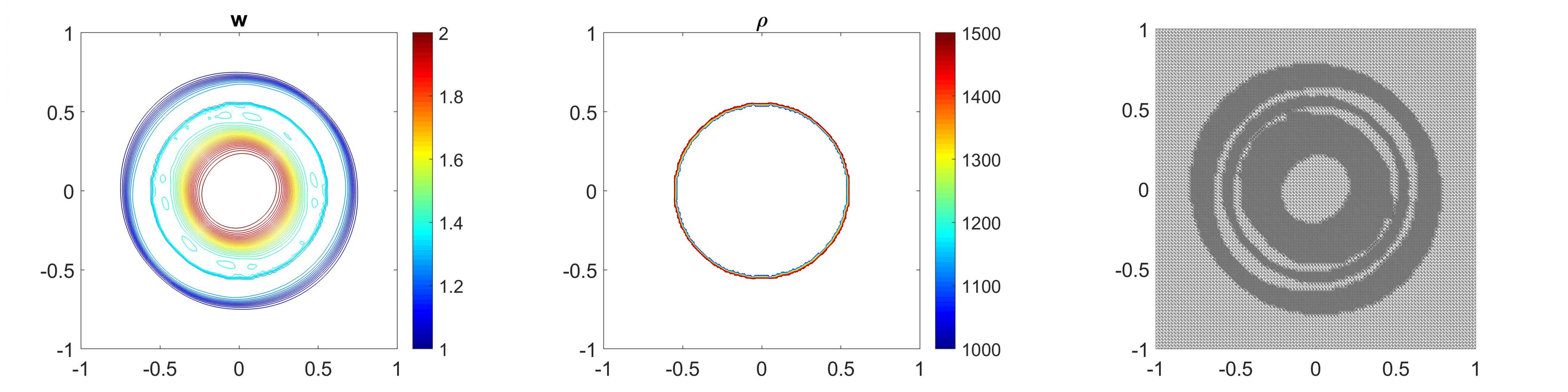}}\label{fig:1c}\\
	\subfigure[Adaptive mesh with $\mathcal{M}=2$.]{\includegraphics[width=1\textwidth]{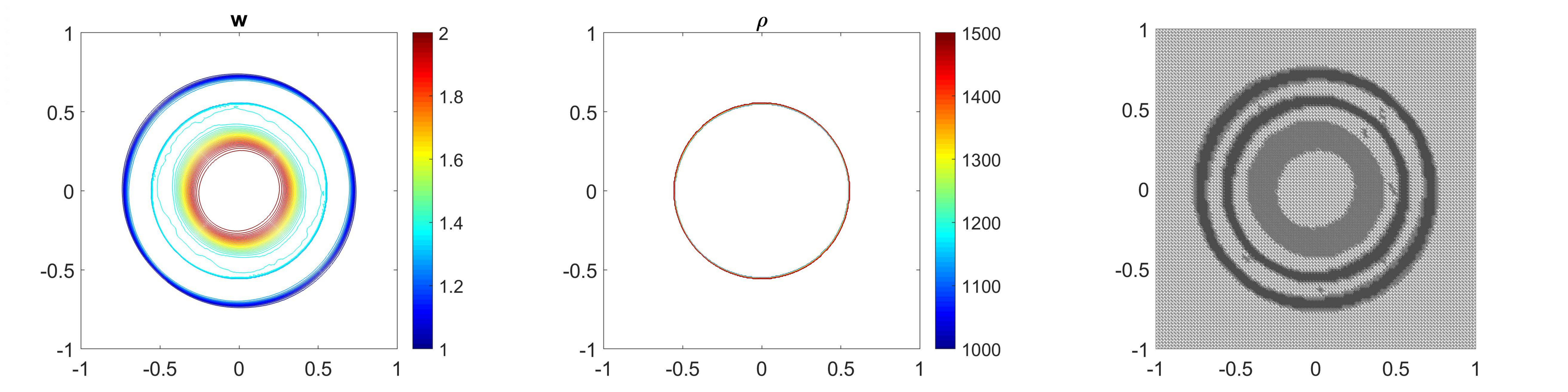}}\label{fig:1d}\\
	\vspace*{2mm}
	\caption{Example 1: Contour plots of computational water surface $w(x,y,0.15)$ (first column) and density $\rho(x,y,0.15)$ (second column) of the IVP \eref{eq:ex1}  with the corresponding meshes (third column).}\label{fig:ex1a}
\end{figure}

\begin{figure}[htp!]
	\centering
	\subfigure[Uniform mesh $2\times 100\times 100$.]{\includegraphics[width=0.7\textwidth]{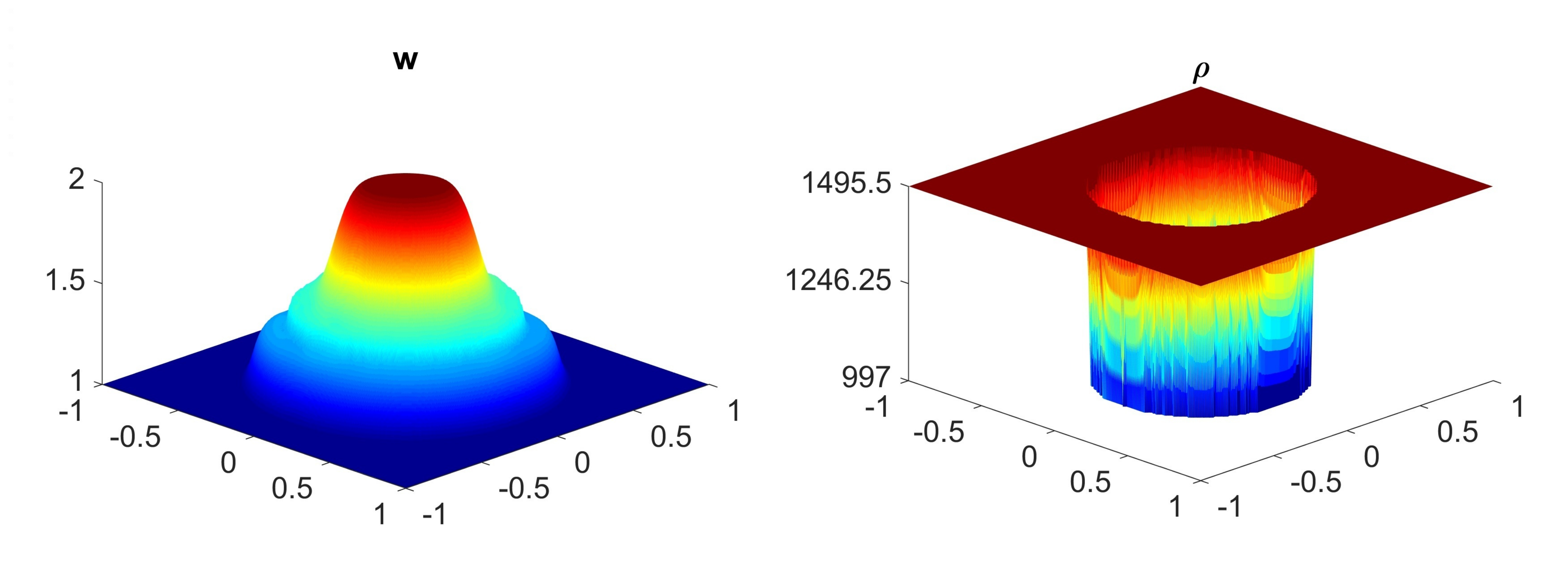}}\label{fig:1as}\\	
	\subfigure[Uniform mesh $2\times 200\times 200$.]{\includegraphics[width=0.7\textwidth]{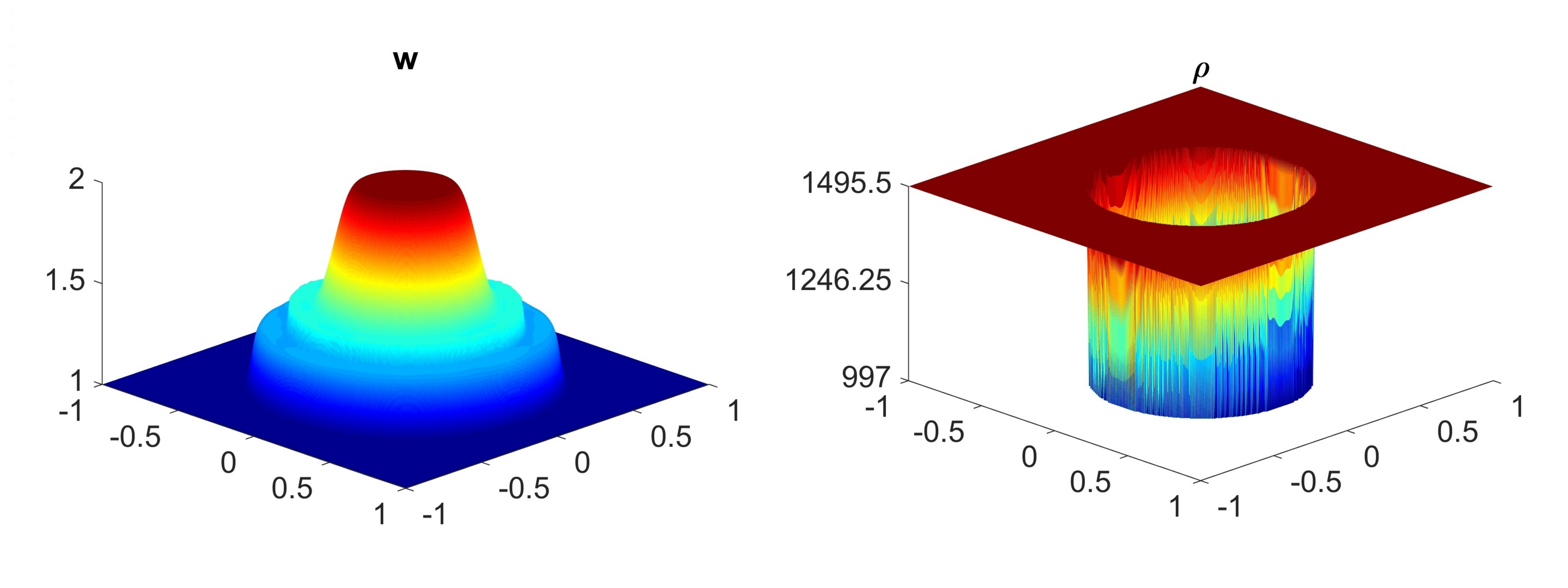}}\label{fig:1bs}\\
	\subfigure[Adaptive mesh with $\mathcal{M}=1$.] {\includegraphics[width=0.7\textwidth]{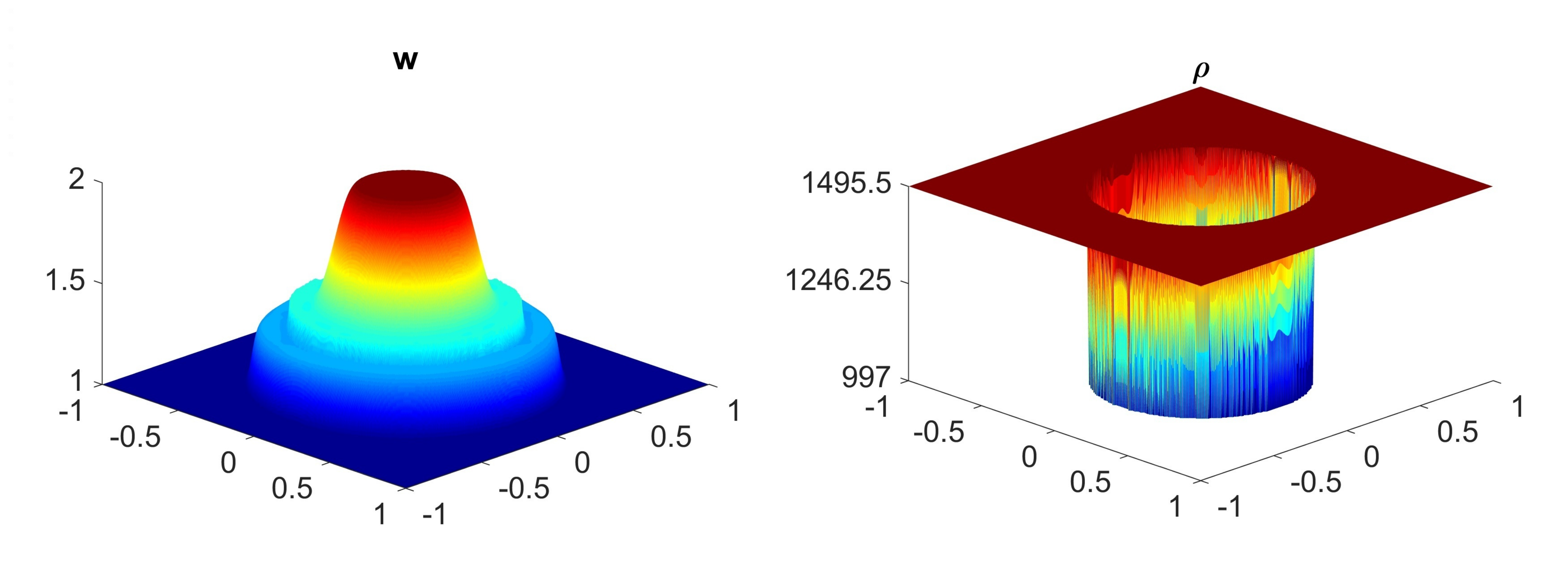}}\label{fig:1cs}\\
	\subfigure[Adaptive mesh with $\mathcal{M}=2$.]{\includegraphics[width=0.7\textwidth]{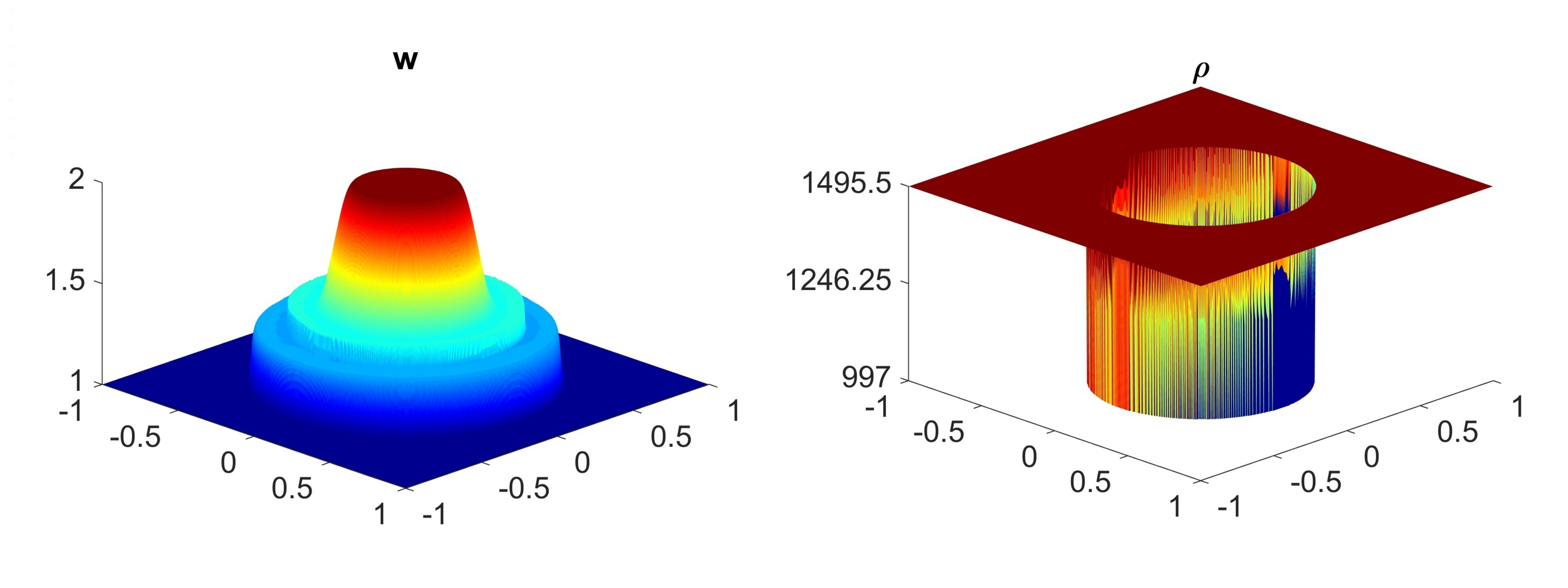}}\label{fig:1ds}\\
	\vspace*{2mm}
	\caption{Example 1: 3-D plots of computational water surface $w(x,y,0.15)$ (first column) and density $\rho(x,y,0.15)$ (second column)  of the IVP \eref{eq:ex1}.}\label{fig:ex1as}
\end{figure}

 In	\fref{fig:ex1a} and \fref{fig:ex1as}, we show the numerical solution of the water surface (first column) and the density (second column) at time $t=0.15$. The solution is calculated by using the central-upwind scheme
	on uniform meshes on \fref{fig:ex1a} (a, b) and  \fref{fig:ex1as} (a,b) and using the
	adaptive central-upwind scheme on \fref{fig:ex1a} (c, d) and  \fref{fig:ex1as} (c,d). The
	adaptive meshes in \fref{fig:ex1a} (third column) are obtained from the
	uniform mesh $2\times 100\times 100$, \fref{fig:ex1a} (a). The
	adaptive mesh  with one level of
	refinement $\mathcal{M}=1$ (as the highest level of
	refinement) is on \fref{fig:ex1a} (c),  and with two levels
	of refinement $\mathcal{M}=2$ (as the highest level of
	refinement) is on \fref{fig:ex1a} (d). As can be observed in  \fref{fig:ex1a}, both $w$ and $\rho$ are much sharper resolved by using the
adaptive central-upwind scheme. We note also,
	that by increasing the level of refinement from $\mathcal{M}=1$ to
	$\mathcal{M}=2$, the number of cells in the mesh increases from
	$40,292$ cells with $\mathcal{M}=1$ to $60,326$ cells with
	$\mathcal{M}=2$, but the accuracy is clearly improved with
	higher resolution as seen in \fref{fig:ex1a} (c) and
	\fref{fig:ex1a} (d). Also from the adaptive meshes in \fref{fig:ex1a} (third column), one can easily see that only cells in the region having steep gradients are refined. This means that the WLR error estimator
accurately detects regions in the domain for adaptive refinement/coarsening.

	Next, in \tref{tab:ex1a} we calculate the $L^1$-errors obtained on the adaptive grids and on the fixed uniform grids. The errors obtained in the uniform meshes are approximate to the errors calculated in the corresponding adaptive meshes (the adaptive meshes have the same size of the smallest cells with the uniform meshes). However, the adaptive scheme uses fewer cells than the central-upwind scheme which does not consider the mesh reconstruction. In addition, in  \tref{tab:ex1ta}, we also compute the $\mathcal{R}_{CPU}$ ratio to compare the
	computational cost of the two methods. The results in \tref{tab:ex1a} and
	\tref{tab:ex1ta} show that the adaptive scheme produces
	similar accuracy as the scheme designed on fixed uniform triangular
	meshes,  but at a  less
	computational cost. 
	
	\begin{table}[htp!]
		\vspace*{5mm}
		\centering
		\begin{tabular}{ |c| c c|c|c c|c|c c|}
			\hline
			\multicolumn{3}{|c|}{algorithm without adaptivity}&\multicolumn{6}{c|}{adaptive algorithm}\\
			\cline{4-9}
			\multicolumn{3}{|c|}{ }
			&\multicolumn{3}{c|}{one level $\mathcal{M}=1$}&\multicolumn{3}{c|}{two levels $\mathcal{M}=2$}\\
			\hline
			cells &$L^1$-error  &rate&cells	&$L^1$-error  &rate&cells	&$L^1$-error  &rate\\
			\hline
			$2\times100\times100$&0.0256&&13,516 &0.0257&&12,800&0.0265&	\\  
			$2\times200\times200$&0.0127	&1.01&40,292 &0.0128&1.00&38,284 & 0.0133&0.99\\
			$2\times400\times400$&0.0047&1.43&153,152 &0.0049&1.39&60,326&0.0050&1.41\\
			\hline	
		\end{tabular}
		\caption{Example 1: $L^1$-errors of the water surface
			$w$  of the IVP \eref{eq:ex1} at $t=0.15$, and the convergence rates of the
			central-upwind scheme without adaptivity (uniform mesh $2\times N\times N, N=100,200,400$) and the
			adaptive scheme (the corresponding adaptive
			mesh is reconstructed from the uniform mesh $2\times N/2^\mathcal{M}\times N/2^\mathcal{M}$).}\label{tab:ex1a}
		\vspace*{2mm}
	\end{table}
	
	\begin{table}[ht!]
		\vspace*{5mm}
		\begin{tabular}{ |c|c| c c|c |c c|}
			\hline
			uniform mesh &\makecell{adaptive mesh\\ $\mathcal{M}=1$}&\multicolumn{2}{c|}{$\mathcal{R}_{CPU}$ with $\mathcal{M}=1$} &\makecell{adaptive mesh\\ $\mathcal{M}=2$}&\multicolumn{2}{c|}{$\mathcal{R}_{CPU}$ with $\mathcal{M}=2$}\\
			\cline{3-4} \cline{6-7}
			(cells)&(cells)&total &\makecell{without\\grid generation} &(cells)&total &\makecell{without\\grid generation}  \\
			\hline
			$2\times100\times100$&13,516&1.81&2.01&12,800&1.83&1.94\\
			$2\times200\times200$&40,292 &1.82&2.14&38,284 & 2.14 &2.30 \\
			$2\times400\times400$&153,152&1.98& 2.29 &60,326 &2.33& 2.50 \\
			\hline
			\multicolumn{2}{|c}{$\mathcal{R}_{CPU}$ average:}&1.87&\multicolumn{1}{c}{2.15 }&\multicolumn{1}{c}{} &2.10& 2.27 \\
			\hline	
		\end{tabular}
		\caption{Example 1: $\mathcal{R}_{CPU}$ ratio for the IVP \eref{eq:ex1} at
			$t=0.15$, where for adaptive central-upwind scheme,
			we consider the total CPU times and CPU times
			without the grid generation.}\label{tab:ex1ta}
	\end{table}
	
	\par
\begin{figure}[ht!]
		\centering
	\subfigure[Adaptive scheme without  interface tracking on adaptive mesh with 51768 cells.]{\includegraphics[width=1\textwidth]{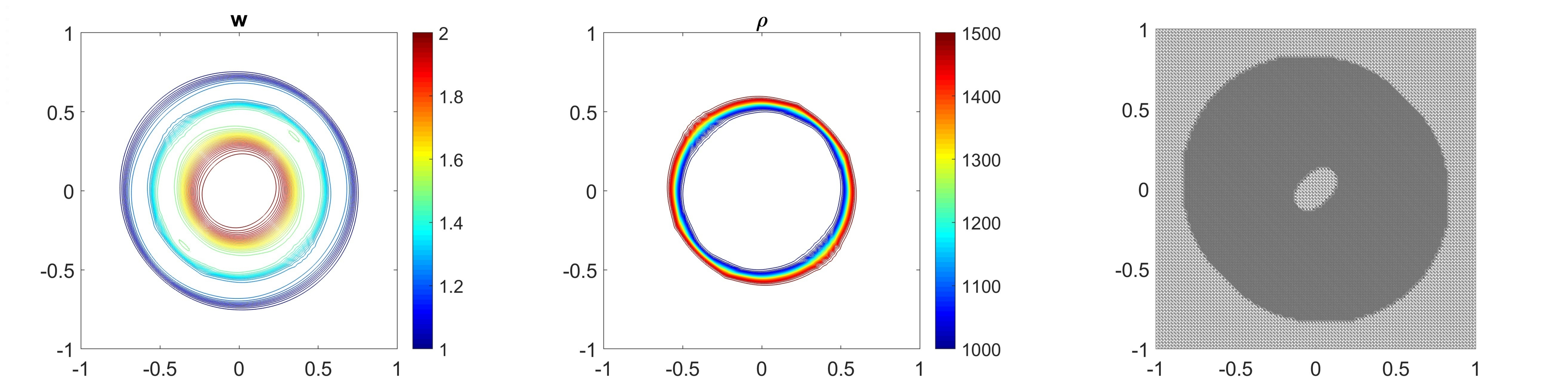}}\label{fig:1anon}\\
		\subfigure[Adaptive scheme with  interface tracking on adaptive mesh with 40940  cells] {\includegraphics[width=1\textwidth]{ex1aconr1l100}}\label{fig:1aint}
	\vspace*{2mm}
			\caption{Example 1: Computational water surface $w(x,y,0.15)$ (first column), density $\rho(x,y,0.15)$ (second column) of the IVP \eref{eq:ex1}, and the corresponding adaptive meshes (third column) obtained by using the proposed adaptive central-upwind scheme (bottom) and using the adaptive scheme  without interface tracking (top).}\label{fig:ex1non}
	\end{figure}

		In addition, we will use this  example to show
          that interface reconstruction presented in Section \ref{sec:interface} plays an important part in preserving the sharpness of the solution as well as improving the accuracy of the adaptive central-upwind scheme. On \fref{fig:ex1non} (a), we plot the water surface $w$ (first) and density $\rho$ (second column) at $t=0.15$ by
          using the adaptive central-upwind method, but implemented
without the interface tracking technique. We then compare the results shown in \fref{fig:ex1non} (a) to  the results calculated by using the adaptive scheme with the interface tracking, \fref{fig:ex1non} (b). The adaptive meshes on \fref{fig:ex1non} (third column) are reconstructed from the uniform mesh $2\times 100\times 100$ with one level of refinement. As can be seen in  \fref{fig:ex1non} (a), both $w$ and $\rho$  are very scattered around the contact wave when we do not track and reconstruct the interface. Meanwhile,  in \fref{fig:ex1non} (b), the proposed adaptive scheme, though using an adaptive mesh with fewer cells (10828 cells fewer), provides more accurate results. 
	
	\par In the next numerical test, we replace the flat bottom with the bottom topography that consists of two Gaussian shaped humps as
	\begin{equation}
	    B(x,y,t)=\begin{cases}
	        0.5e^{-100(x+0.5)^2+(y+0.5)^2}, \quad &\mbox{if} \quad x<0,\\ 
0.6e^{-100(x-0.5)^2+(y-0.5)^2}, \quad &\mbox{if} \quad x\geq 0.
\end{cases}\label{eq:ex1b}
	\end{equation}
			
	\par
\begin{figure}[htp!]
	\centering
	\subfigure[Uniform mesh $2\times 100\times 100$.]{\includegraphics[width=1\textwidth]{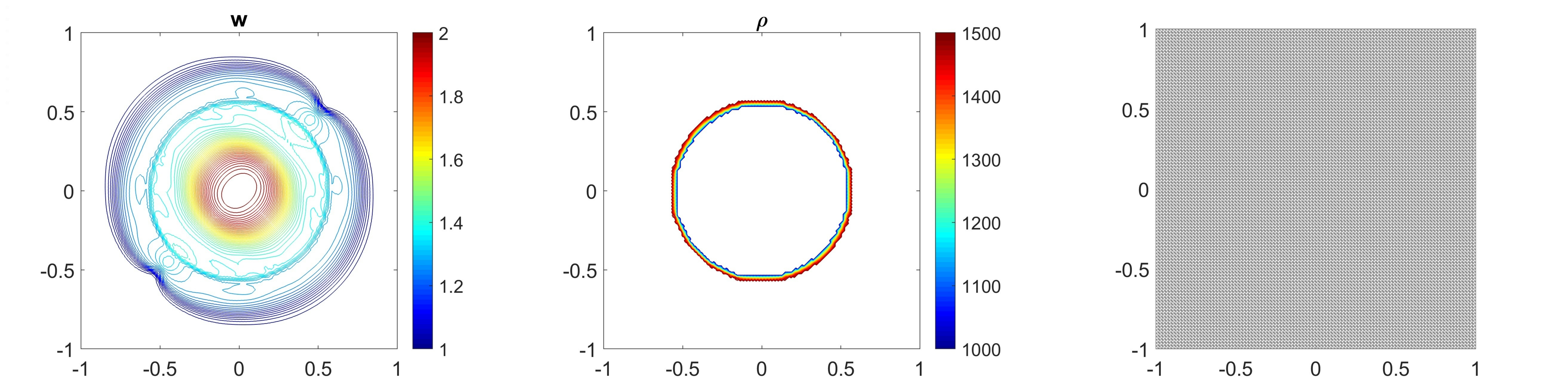}}\label{fig:1ba}\\	
	\subfigure[Uniform mesh $2\times 200\times 200$.]{\includegraphics[width=1\textwidth]{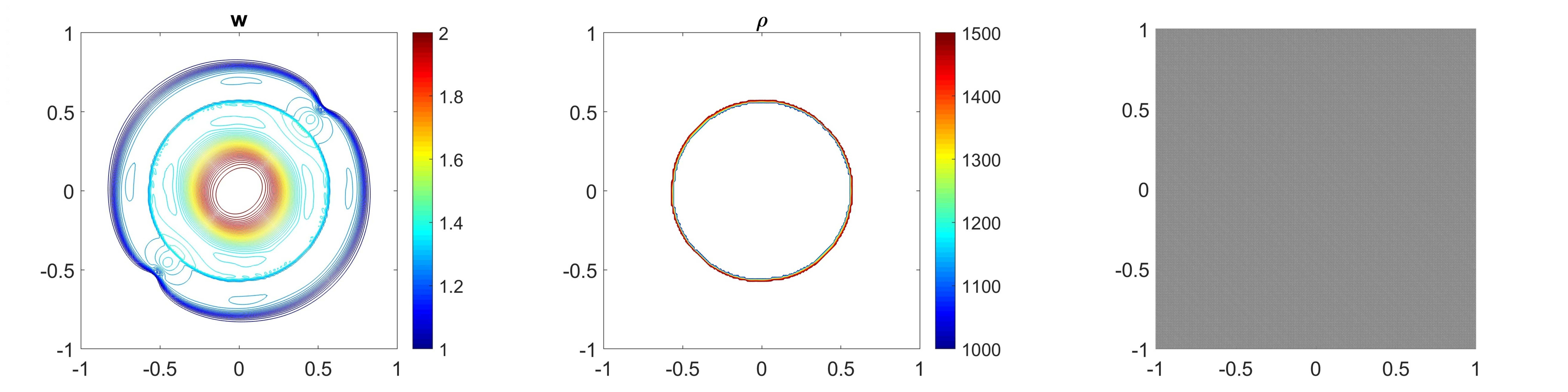}}\label{fig:1bb}\\
	\subfigure[Adaptive mesh with $\mathcal{M}=1$.] {\includegraphics[width=1\textwidth]{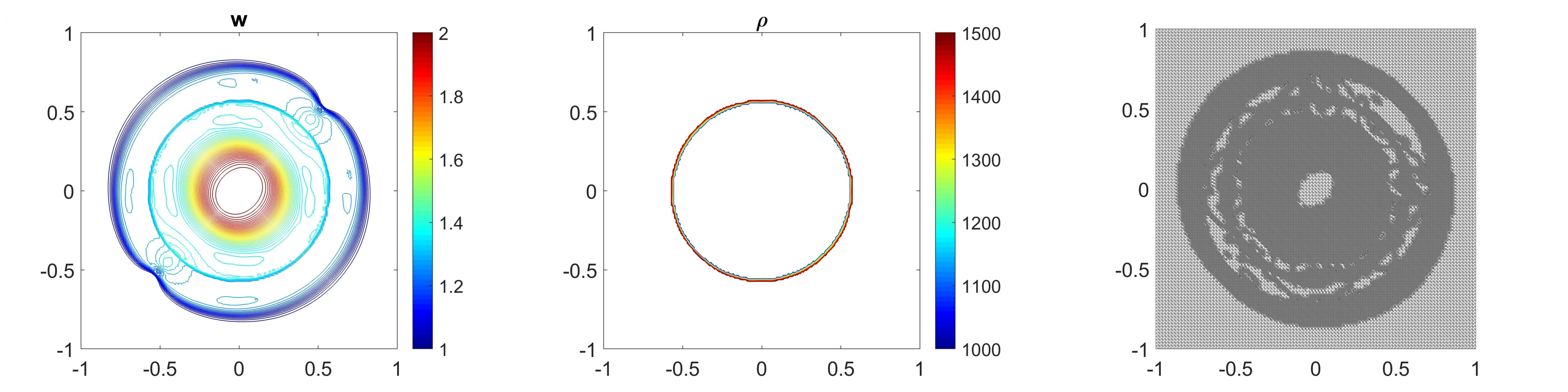}}\label{fig:1bc}\\
	\subfigure[Adaptive mesh with $\mathcal{M}=2$.]{\includegraphics[width=1\textwidth]{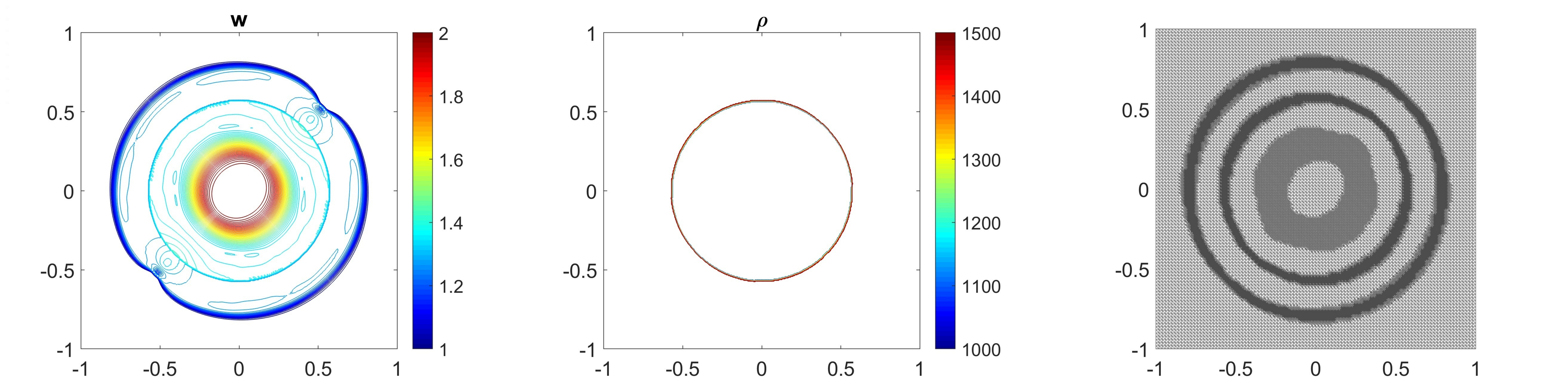}}\label{fig:1bd}\\
	\vspace*{2mm}
	\caption{Example 1: Contour plots of computational water surface $w(x,y,0.2)$ (first column) and density $\rho(x,y,0.2)$ (second column) of the IVP \eref{eq:ex1}-\eref{eq:ex1b} with the corresponding meshes (third column).}\label{fig:ex1b}
\end{figure}

 The purpose of this test is to
	illustrate the performance of the adaptive algorithm in situations having irregular bottom topography. In
	\fref{fig:ex1b}, we show the contour plots of the water surface $w$ (first column) and density $\rho$ (second column) obtained at $t=0.2$ by using the central-upwind scheme with and without adaptivity. The computed solutions of the water surface exhibit reflecting waves where the flow meets the submerged humps. Clearly, from the plots of the adaptive meshes in \fref{fig:ex1b} (third column),  the meshes are adapted to the behavior of
	the flow. Hence, the WLR error estimator is capable to exactly detect the
	location of the steep local gradients in the solution. 

\par	We then recompute the accuracy of the solution for this example \ref{eq:ex1}-\ref{eq:ex1b}, see
        \tref{tab:ex1b}, and the CPU time ratio, see
        \tref{tab:ex1tb}. The results show that the  adaptive scheme uses fewer cells and takes a smaller CPU time to achieve the approximately small $L^1$-errors as computed by the scheme without adaptivity. Therefore, the advantages of the adaptive scheme is maintained in examples with irregular bottom level.
	
		\begin{table}[ht!]
			\vspace*{5mm}
			\centering
			\begin{tabular}{ |c| c c|c|c c|c|c c|}
				\hline
				\multicolumn{3}{|c|}{algorithm without adaptivity}&\multicolumn{6}{c|}{adaptive algorithm}\\
				\cline{4-9}
				\multicolumn{3}{|c|}{ }
				&\multicolumn{3}{c|}{one level $\mathcal{M}=1$}&\multicolumn{3}{c|}{two levels $\mathcal{M}=2$}\\
				\hline
				cells &$L^1$-error  &rate&cells	&$L^1$-error  &rate&cells	&$L^1$-error  &rate\\
				\hline
				$2\times100\times100$&0.0271&&16,252  &0.0273&&16,006 &0.0282&	\\  
				$2\times200\times200$&0.0134	&1.01&57,300  &0.0136&1.00&42,602  & 0.0140&0.99\\
				$2\times400\times400$&0.0050&1.43&213,268 &0.0054&1.40&120,654  &0.0055&1.41\\
				\hline	
			\end{tabular}
			\caption{Example 1: $L^1$-errors of the water surface
				$w$  of the IVP \eref{eq:ex1}-\eref{eq:ex1b} at $t=0.2$, and the convergence rates of the
				central-upwind scheme without adaptivity (uniform mesh $2\times N\times N, N=100,200,400$) and the
				adaptive scheme (the corresponding adaptive
				mesh is reconstructed from the uniform mesh $2\times N/2^\mathcal{M}\times N/2^\mathcal{M}$).}\label{tab:ex1b}
			\vspace*{2mm}
		\end{table}

		\begin{table}[ht!]
			\vspace*{5mm}
			\begin{tabular}{ |c|c| c c|c |c c|}
				\hline
				uniform mesh &\makecell{adaptive mesh\\ $\mathcal{M}=1$}&\multicolumn{2}{c|}{$\mathcal{R}_{CPU}$ with $\mathcal{M}=1$} &\makecell{adaptive mesh\\ $\mathcal{M}=2$}&\multicolumn{2}{c|}{$\mathcal{R}_{CPU}$ with $\mathcal{M}=2$}\\
				\cline{3-4} \cline{6-7}
				(cells)&(cells)&total &\makecell{without\\grid generation} &(cells)&total &\makecell{without\\grid generation}  \\
				\hline
				$2\times100\times100$&16,252 &1.69&1.89&16,006 &1.84&1.94\\
				$2\times200\times200$&57,300  &2.57&2.95&42,602  & 3.41 &3.68 \\
				$2\times400\times400$&213,268 &2.24& 2.60 &120,654&2.63& 2.84 \\
				\hline
				\multicolumn{2}{|c}{$\mathcal{R}_{CPU}$ average:}&2.17&\multicolumn{1}{c}{2.48 }&\multicolumn{1}{c}{} &2.63& 2.82 \\
				\hline	
			\end{tabular}
			\caption{Example 1: $\mathcal{R}_{CPU}$ ratio for the IVP \eref{eq:ex1}-\eref{eq:ex1b} at	$t=0.2$, where for adaptive central-upwind scheme,
				we consider the total CPU times and CPU times
				without the grid generation.}\label{tab:ex1tb}
		\end{table}

	\subsection{Example 2:}

	The second numerical example here was proposed in \citep{chertock2014central} to verify the capability of the adaptive scheme in preserving the steady state solution in ``lake at rest' problems, \eref{eq:lara} and \eref{eq:larb}. In particular, the initial data consists of two “lake at rest” states of type \eref{eq:lara} connected
through the density jump corresponding to the “lake at rest” state of type \eref{eq:larb} as 

\begin{equation}
    	(w,u,v,\rho)^T(x,y,0)=\begin{cases}
	    (3,0,0,\dfrac{4}{3}\rho_0), \quad \mbox{if}\quad  & x^2+y^2<0.25,\\
	    (2,0,0,3\rho_0), \quad & \mbox{otherwise}.
	\end{cases}\label{eq:ex2}
\end{equation}

\par
\begin{figure}[ht!]
	\centering
	\subfigure[Uniform mesh $2\times 100\times 100$.]{\includegraphics[width=1\textwidth]{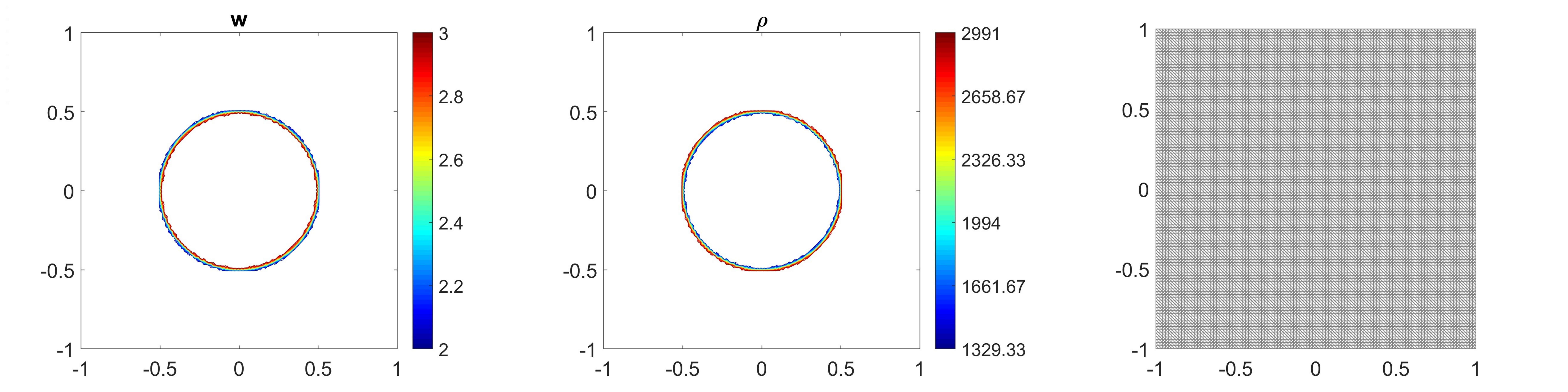}}\label{fig:2a}\\	
	\subfigure[Uniform mesh $2\times 200\times 200$.]{\includegraphics[width=1\textwidth]{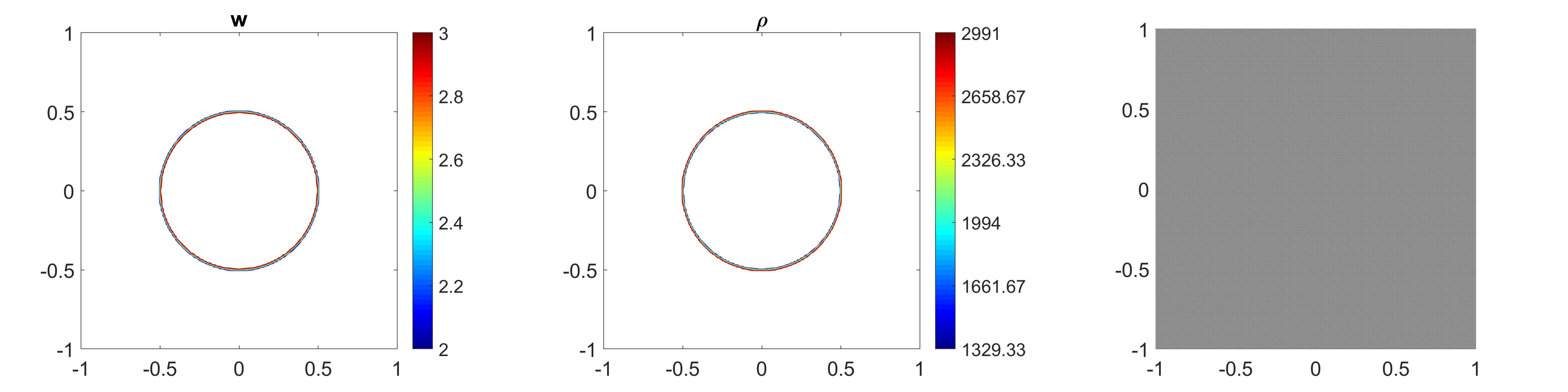}}\label{fig:2b}\\
	\subfigure[Adaptive mesh with $\mathcal{M}=1$.] {\includegraphics[width=1\textwidth]{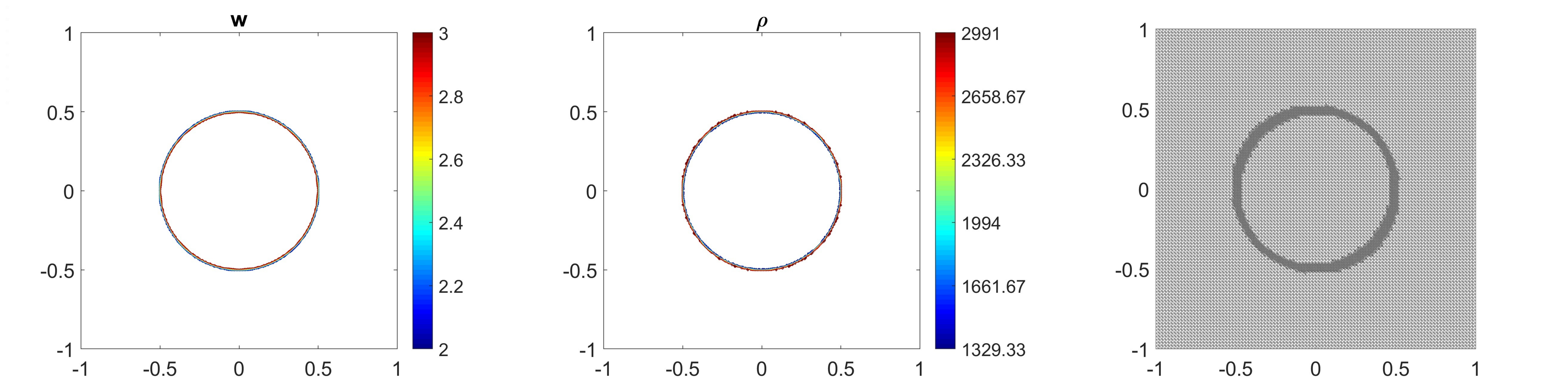}}\label{fig:2c}\\
	\subfigure[Adaptive mesh with $\mathcal{M}=2$.]{\includegraphics[width=1\textwidth]{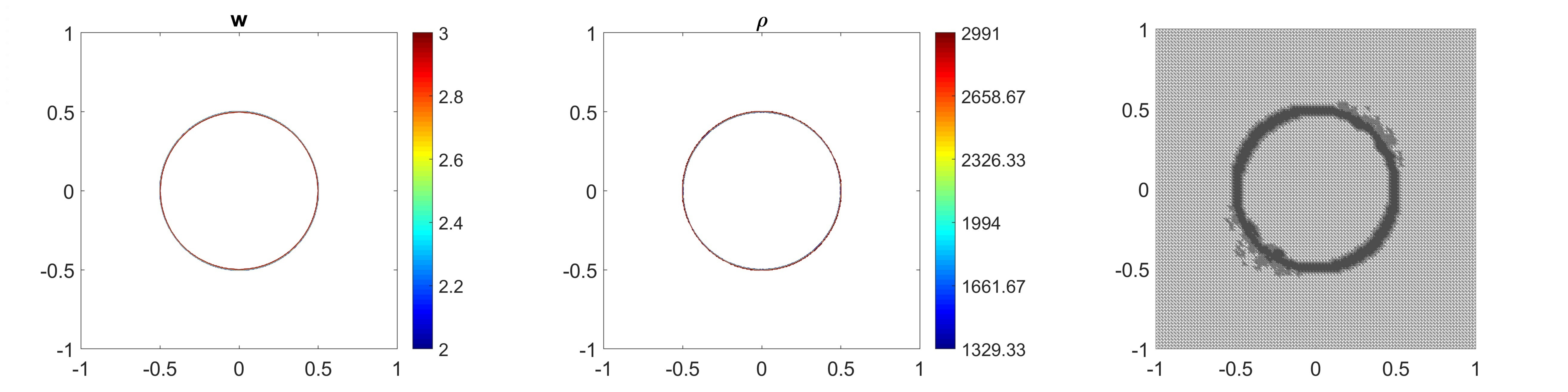}}\label{fig:2d}\\
	\vspace*{2mm}
	\caption{Example 2: Contour plots of computational water surface $w(x,y,0.15)$ (first column) and density $\rho(x,y,0.15)$ with the corresponding meshes (right column).}\label{fig:ex2}
\end{figure}

\begin{figure}[htp!]
	\centering
	\subfigure[Uniform mesh $2\times 100\times 100$.]{\includegraphics[width=0.7\textwidth]{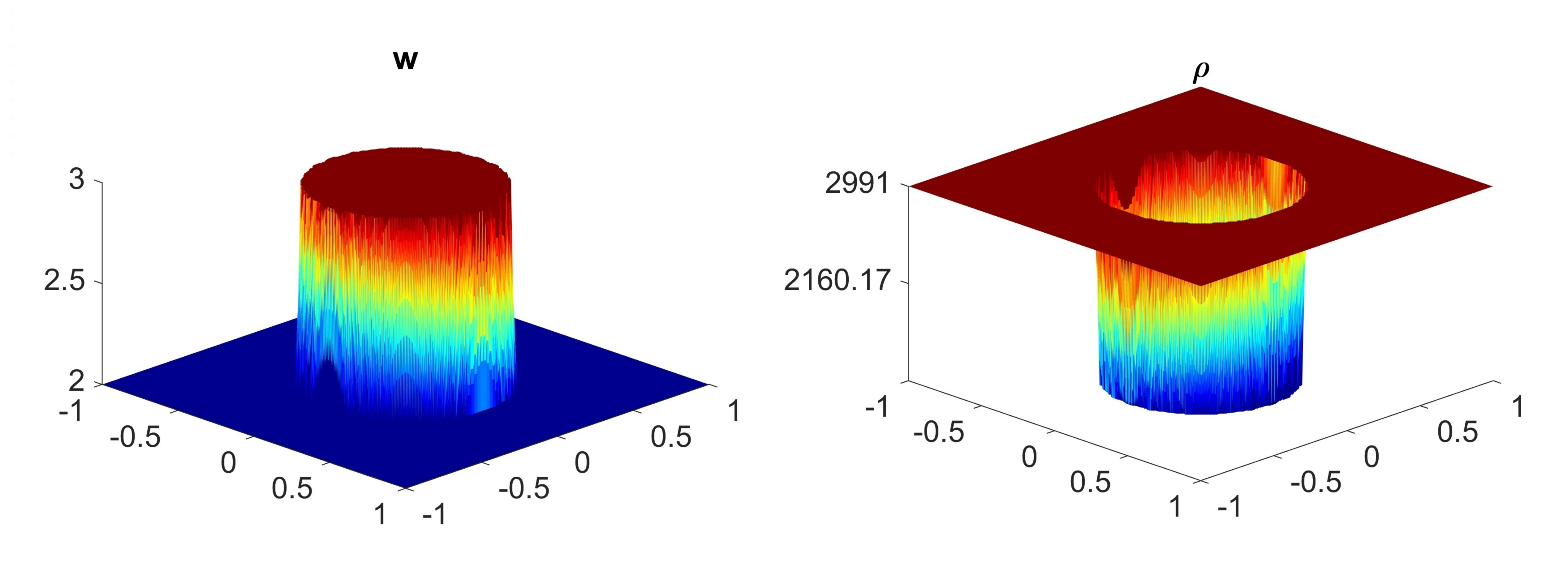}}\label{fig:2as}\\	
	\subfigure[Uniform mesh $2\times 200\times 200$.]{\includegraphics[width=0.7\textwidth]{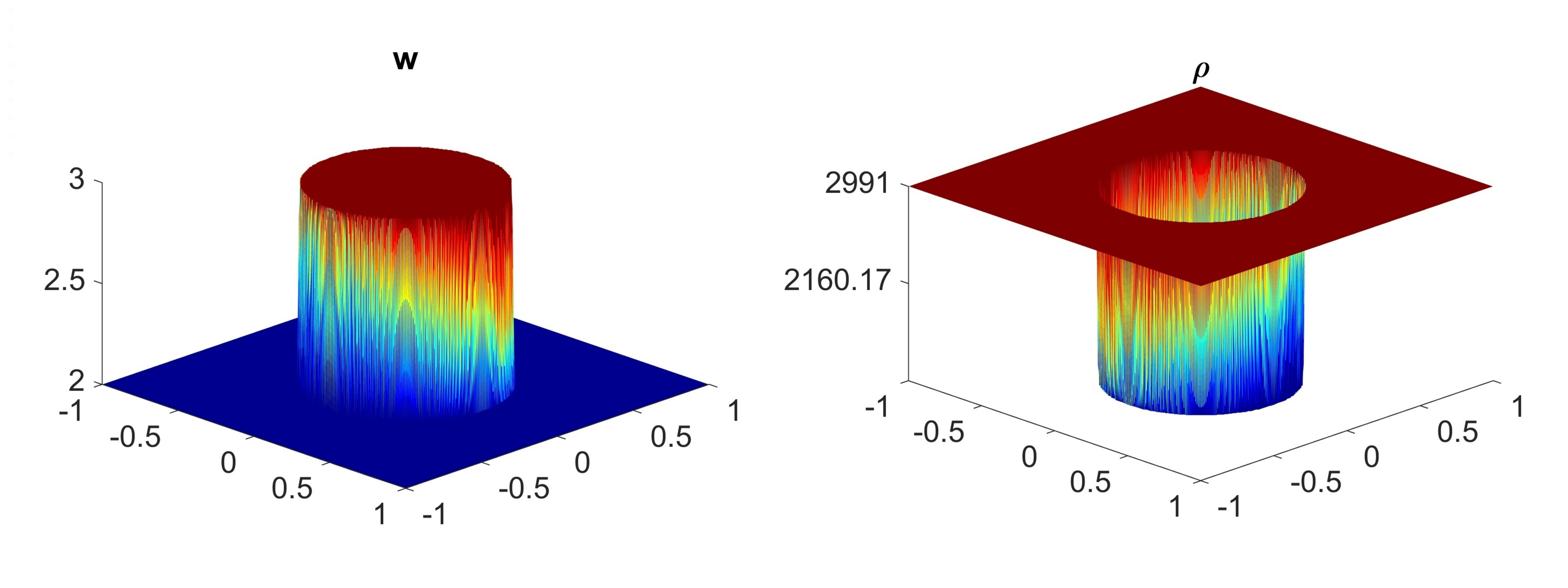}}\label{fig:2bs}\\
	\subfigure[Adaptive mesh with $\mathcal{M}=1$.] {\includegraphics[width=0.7\textwidth]{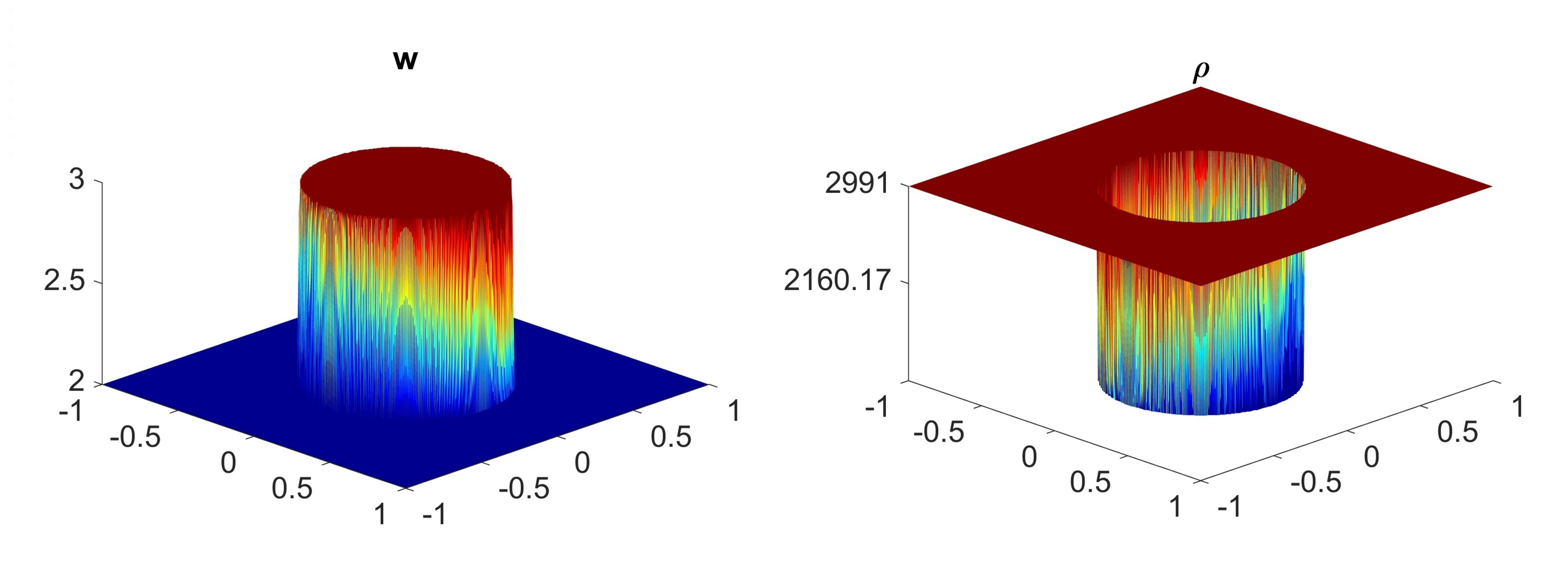}}\label{fig:2cs}\\
	\subfigure[Adaptive mesh with $\mathcal{M}=2$.]{\includegraphics[width=0.7\textwidth]{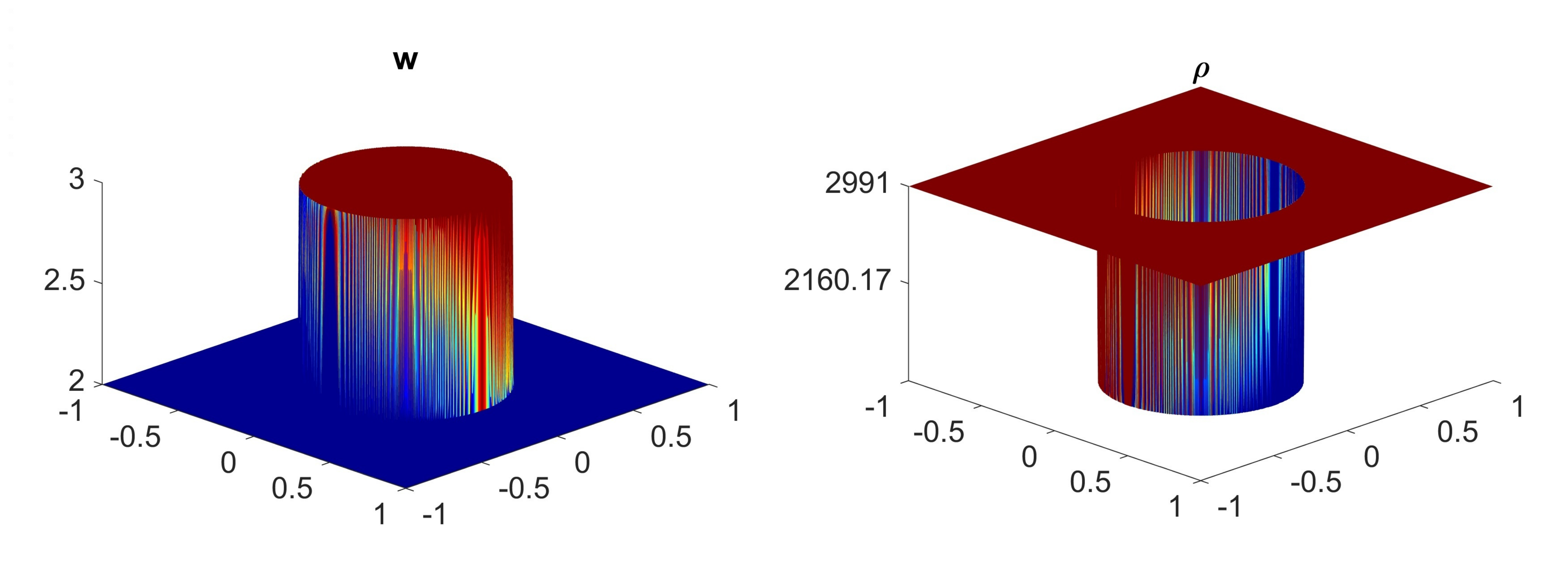}}\label{fig:2ds}\\
	\vspace*{2mm}
	\caption{Example 2: 3-D plots of computational water surface $w(x,y,0.15)$ (left column) and density $\rho(x,y,0.15)$ (right column).}\label{fig:ex2s}
\end{figure}

  In this example, we consider  the bottom topography \eref{eq:ex1b} on a computational domain $[-1, 1]\times[-1, 1]$.  To reconstruct the adaptive meshes, the threshold is set, $\omega = 0.1\max_j(e_j)$. 
In \fref{fig:ex2}, we present the plots of the computed water surface (first column) and density (second column) at $t=0.15$ obtained by using the central-upwind scheme, but without adaptivity, \fref{fig:ex2} (a, b)  and by using the adaptive algorithm  \fref{fig:ex2} (c, d).  The adaptive grids plotted on \fref{fig:ex2} (third column) are generated from the uniform mesh $2\times 100\times 100$ with one
level of refinement $\mathcal{M}= 1$, \fref{fig:ex2} (c), and two levels of refinement $\mathcal{M}= 2$, \fref{fig:ex2} (d). \fref{fig:ex2s} shows the 3D plots of the numerical solution computed by the two methods. As expected, in \fref{fig:ex2} and \fref{fig:ex2s}, the adaptive scheme with interface tracking exactly preserves the steady state. Hence, in \fref{fig:ex2} (third column), the WLR error only marks cells surrounding the circle of density jump for refinement.  In addition,  no pressure oscillations are observed at the interface. 

\par Next, we will illustrate the advantages of the adaptive central-upwind scheme. We compute the $L^1$-error, Table \ref{tab:ex2}, and the CPU ratio,  Table \ref{tab:ex2t}, by using the central-upwind method without adaptivity and using the adaptive algorithm presented in our work. In order to
compare the computational costs and calculate $\mathcal{R}_{CPU}$, we consider uniform and adaptive meshes
with the same size of the smallest cells. Table \ref{tab:ex2} shows that in the adaptive meshes, we achieve  $L^1$-errors as small as the errors obtained in the corresponding uniform meshes. However, the adaptive algorithm uses fewer cells and reduces the CPU times up to eight times. The computational cost is remarkably cut down since as illustrated in \fref{fig:ex2} and \fref{fig:ex2s}, only a few cells in the neighborhood of the density discontinuity have large WLR errors and are therefore marked for refinement. This example has clearly show the efficiency of the proposed scheme for numerically solving the system of multi-fluid flow.

		\begin{table}[ht]
	\vspace*{5mm}
	\centering
	\begin{tabular}{ |c| c c|c|c c|c|c c|}
		\hline
		\multicolumn{3}{|c|}{algorithm without adaptivity}&\multicolumn{6}{c|}{adaptive algorithm}\\
		\cline{4-9}
		\multicolumn{3}{|c|}{ }
		&\multicolumn{3}{c|}{one level $\mathcal{M}=1$}&\multicolumn{3}{c|}{two levels $\mathcal{M}=2$}\\
		\hline
		cells &$L^1$-error  &rate&cells	&$L^1$-error  &rate&cells	&$L^1$-error  &rate\\
		\hline
		$2\times100\times100$&0.0045&&6,284  &0.0045&&5,928 &0.0045&	\\  
		$2\times200\times200$&0.0021&1.10&22,546 &0.0021&1.10&11,260 & 0.0021&1.10\\
		$2\times400\times400$&7.9252e-04&1.41&85,132   &8.5079e-04&1.30&33,904& 8.9818e-04&1.23\\
		\hline	
	\end{tabular}
	\caption{Example 2: $L^1$-errors of the water surface
		$w$ at $t=0.15$, and the convergence rates of the
		central-upwind scheme without adaptivity (uniform mesh $2\times N\times N, N=100,200,400$) and the
		adaptive scheme (the corresponding adaptive
		mesh is reconstructed from the uniform mesh $2\times N/2^\mathcal{M}\times N/2^\mathcal{M}$).}\label{tab:ex2}
	\vspace*{2mm}
\end{table}

\begin{table}[ht]
	\vspace*{5mm}
	\begin{tabular}{ |c|c| c c|c |c c|}
		\hline
		uniform mesh &\makecell{adaptive mesh\\ $\mathcal{M}=1$}&\multicolumn{2}{c|}{$\mathcal{R}_{CPU}$ with $\mathcal{M}=1$} &\makecell{adaptive mesh\\ $\mathcal{M}=2$}&\multicolumn{2}{c|}{$\mathcal{R}_{CPU}$ with $\mathcal{M}=2$}\\
		\cline{3-4} \cline{6-7}
		(cells)&(cells)&total &\makecell{without\\grid generation} &(cells)&total &\makecell{without\\grid generation}  \\
		\hline
		$2\times100\times100$&6,284 &3.01&3.54&5,928&3.48&3.75\\
		$2\times200\times200$&22,546  &3.78&4.76&11,260  & 7.76 &8.66 \\
		$2\times400\times400$&85,132 &4.46& 5.57 &33,904  &11.13&12.86 \\
		\hline
		\multicolumn{2}{|c}{$\mathcal{R}_{CPU}$ average:}&3.75&\multicolumn{1}{c}{4.62 }&\multicolumn{1}{c}{} &7.46& 8.42 \\
		\hline	
	\end{tabular}
	\caption{Example 2: $\mathcal{R}_{CPU}$ ratio at
		$t=0.15$, where for adaptive central-upwind scheme,
		we consider the total CPU times and CPU times
		without the grid generation.}\label{tab:ex2t}
\end{table}

	\subsection{Example 3:}
	
\begin{figure}[H]
		\centering
		\includegraphics[width=0.4\textwidth]{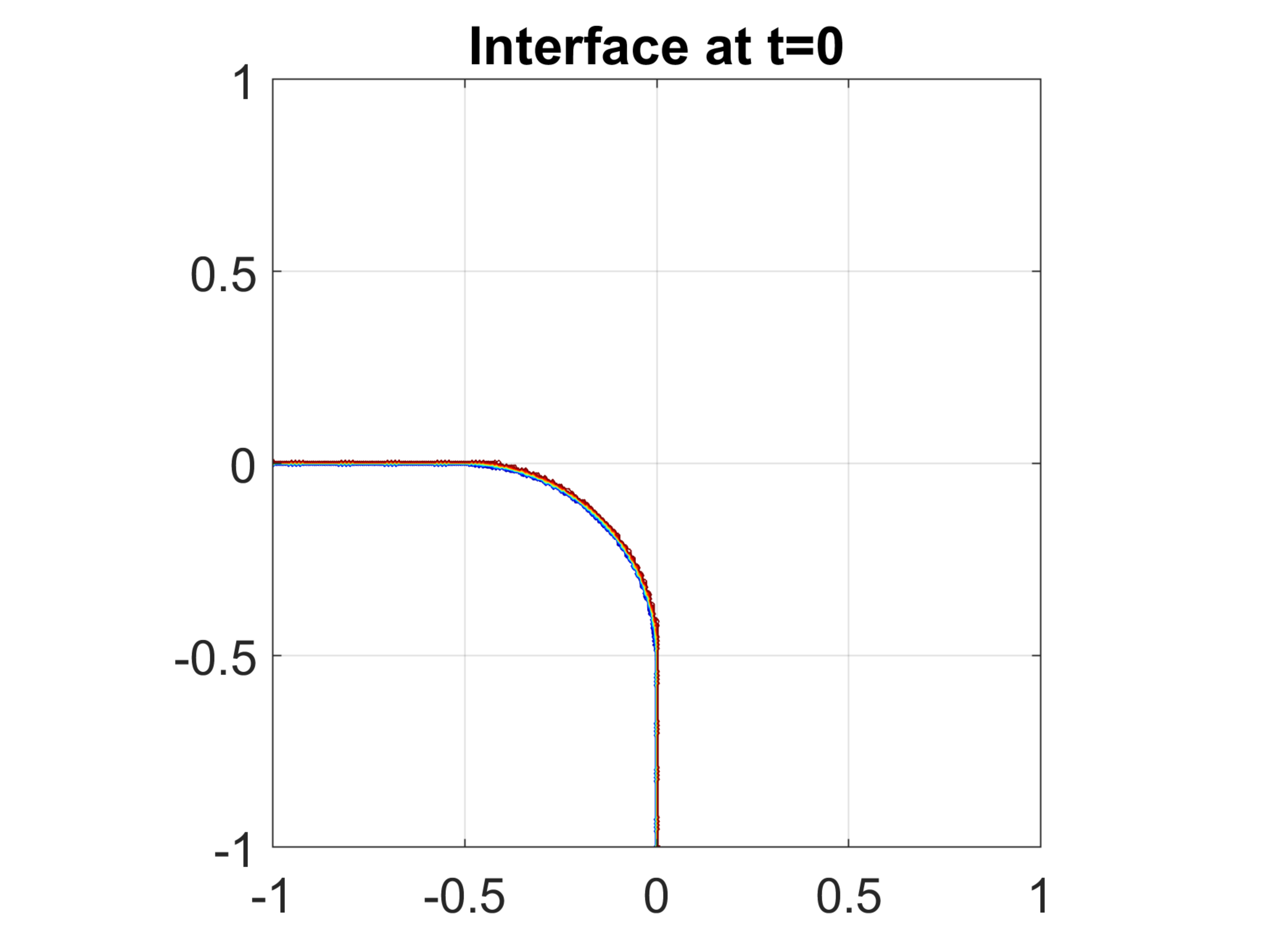}
		\vspace*{2mm}
		\caption{Example 3: The interface at initial time $t=0$.}\label{fig:ex3int}
	\end{figure}

	The last example is designed to illustrate the capability of the proposed adaptive algorithm to handle irregular density interfaces. Hence, in a domain $[-1,1]\times[-1,1]$, the density jump at $t=0$ is given by a curve  which consists of a horizontal segment, a vertical segment, and a quarter of a circle connected at their endpoints as illustrated in \fref{fig:ex3int}. 
		
\begin{figure}[htp!]
		\centering
		\subfigure[Uniform mesh $2\times 100\times 100$.]{\includegraphics[width=1\textwidth]{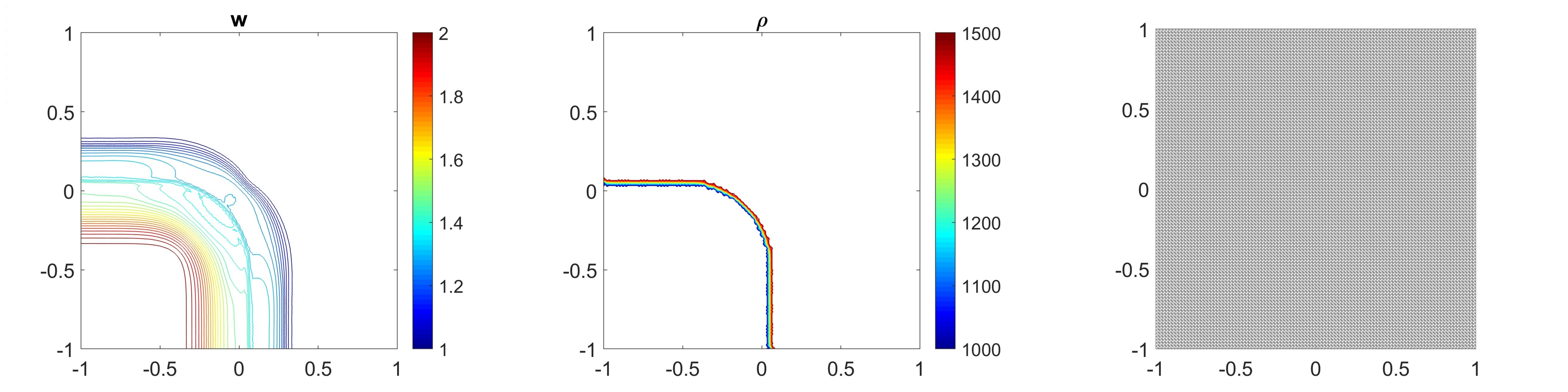}}\label{fig:3a}\\	
		\subfigure[Uniform mesh $2\times 200\times 200$.]{\includegraphics[width=1\textwidth]{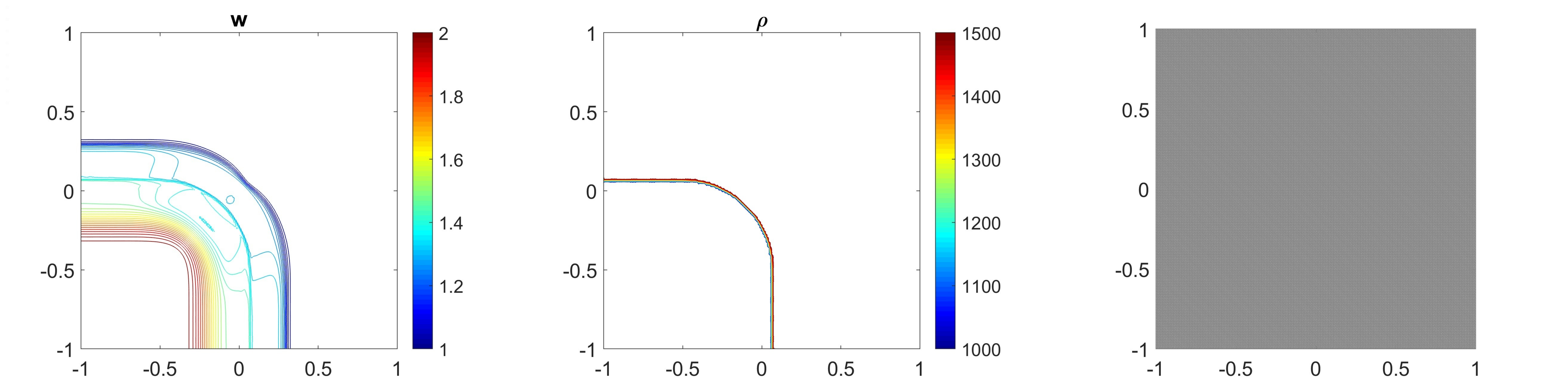}}\label{fig:3b}\\
		\subfigure[Adaptive mesh with $\mathcal{M}=1$.] {\includegraphics[width=1\textwidth]{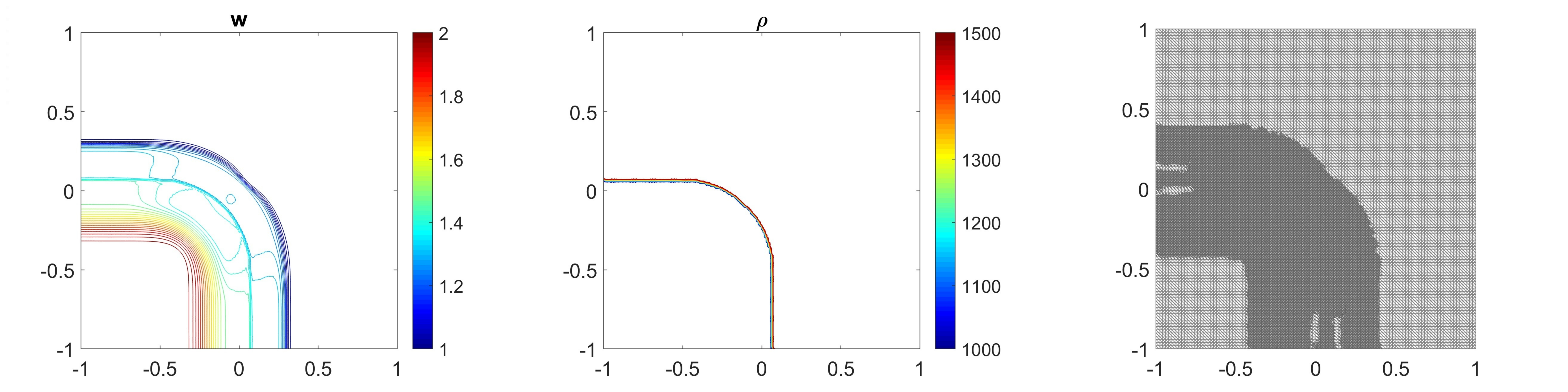}}\label{fig:3c}\\
		\subfigure[Adaptive mesh with $\mathcal{M}=2$.]{\includegraphics[width=1\textwidth]{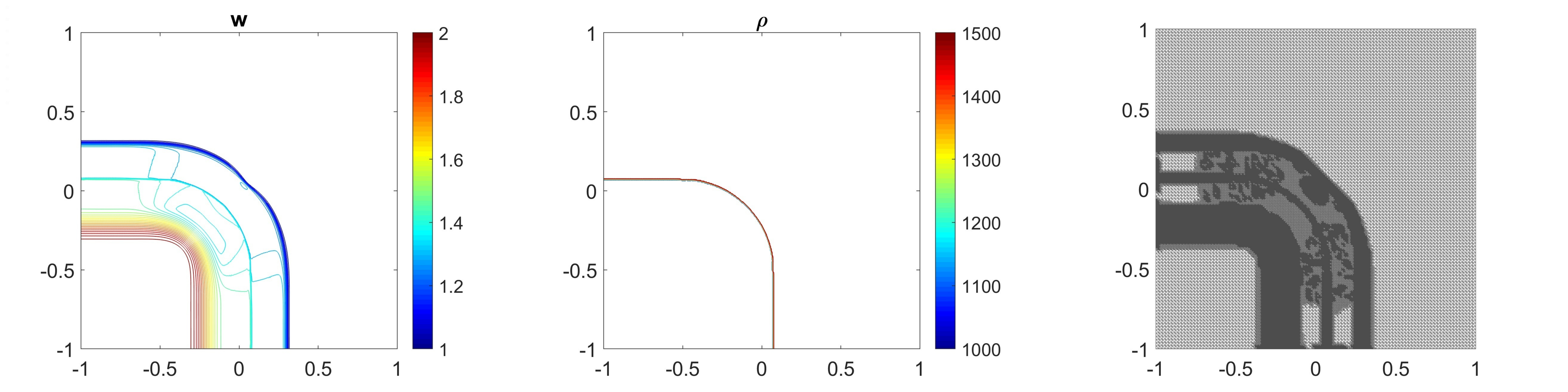}}\label{fig:3d}\\
		\vspace*{2mm}
		\caption{Example 3: Contour plots of computational water surface $w(x,y,0.15)$ (first column) and density $\rho(x,y,0.15)$ (second column) with the corresponding meshes (third column).}\label{fig:ex3}
	\end{figure}

\begin{figure}[htp!]
	\centering
	\subfigure[Uniform mesh $2\times 100\times 100$.]{\includegraphics[width=0.7\textwidth]{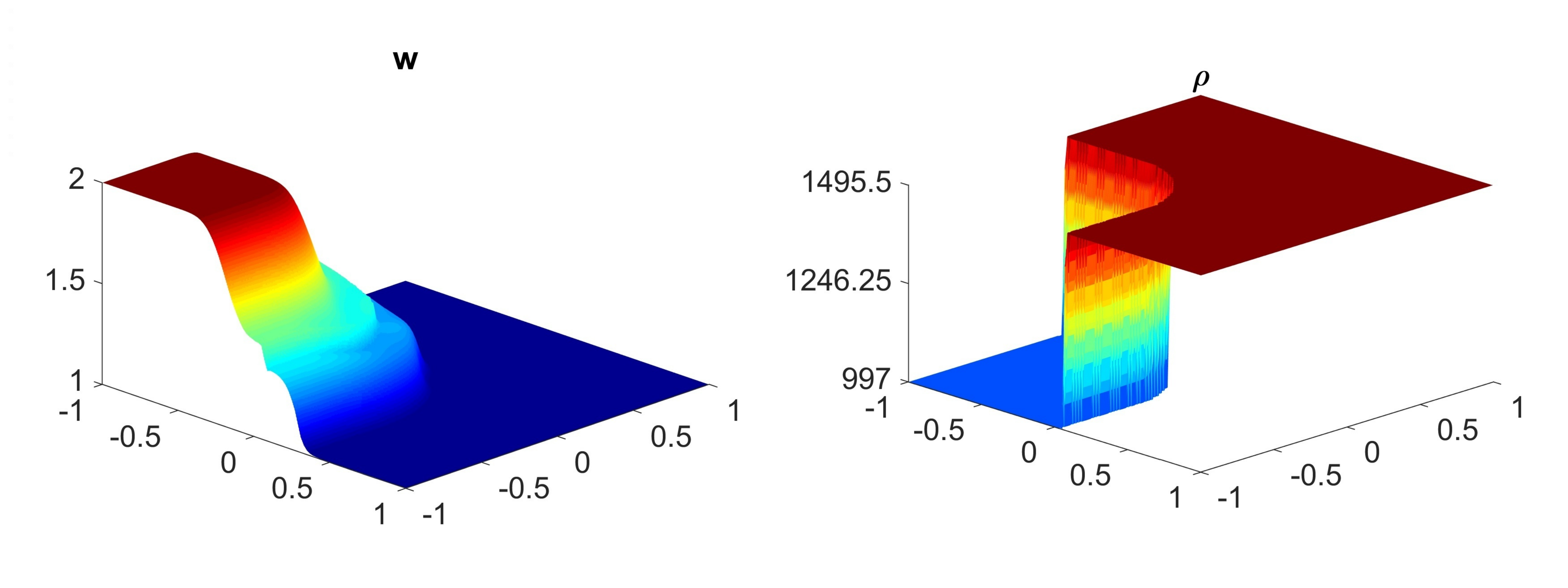}}\label{fig:3as}\\	
	\subfigure[Uniform mesh $2\times 200\times 200$.]{\includegraphics[width=0.7\textwidth]{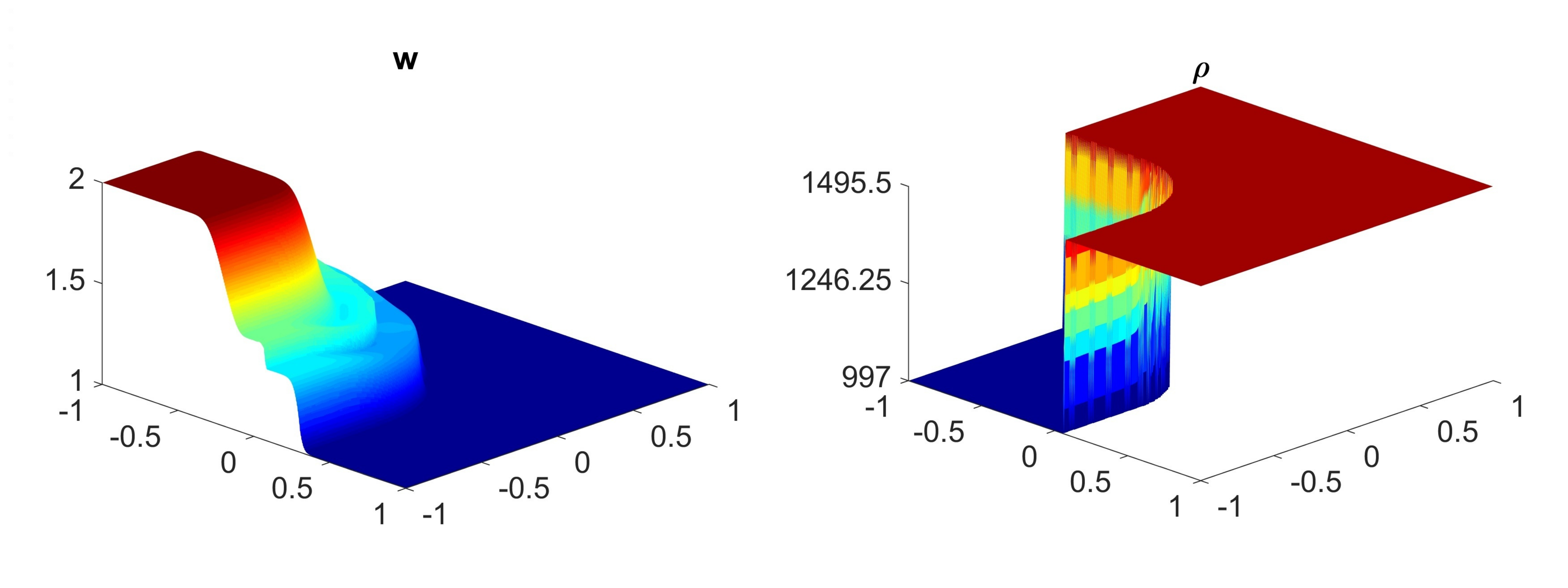}}\label{fig:3bs}\\
	\subfigure[Adaptive mesh with $\mathcal{M}=1$.] {\includegraphics[width=0.7\textwidth]{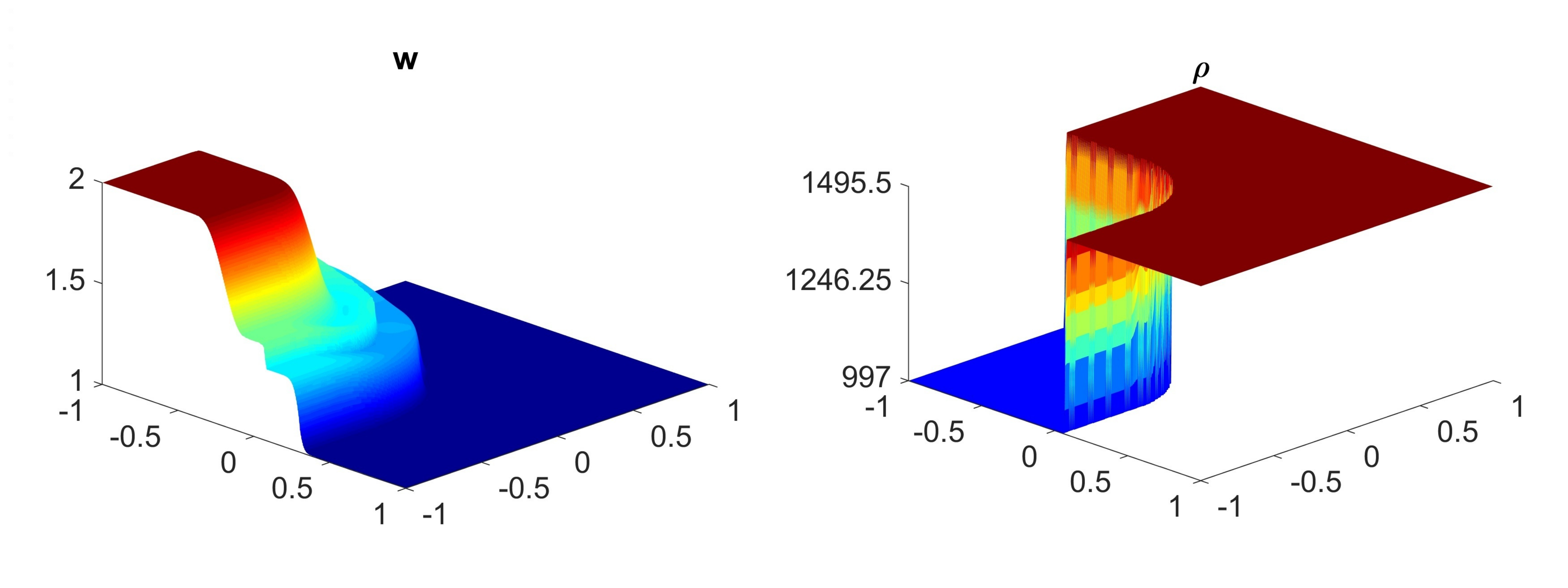}}\label{fig:3cs}\\
	\subfigure[Adaptive mesh with $\mathcal{M}=2$.]{\includegraphics[width=0.7\textwidth]{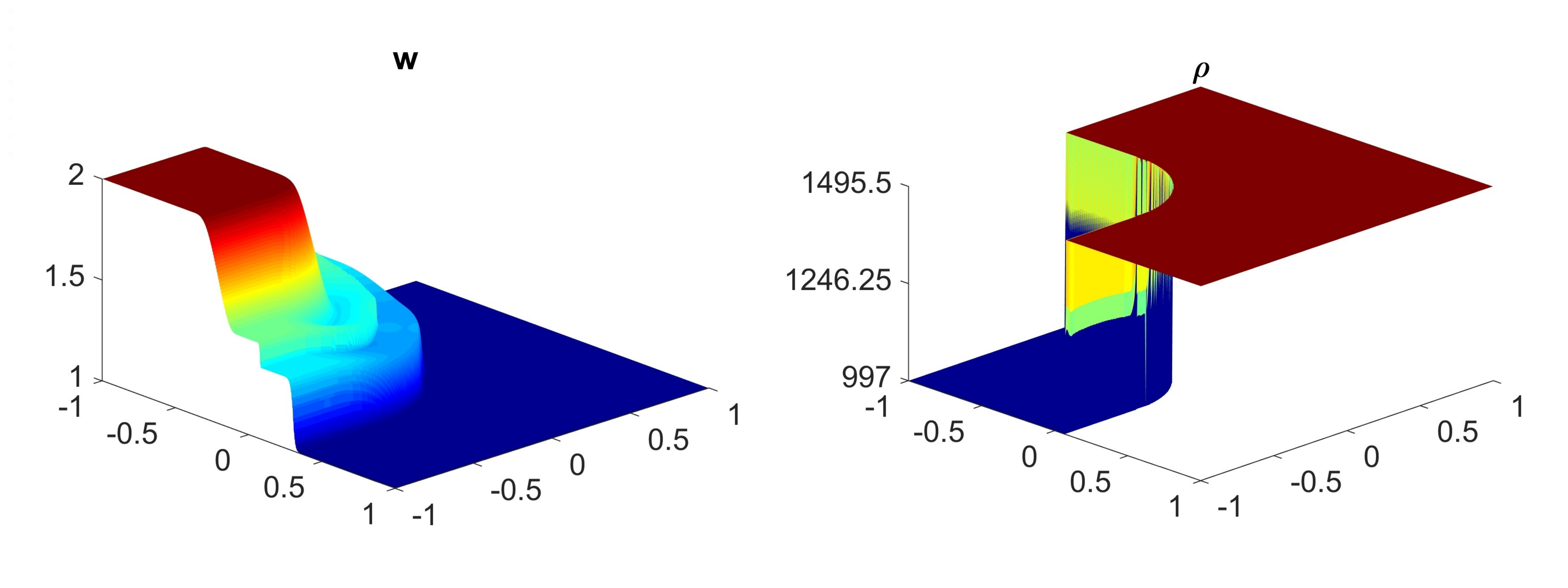}}\label{fig:3ds}\\
	\vspace*{2mm}
	\caption{Example 3: 3-D plots of computational water surface $w(x,y,0.15)$ (left column) and density $\rho(x,y,0.15)$ (right column).}\label{fig:ex3s}
\end{figure}

The initial condition is described by
	
	\begin{equation}
    	(w,u,v,\rho)^T(x,y,0)=\begin{cases}
	    (2,0,0,\rho_0), \quad \mbox{if}\quad  & (x,y)\in \Omega,\\
	    (1,0,0,1.5\rho_0), \quad & \mbox{otherwise}.
	\end{cases}\label{eq:ex2}
\end{equation}
where 	$$\Omega:=\{x<-0.5,y<0\}\cup \{(x+0.5)^2+(y+0.5)^2<0.25\}\cup\{x<0, y<-0.5\}.$$. 
We consider a bottom topography with a surbmerged hump as
\begin{equation*}
B(x,y,t)=0.5e^{-100(x^2+y^2)}.    
\end{equation*}
In this example, we will also perform the same numerical tests which are done in previous examples. Namely, we first calculate the water surface and density at $t=0.15$ using central-upwind scheme without adaptivity and present the results in \fref{fig:ex3} (a, b) and \fref{fig:ex3s} (a, b). In \ref{fig:ex3} (c, d) and \ref{fig:ex3s} (c, d), we
plot the results for $w$ (first column) and $\rho$ (second column) obtained by the adaptive scheme. The adaptive grids in  \ref{fig:ex3} (third column)
are generated from the uniform grid $2\times 100\times 100$ for one level of refinement $\mathcal{M}= 1$ and $\mathcal{M}= 2$ using the threshold $\omega= 0.01 \max_j(e_j)$. As expected, the solutions obtained in the adaptive meshes with high levels of refinement are much sharper than the solutions computed in fixed uniform meshes. The density jump moves Northeast and does not diffuse. There is no non-physical spurious waves generated at the interface. Also, as can be seen in the adaptive meshes, the WLR error indicator captures subtle features of the solution.

Finally, we compute the $L^1$-errors in Table \ref{tab:ex3} and the CPU ratio in Table \ref{tab:ex3t} by using the  adaptive scheme and using the central-upwind method without the mesh generation. By comparing the results presented in Tables  \ref{tab:ex3} and  \ref{tab:ex3t}, one can easily see that at a reduced computational cost, the proposed adaptive central-upwind method is still able to obtain the accurate solutions for this example. In our experiments, we considered only $\mathcal{M} = 1$ and
$\mathcal{M} = 2$, but  to further enhance the accuracy of the numerical solution at the lower computational cost, one can consider higher levels of refinement.

		\begin{table}[ht!]
	\vspace*{5mm}
	\centering
	\begin{tabular}{ |c| c c|c|c c|c|c c|}
		\hline
		\multicolumn{3}{|c|}{algorithm without adaptivity}&\multicolumn{6}{c|}{adaptive algorithm}\\
		\cline{4-9}
		\multicolumn{3}{|c|}{ }
		&\multicolumn{3}{c|}{one level $\mathcal{M}=1$}&\multicolumn{3}{c|}{two levels $\mathcal{M}=2$}\\
		\hline
		cells &$L^1$-error  &rate&cells	&$L^1$-error  &rate&cells	&$L^1$-error  &rate\\
		\hline
		$2\times100\times100$&0.0155&&12,164  &0.0145&&8,994 &0.0166&	\\  
		$2\times200\times200$&0.0074&1.07&40,814  &0.0074&0.97&30,263 & 0.0071&1.23\\
		$2\times400\times400$&0.0027&1.45&148,473  &0.0028&1.40&87,214 &0.0028&1.34\\
		\hline	
	\end{tabular}
	\caption{Example 3: $L^1$-errors of the water surface
		$w$ at $t=0.15$, and the convergence rates of the
		central-upwind scheme without adaptivity (uniform mesh $2\times N\times N, N=100,200,400$) and the
		adaptive scheme (the corresponding adaptive
		mesh is reconstructed from the uniform mesh $2\times N/2^\mathcal{M}\times N/2^\mathcal{M}$).}\label{tab:ex3}
	\vspace*{2mm}
\end{table}

\begin{table}[ht!]
	\vspace*{5mm}
	\begin{tabular}{ |c|c| c c|c |c c|}
		\hline
		uniform mesh &\makecell{adaptive mesh\\ $\mathcal{M}=1$}&\multicolumn{2}{c|}{$\mathcal{R}_{CPU}$ with $\mathcal{M}=1$} &\makecell{adaptive mesh\\ $\mathcal{M}=2$}&\multicolumn{2}{c|}{$\mathcal{R}_{CPU}$ with $\mathcal{M}=2$}\\
		\cline{3-4} \cline{6-7}
		(cells)&(cells)&total &\makecell{without\\grid generation} &(cells)&total &\makecell{without\\grid generation}  \\
		\hline
		$2\times100\times100$&12,164  &1.75&1.98&8,994 &2.69&2.87\\
		$2\times200\times200$&40,814 &2.71&3.14&30,263  & 3.57 &3.86 \\
		$2\times400\times400$&148,473  &2.86& 3.38 &87,214&4.60&5.03 \\
		\hline
		\multicolumn{2}{|c}{$\mathcal{R}_{CPU}$ average:}&2.44&\multicolumn{1}{c}{2.83 }&\multicolumn{1}{c}{} &3.62&3.92 \\
		\hline	
	\end{tabular}
	\caption{Example 3: $\mathcal{R}_{CPU}$ ratio at
		$t=0.15$, where for adaptive central-upwind scheme,
		we consider the total CPU times and CPU times
		without the grid generation.}\label{tab:ex3t}
\end{table}

	\section{Conclusion}\label{8sect05}
We have developed a new adaptive well-balanced and positivity preserving central-upwind scheme on unstructured traingular meshes for
shallow water equations with variable density. The scheme is designed as an extension the scheme in \citep{epshteyn2020adaptive} by utilizing the interface tracking method in \citep{chertock2014central} and the interface reconstruction in \citep{GHAZIZADEH2020104633}. The proposed scheme is capable to preserve the steady state solutions \eref{eq:lara} and \eref{eq:larb} and prevent the oscillation at the density jumps. In addition, to achieve an efficient strategy for the
adaptive mesh reconstruction, we also obtain a robust local error indicator. We performed several challenging numerical
tests for multi-fluid models and we demonstrated that the new
adaptive central-upwind scheme maintains well-balanced and positivity-preserving properties and obtains high-accuracy at a reduced computational cost.
			\section*{Acknowledgements}
			I would like to sincerely thank my advisor, Prof. Yekaterina Epshteyn, for her suggestion of the problem to investigate, her encouragement and help. The proposed approach in this paper is an extension of the adaptive method developed in \citep{epshteyn2020adaptive} which is my joint work with her. Without her support, this work would not have started, progressed, or ended. 
	\section*{References} 
	
	\bibliography{ref}
	\bibliographystyle{siam}
\end{document}